




\input amssym.def
\input amssym.tex
 \font\newrm =cmr10 at 24pt
\def\bul{\raise .9pt\hbox{\newrm .\kern-.105em } }

 \def\fr{\frak}

 \def\h{\hbox{ }}
 
 \def\p{{\fr p}}
 \def\u{{\fr u}}
 \def\r{{\fr r}}
 
 \def\m{{\fr m}}
 \def\n{{\fr n}}
 \def\a{{\fr a}}
 \def\d{{\fr d}}
 
 \def\ss{{\fr s}}
 \def\k{{\fr k}}
 \def\b{{\fr b}}
 \def\cc{{\fr c}}
 \def\hh{{\fr h}}
 \def\tt{{\fr t}}
 \def\ee{{\fr e}}

 \def\g{{\fr g}}
 \def\v{{\fr v}}
 
 \def\q{{\fr q}}

 \def\<{\le}
 \def\>{\ge}
 
 \def\i{^{-1}}

 \def\s{{\h\subset\h}}

 \def\mapright#1
  {\smash{\mathop
  {\longrightarrow}
  \limits^{#1}}}

 \def\kk#1{{\kern .4 em} #1}

\hoffset=2pc
\voffset=2pc
\hsize = 31pc
\vsize = 45pc
\overfullrule = 0pt

\centerline{\bf A Branching Law for Subgroups Fixed by an
Involution}

\centerline{\bf and 
  a Concompact Analogue of the Borel-Weil Theorem}

\vskip 1pc 

\centerline{\bf Bertram Kostant}

\vskip 1pc 

{\bf Abstract.} We give a branching law for subgroups fixed by
an involution. As an application we give a generalization of
the Cartan-Helgason theorem and a noncompact analogue of the
Borel-Weil theorem.

\centerline{\bf 0. Introduction}\vskip 1pc
{\bf 0.1.} Let $G_{\Bbb C}$ be a simply-connected complex semisimple Lie
group and let $\g= Lie\,G_{\Bbb C}$. Let $\g_o$ be a real form of
$\g$ and let $G$ be the real semisimple Lie subgroup of $G_{\Bbb
C}$ corresponding to $\g_o$. Using standard notation let $G =
K\,A\,N$ be an Iwasawa decomposition of $G$ and let $M$ be the
centralizer of $A$ in the maximal compact subgroup $K$ of $G$.
The complexified Lie algebras of $K,M,A,$ and $N$ are denoted
respectively by $\k,\m,\a$ and $\n$. Let $\hh_{\m}$ be a
Cartan subalgebra of $\m$ so that $\hh = \hh_{\m} + \a$ is a
Cartan subalgebra of $\g$. Let $\m = \m_- + \hh_{\m} + \m_+$
be a triangular decomposition of $\m$ so that $\b_{\m} =
\hh_{\m} + \m_+$ is a Borel subalgebra of $\m$ and $$\b =
\b_{\m} + \a + \n$$ is a Borel subalgebra of $\g$.

Let $\Lambda\s \hh^*$ be the set of dominant integral linear forms
on $\hh$, with respect to $\b$ and, for each $\lambda\in \Lambda$, let
$\pi:\g\to End\,V_{\lambda}$ be an irreducible representation with
highest weight $\lambda$. Let $0\neq v_{\lambda}\in V_{\lambda}$ be
a highest weight vector. It follows from the Lie algebra Iwasawa
decomposition $\g = \k + \a +\n$ that $U(\g) = U(\k)U(\a+\n)$ using
the standard notation for universal enveloping algebras. However if
$\lambda\in \Lambda$, then $\Bbb C v_{\lambda}$ is stable under
$U(\a+\n)$ and hence $V_{\lambda}$ is a cyclic $U(\k)$-module. In
fact $$V_{\lambda} = U(\k)\,v_{\lambda}$$ Let $L_{\lambda}(\k)$ be
the (left ideal) annihilator of $v_{\lambda}$ in $U(k)$. It is then
an elementary fact that if $Z$ is any irreducible $\k$-module, one
has $$\hbox{ multiplicity of $Z$ in $V_{\lambda} = $ dim
$Z^{L_{\lambda}(\k)}$}\eqno (0.1)$$ where $Z^{S}= \{w\in Z\mid
S\cdot w = 0\}$ for any subset $S\s U(\k)$. 

The equation (0.1) becomes a useful branching law as soon as one can
explicitly determine generators of $L_{\lambda}(\k)$. It is one of
the main results of this paper to explicitly exhibit such
generators. 
 
Let $\Delta\s \hh^*$ be the set of
roots for $\g$ and let 
$\Delta_+$ be the set of positive roots corresponding to $\b$. For
any
$\varphi\in \Delta$ let $e_{\varphi}\in \g$ be a corresponding root
vector. Let $\ell = rank\,\g$ and let $\Pi =
\{\alpha_1,\ldots,\alpha_{\ell}\}$. Let $I = \{1,\ldots,\ell\}$ and
let $\{h_i\},\,i\in I$, be the basis of $\hh$ such that $h_i\in
\g_i$ where $\g_i$ is the TDS generated by $e_{\alpha_i}$ and
$e_{-\alpha_i}$ and 
$[h_i,e_{\alpha_i}]= 2 e_{\alpha_i}$. Then there uniquely exists
$\lambda_i\in
\Lambda,\,i\in I$, such that $\lambda_i(h_j) = \delta_{ij}$ for
$i,j\in I$. In addition, every element $\lambda\in \Lambda$ may be
uniquely written $$\lambda =\sum_{i\in I}n_i(\lambda)\,\lambda_i$$
where
$n_i(\lambda)\in \Bbb Z_+$. 

{\bf 0.2.} Now let $I = I_{\n}\cup I_{\m}$ be the partition defined so that
$i$ is in $I_{\n}$ or $I_{\m}$ according as $e_{\alpha_i}\in \n$ or
$e_{\alpha_i}\in \m_+$. Let $\theta$ be the complex Cartan
involution of $\g$ so that $\k$ is the space of $\theta$ fixed
points. For any $i\in I_{\n}$ let $z_i = e_{-\alpha_i} +
\theta(e_{-\alpha_i})$. One of course has that $z_i\in \k$. The set
$I_{\n}$ partitions into a disjoint union $I_{\n} = I_{s} \cup
I_{nil}$ where $I_{s} = \{i\in I_{\n}\mid \theta(e_{-\alpha_i})\in
\Bbb C e_{\alpha_i}\}$ and $I_{nil}$ is the complement of $I_s$ in
$I_{\n}$. One readily notes that $z_i$ is semisimple if $i\in I_s$
and $z_i$ is nilpotent if $i\in I_{nil}$. In the latter case one
easily proves that $[e_{-\alpha_i}, \theta(e_{-\alpha_i})]= 0$. For
any $\lambda\in \Lambda$ and  $i\in I_{\n}$ let
$q_{\lambda,i}(t)\in \Bbb C[t]$ be the polynomial of degree
$n_i(\lambda)+ 1$ defined by putting
$$q_{\lambda,i}(t)= (t-n_j(\lambda))(t-n_j(\lambda)+2)\cdots (t+
n_i(\lambda))\eqno (0.2)$$ if $i\in I_s$ and $$q_{\lambda,i}(t)
= t^{n_i(\lambda)+ 1}\eqno (0.3)$$ if $i\in I_{nil}$. One readily
proves that if $i\in I_s$, then $e_{-\alpha_i}$ may be normalized
so that, with respect to the Killing form inner product, $(z_i,z_i) =
(h_i,h_i)$. In that case $q_{\lambda,i}(t)$ is the characteristic
polynomial of
$\pi_{\lambda}(z_i)$ in the irreducible $\g_i$-module
$U(\g_i)\,v_{\lambda}$. We assume that $e_{-\alpha_i}$ is so
normalized. No normalization is needed if $i\in I_{nil}$. The
following is our first main result. \vskip 1pc {\bf Theorem 0.1.}
{\it Let $\lambda\in \Lambda$. Let $\{y_j\},\,j=1,\ldots,k$, be
a basis of $\hh_{\m}$ and let ${\cal G}_{\lambda}(m)$ be the finite
subset of $U(\m)$ defined so that $${\cal G}_{\lambda}(\m) =
\{e_{-\alpha_i}^{n_i(\lambda)+ 1}\mid i\in I_{\m}\}\,\,\cup\,\,
\{y_j-\lambda(y_j)\mid j=1,\ldots,k\}\,\,\cup\,\, \{e_{\alpha_i}\mid
i\in I_{\m}\}$$ and if $\n_*$ is the Killing form complement of $\m$
in
$\k$, let ${\cal G}_{\lambda}(\n_*)$ be the finite set of polynomials
of elements $\n_*$ defined so that $${\cal G}_{\lambda}(\n_*) =
\{q_{\lambda,i}(z_i)\mid i\in I_{\n}\}$$ Put ${\cal G}_{\lambda} =
{\cal G}_{\lambda}(\m)\cup {\cal G}_{\lambda}(\n_*)$. Then the left
ideal annihilator $L_{\lambda}(\k)$ in $U(\k)$ of the highest weight
vector $v_{\lambda}$ of $V_{\lambda}$ is the left ideal generated
by ${\cal G}_{\lambda}$. In particular if $Z$ is any irreducible
$\k$-module, then $$\hbox{ multiplicity of $Z$ in $V_{\lambda} = $ 
dim
$Z^{{\cal G}_{\lambda}}$}\eqno (0.4)$$} \vskip 1pc  By separating the
conditions imposed by ${\cal G}_{\lambda}(\m)$ and ${\cal
G}_{\lambda}(\n_*)$ we may express (0.4) in a simpler form (see
Theorem 0.5 below). Let $\Lambda_{\m}\s
\hh_{\m}^*$ be the set of all dominant integral linear forms on
$\hh_{\m}$ with respect to
$\b_{\m}$ (see (178)). If
$\lambda\in
\Lambda$, then clearly
$\Lambda|\hh_{\m}\in \Lambda_{\m}$. For any $\nu\in \Lambda_{\nu}$
let $\tau_{\nu}:\m\to Aut \,Y_{\nu}$ be an irreducible
$\m$-module with highest weight $\nu$. If $Z$ is a finite
dimensional irreducible $\k$-module then, as a $\m$-module, let 
$Z[\nu]$ be the primary $\tau_{\nu}$ component of $Z$. The
subspace of $\b_{\m}$ highest weight vectors in $Z[\nu]$ is just
$(Z[\nu])^{m_+}$. \vskip 1pc  {\bf Theorem 0.5.} {\it Let
$\lambda\in \Lambda$ and let $Z$ be any finite dimensional
irreducible $\k$-module. Then $$\hbox{ multiplicity of $Z$ in
$V_{\lambda} = $ dim $\,\{w\in (Z[\lambda |\hh_{\m}])^{m_+}\mid
q_{\lambda,i}(z_i)\,w = 0,\forall i\in I_{\n}$}\}\eqno (0.5)$$}

{\bf 0.3.} Another consequence of Theorem 0.1 is a generalization of the Cartan-Helgason
theorem. In the generalization an arbitrary irreducible representation of $M$
replaces the trivial representation of $M$. For $j\in I_s$ let
$\varepsilon_j = exp\,\pi i\, h_j$ and let $F_s$ be the (abelian)
group generated by
$\{\varepsilon_j\},\,j\in I_s$. The following structure
theorem concerning the disconnectivity of $M$ does not appear to be
in the literature. It however was known to David Vogan who points out
that it is implicit in [V]. Theorem 0.6 follows from Theorem 0.1. Let
$\ell_s = card\,I_s$ and let $M_e$ be the identity component of $M$. \vskip 1pc {\bf Theorem
0.6.} {\it One has a group isomorphism $$ F_s \cong \Bbb Z_2^{\ell_s}\eqno (0.6)$$
Furthermore $F_s\s M$ and $M$ has the following product structure:
$$M = F_s\times M_e\eqno (0.7)$$}\vskip 1pc Let $\widehat {F}_s$ be
the character group of $F_s$. Let
$\Lambda_{M_e}$ be the set of all
$\nu\in \Lambda_{\m}$ such that $\tau_{\nu}|\m_o$ exponentiates to an
(irreducible ) representation of
$M_e$. If $\lambda\in \Lambda$ then clearly $\lambda|\hh_{\m} \in \Lambda_{M_e}$. Extending
its previous use, let
$\tau_{\nu}$, for
$\nu\in
\Lambda_{M_e}$, also denote the representation of $M_{e}$ which arises by
exponentiating $\tau_{\nu}|\m_o$. Then $\widehat{F}_{s}\,\times\, \Lambda_{M_e}$
parameterizes the set
$\widehat {M}$ of all equivalences classes of irreducible representations (the unitary
dual) of
$M$. For each
$(\zeta,\nu)\in \widehat{F}_{s}\,\times\, \Lambda_{M_e}$ let
$$\tau_{\zeta,\nu}:M_e\to Aut\, Y_{\zeta,\nu}$$ be the irreducible representation where,
as a vector space, $Y_{\zeta,\nu} = Y_{\nu}$, but with the $M$-module structure defined so
that
$$(\varepsilon,a)\,y = \zeta(\varepsilon)\,\tau_{\nu}(a)\,y$$ for
$(\varepsilon,a)\in F_{s}\,\times\, M_e = M$ and $y\in Y_{\zeta,\nu}$. By abuse of
notation we will take $$\widehat {M} =\{\tau_{\zeta,\nu}\},\,(\zeta,\nu)\in
\widehat{F}_{s}\,\times\,
\Lambda_{M_e}$$ For any $\lambda\in \Lambda$ let 
$\lambda|F_{s}\in
\widehat{F}_{s}$ be (uniquely) defined so that, for $j\in I_{s}$, one has  
$(\lambda|F_{s})(\varepsilon_j) = 1$ if $n_j(\lambda)$ is even and
$(\lambda|F_{s})(\varepsilon_j) = -1$ if
$n_j(\lambda)$ is odd. 

If $\lambda\in \Lambda$ then $U(\m)\,v_{\lambda}$ is an irreducible $M$-module with respect
to $\pi_{\lambda}|M$. The significance of the pair $(\lambda|F_s,\lambda|\hh_{\m}) \in
\widehat{F}_{s}\,\times\, \Lambda_{M_e}$ is that one has a module equivalence
$$Y_{\lambda|F_s,\lambda|\hh_{\m}}\cong U(\m)\,v_{\lambda}\eqno (0.8)$$ of irreducible
$M$-modules. One readily shows that the map $$\Lambda\to \widehat{F}_{s}\,\times\,
\Lambda_{M_e},\qquad \lambda\mapsto (\lambda|F_s,\lambda|\hh_{\m})\eqno (0.9)$$ is
surjective. For any $(\zeta,\nu)\in \widehat{F}_{s}\,\times\, \Lambda_{M_e}$, let ${\cal
F}_{\zeta,\nu}\s \Lambda$ be the fiber of (0.9) over $(\zeta,\nu)$. If $\zeta_1$ is the
trivial character of $F_s$, so that $\tau_{\zeta_1,0}$ is the trivial representation
of $M$, then the Cartan-Helgason theorem asserts that $${\cal F}_{\zeta_1,0} =
\Lambda_{sph}\eqno (0.10)$$ where $$\Lambda_{sph}= \{\lambda\in \Lambda\mid
\pi_{\lambda}|\k\,\hbox{contains the trivial representation of $\k$}\}\eqno (0.11)$$ One may
introduce a partial order in $\Lambda$ by declaring that $\lambda^b \gg \lambda^a$ if there
exists a (necessarily unique and surjective) $\k$-map $s:V_{\lambda^b}\to V_{\lambda^a}$
such that $s(v_{\lambda^b}) = v_{\lambda^a}$. In such a case one necessarily has
$$\hbox{multiplicity of $Z$ in $V^a \leq\,\,$  multiplicity of
$Z$ in $V^b$}\eqno (0.12)$$ for all irreducible $K$-modules $Z$. We may express (0.11) by
saying that ${\cal F}_{\zeta_1,0}= \{\lambda\in \Lambda\mid \lambda\gg 0\}$. Noting
that $\lambda_{\zeta_1,0} = 0$ in the next theorem the generalization of the
Cartan-Helgason theorem is the statement (see Theorem 3.12)
\vskip 1pc {\bf Theorem 0.7.} {\it Let $(\zeta,\nu)\in \widehat{F}_{s}\,\times\,
\Lambda_{M_e}$. Then there exists a unique minimal element $\lambda_{\zeta,\nu}\in
{\cal F}_{\zeta,\nu}$. In fact ${\cal F}_{\zeta,\nu}= \{\lambda\in \Lambda\mid \lambda \gg
\lambda_{\zeta,\nu}\}$. In particular $$\hbox{multiplicity of $Z$ in
$V_{\lambda_{\zeta,\nu}}\leq\,\,$  multiplicity of
$Z$ in $V_{\lambda}$}\eqno (0.12)$$ for all irreducible $K$-modules $Z$ and all
$\lambda\in {\cal F}_{\zeta,\nu}$. Furthermore $${\cal F}_{\zeta,\nu} =
\lambda_{\zeta,\nu} + \Lambda_{sph}\eqno (0.13)$$}\vskip 1pc Explicit formulas for
$\lambda_{\zeta,\nu}$ and the elements in $\Lambda_{sph}$ are given in (193) and (187).
\vskip 1pc {\bf 0.4.} Let $D =  \widehat{F}_{s}\,\times\,\Lambda_{M_e}$ so that $\widehat
{M} =
\{\tau_{\delta}\mid \delta\in D\}$. If $\phi$ is a smooth scalar or vector valued function 
on $G$ and  $g\in G,\,z\in \g$, then $g\cdot \phi$ and $z\cdot \phi$ is again such a
function where the action is defined by (200) and (201). The set of principal series
representations of $G$ are parameterized by $D\times \a^*$. If $(\delta,\xi)\in D\times
\a^*$, then the corresponding Harish-Chandra module, $H(Y_{\delta},\xi)$, is a 
$(\g,K)$-module of $Y_{\delta}$ valued functions on $G$ satisfying (202) and (203). Rather
than dealing with vector valued functions on $G$ we can, without loss, use scalar valued
functions on $G$, by use of the Borel-Weil theorem applied to $M$. The space
$C^{\infty}(G)$ is a (left) module for $U(\g)$ when the latter operates as left invariant
differential operators on $G$. The action is denoted by $f\cdot u$ for $f\in
C^{\infty}(G)$ and $u\in U(\g)$. If $u\in \g_o$ then $f\cdot u$ is given by (205). If
$(\delta,\xi)\in D\times \a^*$, then the $(\g,K)$-module $H(\delta,\xi)$ of scalar valued
functions on $G$ is defined in (218). A $(\g,K)$-equivalence of $H(\delta,\xi)$ and
$H(Y_{\delta},\xi)$ is established in Proposition 4.4.

A result of Casselman (see [C] or [BB]) asserts that any irreducible Harish-Chandra module
$H$ embeds in $H(\delta,\xi)$ for some $(\delta,\xi)$. If $H = V_{\lambda}$ for $\lambda\in
\Lambda$, then Wallach has shown (see [W], \S 8.5) that there exists a unique
$(\delta,\xi)\in  D\times \a^*$ for which one has a $(\g,K)$ injection $$V_{\lambda}\to
H(\delta,\xi)\eqno (0.14)$$ and the map (0.14) is unique up to scalar multiplication. For
completeness and to simplify the proof of Theorem 0.8 below we reprove, using our scalar
valued functions instead of vector valued functions, the Wallach result. The explicit
determination of
$(\delta,\xi)$ is given in Propositions 4.16 and 4.17. Let
${\cal V}_{\lambda}$ be the image of (0.14). The functions in $H(\delta,\xi)$ are of
course determined by their restrictions to $K$. Theorem 0.8 (proved in this paper as
Theorem 4.19) is a noncompact analogue of the Borel-Weil theorem in the sense that
${\cal V}_{\lambda}$ will turn out be a finite dimensional space in $H(\delta,\xi)$ of
solutions of certain differential equations on $G$ arising from certain elements 
in $U(\k)$ operating as left invariant differential operators. If one was dealing with the
compact form of
$G_{\Bbb C}$, the Borel-Weil theorem asserts the equations would be of Cauchy-Riemann type.
\vskip 1pc {\bf Theorem 0.8.} {\it Let
$\lambda
\in
\Lambda$ and let
$\delta
\in D$ and
$\xi\in a^*$ be defined as in (0.14) (see (246) and (248)). Then ${\cal
V}_{\lambda}$ may be given by $${\cal V}_{\lambda} = \{f\in H(\delta,\xi)\mid f\cdot
q_{\lambda^c,i}(z_i)= 0,\,\,\forall i\in I_{\n}\}\eqno (0.15)$$ where $z_i\in
\k$ and $q_{\lambda^c,i}(t)\in \Bbb C[t]$, for $i\in I_{\n}$, are defined as in \S 0.2, and
$\lambda^c\in \Lambda$ is the highest weight of the dual module $V_{\lambda}^*$.}
 \vskip 1.5pc

\centerline{\bf 1. The generators of the left ideal $L_{\lambda}(\k) $}\vskip
1pc
 
{\bf 1.1.} Let $\g$ be a complex semisimple Lie algebra and let $(x,y)$ be
the Killing form $B_{\g}$ on $\g$. If $\r$ is any Lie subalgebra of
$\g$, let $U(\r)$ be the universal enveloping algebra of $\r$. Of
course we can regard $U(\r)\s U(\g)$.

 Let
$\hh\s
\b$ be, respectively, a Cartan subalgebra and a Borel subalgebra of
$\g$. The restriction $B_{\g}|\hh$ defines an isomorphism
$$\eta:\hh\to\hh^*\eqno (1)$$ For $\mu,\nu\in \hh^*$, the dual space to $\hh$, let
$(\mu,\nu)$ be the bilinear form on $\hh^*$ defined so that $(\eta(x),\eta(y))
= (x,y)$ for any $x,y\in \hh$. Let
$\ell = dim\,\hh$ and let
$\Delta\s\hh^*$ be the set of roots for the pair $(\hh,\g)$. Let a
complete set of root vectors
$\{e_{\varphi}\},\,\varphi\in \Delta$, be chosen. If $\ee\s\g$ is
any subspace which is stable under the action of
$ad\,\hh$, let $\Delta(\ee)=
\{\varphi\in \Delta\mid e_{\phi}\in
\d\}$. A choice of a set $\Delta_+$ of positive roots is made by
putting $\Delta_+ = \Delta(\b)$. Let $\Pi=
\{\alpha_1,\ldots,\alpha_{\ell}\}\s \Delta_+$ be the set of simple
positive roots. 

Let $\Lambda\s \hh^*$ be the semigroup
of dominant integral weights with respect to $\b$. One knows there uniquely
exists $\lambda_i\in\Lambda,\,i=1,\ldots,\ell$, such that
$$2(\lambda_i,\alpha_j)/(\alpha_j,\alpha_j)= \delta_{ij}\eqno(2)$$
for
$i,j\in \{1,\ldots,\ell\}$ and that any $\lambda\in \Lambda$ can
be uniquely written $$\lambda = \sum_{i=1}^{\ell}
n_i(\lambda)\,\lambda_i$$ where $n_i(\lambda)\in \Bbb
Z_+$. For $i=1,\ldots,\ell$, let $\g_i$ be the TDS generated by
$e_{\alpha_i}$ and $e_{-\alpha_i}$. Then there exists a unique
scalar multiple, $h_i$, of $[e_{\alpha_i},e_{-\alpha_i}]$ such
that $[h_i,e_{\alpha_i}] = 2e_{\alpha_i}$ and $[h_i,e_{-\alpha_i}] =
-2e_{-\alpha_i}$. One has that $\Bbb C h_i = \g_i\cap \hh$ and
$$\eta(h_i) = 2\alpha_i/(\alpha_i,\alpha_i)\eqno (3)$$ With respect to the
pairing of $\hh$ and $\hh^*$, the equations (2) imply that $$\langle
\lambda_i,h_j\rangle = \delta_{ij}\eqno (4)$$   

 Let $\b_-$ be the Borel subalgebra of $\g$ which contains $\hh$
and is such that $\Delta(\b_-) = \Delta_-$ where $\Delta_- =
-\Delta_+$. Let $\u$ and $\u_-$, respectively, be the nilradicals
of $\b$ and $\b_-$. One then has the linear space direct sum 
$\g = \u_- + \b$. The Poincar\'e-Birkhoff-Witt theorem then implies
that
$$U(\g) = U(\u_-)U(\b)\eqno (5)$$ 

For each
$\lambda\in \Lambda$ let $$\pi_{\lambda}:U(\g)\to End\,
V_{\lambda}$$ be an irreducible representation with highest weight
$\lambda$. Also let $0\neq v_{\lambda}\in V_{\lambda}$ be a highest
weight vector. If $u\in U(\g)$ and $v\in V_{\lambda}$ and there is
no danger of confusion, we will occasionally write $u\,v$ or $u\cdot v$ in
$V_{\lambda}$ for $\pi_{\lambda}(u)v$. 

Let $\lambda\in \Lambda$. Since the 1 dimensional subspace $\Bbb C
v_{\lambda}$ is stable under $U(\b)$, it follows from (5) that $V_{\lambda}$ is a
cyclic module for $U(\u_-)$ and in fact $$V_{\lambda} = U(\u_-)v_{\lambda}\eqno
(6)$$ Let $L_{\lambda}(\u_-)$ be the annihilator of $v_{\lambda}$ in $U(\u_-)$ so
that
$L_{\lambda}(\u_-)$ is a left ideal in $U(\u_-)$ and, as $U(\u_-)$-modules, one has the
isomorphism $$V_{\lambda} = U(\u_-)/L_{\lambda}(\u_-)\eqno (7)$$ Now let $i\in
\{1,\ldots,\ell\}$. Obviously $e_{\alpha_i}v_{\lambda}=0$ and $h_iv_{\lambda} =
n_i(\lambda)v_{\lambda}$ by (4) using the notation of \S1.1. The following
is well known and is a consequence of the representation theory of $Sl(2,\Bbb
C)$. \vskip 1pc {\bf Proposition 1.1.} {\it The cyclic $U(\g_i)$-submodule of
$V_{\lambda}$ generated by $v_{\lambda}$ is $U(\g_i)$-irreducible and has
dimension $n_i(\lambda) + 1$. In addition $n_i(\lambda)$ is the smallest
nonneqative integer $k$ such that $e_{-\alpha_i}^{k +1}v_{\lambda}=0$.}\vskip
1pc It follows from Proposition 1 that $e_{-\alpha_i}^{n_i(\lambda)+1}\in
L_{\lambda}(\u_-)$. A proof of the following theorem may be found in [PRV]. See Lemma
2.1 in [PRV]. For more details about this result see Remark 44 (especially
(170)), p. 315 in [K3]. \vskip 1pc {\bf Proposition 1.2.} {\it Let
$\lambda\in
\Lambda$. Then the elements $e_{-\alpha_i}^{n_i(\lambda)+1},\,i=1,\ldots,\ell$,
generate $L_{\lambda}(\u_-)$. That is, $$L_{\lambda}(\u_-) = \sum_{i=1}^{\ell}
U(\u_-)e_{-\alpha_i}^{n_i(\lambda)+1}\eqno (8)$$}\vskip 1pc {\bf 1.2.} Let $\g_o$
be a real form of $\g$ and let $$\g_o = \k_o + \p_o\eqno (9)$$ be a Cartan
decomposition of
$\g_o$. Thus $Ad_{\g_o}\k_o$ is a maximal compact subgroup of $Ad_{\g_o}\g_o$
and $\p_o$ is the $B_{\g}$-orthocomplement of $\k_o$ in $\g_o$. Let $\theta$ be
the involution on $\g$ such that $\theta |\g_o$ is the Cartan involution of
$\g_o$ corresponding to (9). Let $\a_o$ be a maximal abelian subalgebra
contained in $\p_o$ and let $\m_o$ be the centralizer of $\a_o$ in $\k_o$. Then
one knows there exists $w\in \a_o$ whose centralizer in $\g_o$ is $\m_o + \a_o$.
The operator $ad\,w \vert \g_o$ is diagonalizable (and hence has real
eigenvalues). Let $\n_o$ be the span of eigenvectors of $ad\,w$ in $\g_o$
corresponding to positive eigenvalues. Then $\g_o = \k_o + \a_o +\n_o$ is,
infinitesimally, an Iwasawa decomposition of $\g_o$. We will denote the
complexification of the Lie algebras introduced above by dropping the subscript
$o$. In particular,
$$\g = \k + \a + \n\eqno (10)$$ is, infinitesimally, a complexified Iwasawa
decomposition of $\g$. Also $\q = \m + \a + \n$ is a parabolic Lie
subalgebra of $\g$ and $\m +\a, \n$ are, respectively, a Levi factor and
nilradical of $\q$. Also if $\g^w$ is the centralizer of $w$ in $\g$, then
$$\g^w = \m + \a\eqno (11)$$ Let
$(\hh_{\m})_o$ be a Cartan subalgebra of the reductive Lie algebra $\m_o$ and
let $\hh_{\m}$ be the complexification of $(\hh_{\m})_o$. Then
$\hh_{\m} + \a$ is a Cartan subalgebra of
$\g$. We will fix the Cartan subalgebra $\hh$ of $\g$ in \S1.1 so that $$\hh=
\hh_{\m} + \a\eqno (12)$$ Next the Borel subalgebra $\b$ of \S1.1 will be fixed so
that
$\hh\s\b $ (an assumption in \S1) and such that $\b\s \q$. It follows
immediately that $$\a + \n\s \b\eqno (13)$$ But then as a consequence of (10),
(13) and the Poincar\'e-Birkhoff-Witt theorem one has $$U(\g) =
U(\k)U(\b)\eqno (14)$$

A Lie subalgebra of $\g$ is called symmetric if it is the set of fixed
vectors for an involutory automorphism of $\g$. In particular $\k$ is a
symmetric subalgebra of $\g$. In fact $\k$ is the most general symmetric
subalgebra of $\g$. Let
$\lambda\in
\Lambda$ and let
$L_{\lambda}(\k)$ be the (left ideal) annihilator in $U(\k)$ of $v_{\lambda}$ in
$V_{\lambda}$. The equality (14) readily implies
\vskip 1pc {\bf Proposition 1.3.} {\it Any finite dimensional irreducible
$\g$-module, where $\g$ is a complex semisimple Lie algebra, is a cyclic module
for any symmetric subalgebra of $\g$. More specifically, in the notation above,
let
$\lambda\in \Lambda$. Then $V_{\lambda}$ is a cyclic $U(\k)$-module and
$v_{\lambda}$ is a cyclic generator. In particular, as $U(\k)$-modules, one has
an isomorphism $$V_{\lambda} \cong U(\k)/L_{\lambda}(\k)\eqno (15)$$}

 \vskip .1pc
{\bf 1.3.} If
$Z$ is any $U(\k)$-module and $S\s U(\k)$ is any subset, let $Z^S = \{z\in Z\mid
Sz= 0\}$. Let $Z$ be any $U(\k)$-module and let $\lambda\in \Lambda$. If
$\sigma\in Hom_{\k}(V_{\lambda},Z)$, then obviously $\sigma(v_{\lambda}) \in
Z^{L_{\lambda}(\k)}$. Conversely if $z\in Z^{L_{\lambda}(\k)}$ then by (15)
there clearly exists a unique element $\sigma\in Hom_{\k}(V_{\lambda},Z)$ such
that $\sigma(v_{\lambda})= z$. That is, one has \vskip 1pc {\bf Lemma 1.4.} {\it
Let $\lambda\in \lambda$ and let $Z$ be a $U(\k)$-module. Then the map
$$Hom_{\k}(V_{\lambda},Z)\to Z^{L_{\lambda}(\k)},\quad \sigma\mapsto
\sigma(v_{\lambda})\eqno (16)$$ is a linear isomorphism.} \vskip 1pc  Of
course $\k$ is a reductive Lie subalgebra of $\g$. The assumption that a
finite dimensional $\k$-module $Z$ is completely reducible reduces to
the assumption that $Z$ is completely reducible as a $Cent\,\k$-module. Lemma
1.4. leads to the following branching principle for symmetric subalgebras of
semisimple Lie algebras.
\vskip 1pc {\bf Proposition 1.5.} {\it Let $Z$ be a finite dimensional completely
reducible $\k$-module. Let $\lambda\in \Lambda$. Then there is a nonsingular
pairing of $Hom_{\k}(V_{\lambda},Z)$ and $Hom_{\k}(Z,V_{\lambda})$ so that 
(recalling (16)) one has $$dim\,\,Hom_{\k}(Z,V_{\lambda}) =
dim\,\,Z^{L_{\lambda}(\k)}\eqno (17)$$ In particular if $Z$ is $\k$-irreducible,
then the $$\hbox{multiplicity of $Z$ in $\pi_{\lambda}|\k =
dim\,Z^{L_{\lambda}(\k)}$}\eqno (18)$$}

\vskip .5pc {\bf Proof.} We may identify
$Hom_{\k}(V_{\lambda},Z)$ with $(V_{\lambda}^*\otimes Z)^{\k}$. On the other
hand, the dual space to $V_{\lambda}^*\otimes Z$ may be written as
$Z^*\otimes V_{\lambda}$. The latter however identifies with
$Hom\,(Z,V_{\lambda})$. Since the tensor product action of $\k$ is
clearly completely reducible, the nonsingular pairing of $V_{\lambda}^*\otimes
Z$ and $Z^*\otimes V_{\lambda}$ restricts to a nonsingular pairing of
$(V_{\lambda}^*\otimes Z)^{\k}$ with
$(Z^*\otimes V_{\lambda})^{\k} = Hom_{\k}(Z,V_{\lambda})$. QED \vskip 1pc
{\bf Remark 1.6.} Proposition 1.5 is useful to establish branching laws
for symmetric Lie subalgebras $\k$ of $\g$ if one can determine explicit
generators of $L_{\lambda}(\k)$ for any $\lambda\in \Lambda$. Such a determination
will be the one of the main results in this paper. \vskip 1pc {\bf 1.4.} The
bilinear form
$B_{\g}|(\m +
\a)$ is nonsingular since $\m +
\a$ is a Levi factor of a parabolic subalgebra ($\q$) of $\g$. Let $\ss$ be the
$B_{\g}$-orthocomplement of $\m +\a$ in $\g$ so that $$\g = \m + \a +
\ss\eqno (19)$$ Obviously
$[\m +\a,\ss]\s
\ss$. In particular $[\hh,\ss]\s \ss$, since $\hh\s \m +\a$, so that $$\ss =
\sum_{\varphi\in \Delta(\ss)}\Bbb C e_{\varphi}$$ It follows from (11)
that $$\Delta(\ss) = \{\varphi\in \Delta\mid \langle \varphi,w\rangle \neq
0\}\eqno (20)$$ In particular $$\Delta(\ss) = -\Delta(\ss)\eqno (21)$$
Obviously $\n = \ss\cap \q$. But in fact $$\n = \ss\cap \b\eqno
(22)$$ This follows easily from the fact that
$\n$ must be contained in $\b$ since $\b$ is maximal solvable and $\b$
normalizes $\n$ because $\b\s\q$. It follows from (22) that $$\n =
\sum_{\varphi\in \Delta_+\cap \Delta(\ss)}\Bbb C e_{\varphi}\eqno (23)$$ But (21)
and (22) also imply that one has a linear direct sum
$$\ss =\n +\n_-\eqno (24)$$ where $$\n_- = \sum_{\varphi\in \Delta_+\cap
\Delta(\ss)}\Bbb C e_{-\varphi}\eqno (25)$$ Let $\q_- = \m +\a + \n_-$. \vskip
1pc {\bf Proposition 1.7.} {\it The subspace $\q_-$ is a parabolic subalgebra of
$\g$ with nilradical $\n_-$ and Levi factor $\m +\a$.}

\vskip 1pc {\bf
Proof.} Recalling the definition of $w\in \a_o$ in \S 1.2 it follows from
(11), (21) and (25) that (a) $\q$ is the $\Bbb C$-span of all eigenvectors of
$ad\,w$ corresponding to nonpositive eigenvalues, (b) $\m +\a$ is the kernel of
$ad\,w$ and (c) $\n_-$ is the $\Bbb C$-span of all eigenvectors of $ad\,w$
corresponding to negative eigenvalues. QED\vskip 1pc The inclusion $\b\s \q$
implies of course that $\u\s \q$ using the notation of \S1.2. But $\n\s \u$ by
(23). However $\u\cap (\m + \a) = \u\cap \m$ since $\a$ is central in the Levi
factor
$\m +\a$. Let $\m_+ = \u\cap \m$ so that one has the Lie algebra semi-direct sum 
$$\u = \m_+ +  \n\eqno (26)$$ and hence one has the disjoint union $$\Delta_+ =
\Delta(\m_+) \cup \Delta(\n)\eqno (27)$$ On the other hand, $\Delta(\m) =
-\Delta(\m)$ since $B_{\g}|\m$ is clearly nonsingular. Thus if $\m_- =
\u_-\cap\m$ one has $$\Delta(\m_-) = - \Delta(\m_+)\eqno (28)$$ In addition $$\m
= \m_- + \hh_{\m} + \m_+\eqno (29)$$ is a ``triangular decomposition" of the
reductive Lie algebra $\m$. Indeed $$\g = \u_- + \hh +\u \eqno (30)$$ is a
``triangular decomposition" of $\g$ and the three components on the right side
of (29) are the respective intersections of $\m$ with the three components on
the right side of (30). One notes that $\m$ is the sum of the respective three
intersections since $[\a,\m] = 0$ and all three components on the right side of
(30) are stable under $ad\,\a$. Let $\b_{\m} = \hh_{\m} + \m_+$. Then $\b_{\m}$
is a Borel subalgebra of $\m$ since (29) is a triangular decomposition of $\m$.
Furthermore $$\b_{\m} = \b\cap \m $$ since (1) $\b_{\m}$ is obviously contained
in $\b\cap \m$ and (2) a Borel subalgebra is maximal solvable and $\b\cap \m$ is
solvable. 
Taking the negative of the roots in (27) it follows from (23), (25) and (28) that
one has the disjoint union $$\Delta_- = \Delta(\m_-) \cup \Delta(\n_-)\eqno
(31)$$  and consequently the linear space direct sum $$\u_- = \m_- + \n_-\eqno
(32)$$ Recalling Proposition 1.7 let
$\v$ be the Lie subalgebra of $\q_-$ defined by putting $\v = \m + \n_-$. Then by
(24) one has the linear direct sum $$\g = \v +\a +\n\eqno (33)$$ On the other hand
\vskip 1pc {\bf Proposition 1.8.} {\it One has the linear space direct sum
$$\v = \u_- + \b_{\m}\eqno (34) $$ and the map $$U(\u_-)\otimes U(\b_{\m}) \to
U(\v), \qquad u\otimes v \mapsto u\,v\eqno (35)$$ is a linear
isomorphism.}

\vskip 1pc {\bf Proof.} The fact that the right side of (34) is a
direct sum is immediate since $\b_{\m}\s \b$. On the other hand, by (29) and (32)
$$\eqalign{\u_- + \b_{\m} &= \n_- + \m_- + \b_{\m}\cr
&= \n_- + \m\cr
&= \v\cr}$$ The statement (35) is an immediate consequence of the 
Poincar\'e-Birkhoff-Witt theorem. QED\vskip 1pc {\bf 1.5.} Let $\lambda\in
\Lambda$. Comparing (10) and (33) the argument leading to Proposition 1.3 applies
as well to
$U(\v)$. Namely, one has that $V_{\lambda}$ is a cyclic module for $U(\v)$ with
cyclic generator $v_{\lambda}$. Let $L_{\lambda}(\v)$ be the annihilator of
$v_{\lambda}$ in $U(\v)$. As a preliminary step in determining the left ideal
$L_{\lambda}(\k)$ in $U(\k)$ we will now determine $L_{\lambda}(\v)$. It follows from
(27) that the set, $\Pi$, of simple positive roots admits the following
partition:
$$\Pi = \Pi_{\n}\cup \Pi_{\m}\eqno (36)$$ where $\Pi_{\n}=\Pi\cap \Delta(\n)$
and $\Pi_{\m}= \Pi\cap \Delta(\m_+)$. Let $\{y_j\},\,j=1,\ldots,k,$ be a basis
of
$\hh_{\m}$. Let $$\omega_{\lambda}:U(\b_{\m})\to \Bbb C$$ be the character on
$U(\b_{\m})$ defined so that $u\,v_{\lambda} = \omega_{\lambda}(u) v_{\lambda}$ for
all
$u\in U(\b_{\m})$. It is clear that the ideal $Ker\,\omega_{\lambda}$ of codimension
one in $U(\b_{\m})$ is generated by $\{y_j-\lambda(y_j)\},\,j=1,\ldots,k,$ and
$\{e_{\alpha}\},\,\alpha\in \Pi_{\m}$. 
\vskip 1pc {\bf Proposition 1.9.} {\it Let $\lambda\in \Lambda$. Then
$$L_{\lambda}(\v) =
\sum_{\i=1}^{\ell} U(\v)e_{-\alpha_i}^{n_i(\lambda) +1} + \sum_{j=1}^k
U(\v)(y_j-\lambda(y_j)) + \sum_{\alpha\in \Pi_{\m}}U(\v)e_{\alpha}\eqno
(37)$$}\vskip 1pc {\bf Proof.} Now one has the direct sum $$U(\b_{\m}) =
Ker\,\omega_{\lambda}
\oplus \Bbb C\, 1\eqno (38)$$ But, by abuse of notation in Proposition 1.8, we
can then write $$U(\v) = U(\u_-)\otimes U(\b_{\m})\eqno (39)$$ Thus one has the
direct sum
$$U(\v) = U(\u_-)
\oplus U(\u_-)\,Ker\,\omega_{\lambda}\eqno (40)$$ But obviously 
$U(\u_-)\,Ker\,\omega_{\lambda}\s J_{\lambda}$ so that $$L_{\lambda}(\v) =
(U(\u_-)\cap L_{\lambda}(\v)) \oplus U(\u_-)\,Ker\,\omega_{\lambda}\eqno (41)$$ On
the other hand, clearly $U(\u_-)\cap L_{\lambda}(\v) = L_{\lambda}(\u_-)$ using
the notation of Proposition 1.2 so that $$\eqalign{U(\u_-)\cap
L_{\lambda}(\v)&=\sum_{i=1}^{\ell} U(\u_-)e_{-\alpha_i}^{n_i(\lambda) +1}\cr
&\s \sum_{i=1}^{\ell} U(\v)e_{-\alpha_i}^{n_i(\lambda) +1}\cr}\eqno (42)$$ 
However, as noted above, 
$Ker\,\omega_{\lambda}$ is the ideal in
$U(\b_{\m})$ generated by $\{y_j-\lambda(y_j)\},\,j=1,\ldots,k,$ and
$\{e_{\alpha}\},\,\alpha
\in \Pi_{\m}$. But then by (39) the left ideal in $U(\v)$ generated by these
elements is the second summand in (41). Hence the second summand on the right
side of (41) is the same as the sum of the last two sums on the right side of
(37). But then if
$L_{\lambda}'(\v)$ is the left ideal given by the right side of (37), one has
$L_{\lambda}(\v)\s L_{\lambda}'(\v)$ by (42). But obviously  $L_{\lambda}'(\v)\s
L_{\lambda}(\v)$. QED\vskip 1pc {\bf 1.6.} Now if $\r\s\g$ is a subspace stable
under the involution $\theta$ let $\r^{\theta}$, (see \S 1.2) be the space of
$\theta$ fixed vectors in $\r$. Of course $\g^{\theta} = \k$ and $\r^{\theta} =
\r\cap \k$. It is immediate from the orthogonal direct sum (19) that $\m +\a$ and
$\ss$ are stable under $\theta$, and hence if $\n_*= \ss^{\theta}$ then (19) implies
 $$\k = \m + \n_*\eqno (43)$$ One has $\theta(w) = -w$ for the element $w\in \a$
which defines $\n$ and also $\n_-$ in terms of the spectrum of $ad\,w$ (see \S
1.2, (23), and (25)). But then clearly $$\theta(\n_-) = \n\eqno (44)$$ But then, by
(24),
$$\n_* = \{ z + \theta(z)\mid z\in \n_-\}\eqno (45)$$ But now $\hh = \hh_{\m} + \a$
is stable under $\theta$. Hence $\theta$ carries root spaces to root spaces.
Thus we may define an involution $\theta:\Delta\to \Delta$, preserving root
addition when the sum is a root, such that for any
$\varphi\in \Delta$ one has $\Theta (\Bbb C e_{\varphi}) = \Bbb
Ce_{\theta(\varphi)}$. More explicitly, since $\theta = \theta^{-1}$, one readily
has
$$\langle \theta(\varphi),x\rangle = \langle \varphi,\theta(x)\rangle \eqno
(46)$$ for any
$\varphi\in \Delta$ and $x\in \hh$.  With respect to the action of
$\theta$ on
$\Delta_-$, it follows from (31) and (44) that $$\eqalign{&\theta(\Delta(\n_-)) =
\Delta(\n)\cr &\hbox{$\theta$ reduces to the identity map on
$\Delta(\m_-)$}\cr}\eqno (47)$$

\vskip 1pc{\bf Remark 1.10.} 
If $\varphi\in
\Delta(\n_-)$ then
$\theta(\varphi)\in \Delta(\n)$ by (47). However it is not necessarily true that
$\theta(\varphi) = -\varphi$. A condition that $\theta(\varphi) = -\varphi$ is given
in the next proposition. \vskip 1pc For any $\varphi\in \Delta$ let
$h_{\varphi}$ be the unique element in the TDS generated by $e_{\varphi}$ and
$e_{-\varphi}$ such that $[h_{\varphi},e_{\varphi}] = 2e_{\varphi}$ and
$[h_{\varphi},e_{-\varphi}] = -2e_{-\varphi}$. One of course has $h_{\varphi}\in
\hh$. \vskip 1pc {\bf Proposition 1.11.} {\it Let $\varphi\in \Delta_-$. Then
$$[e_{\varphi},\theta( e_{\varphi})]\in \a\eqno (48)$$ Furthermore
$[e_{\varphi},\theta( e_{\varphi})]\neq 0$ if and only if $\theta(\varphi)=
-\varphi$, in which case $[e_{\varphi},\theta( e_{\varphi})]$ is a nonzero
multiple of $h_{\varphi}$, noting that $h_{\varphi}\in \a$ by (48). In any case
$\varphi +
\theta(\varphi)$ is never a root.} \vskip 1pc

{\bf Proof.} By the definition of 
$\theta(\varphi)$ one has $0\neq \theta( e_{\varphi})\in \Bbb C
e_{\theta(\varphi)}$. Let $x\in \a$. Then $[x,\theta( e_{\varphi})] =
-\langle\varphi,x\rangle\, \theta( e_{\varphi})$ by (46) since
$\theta(x) = - x$. Thus if we put $z = [e_{\varphi},\theta( e_{\varphi})]$, then
$$[x,z]= 0\eqno (49)$$  On the other hand
$\theta(z) = - z$ since $\theta$ is involutory. Thus 
$z\in \p$. But then (49) proves (48) since $\a$ is maximally commutative in $\p$
(because $\a_o$ is maximally commutative in $\p_o$ and $\g_o$ is a real form of
$\g$). Assume  $\theta(\varphi)\neq -\varphi$. But then $z\neq 0$ if and only if 
$\varphi + \theta(\varphi)$ is a root. Furthermore in this case $z$ is a nonzero
multiple of
$e_{\varphi + \theta(\varphi)}$. This is impossible since all elements in
$\a$ are semisimple. Thus $z = 0$ and $\varphi + \theta(\varphi)$ is not a root.
If $\theta(\varphi) = - \varphi$, then from the structure of a TDS, it follows
that $z$ is a nonzero multiple of $h_{\varphi}$. In addition $h_{\varphi}\in \a$
by (48). QED\vskip 1pc Let $\u_* = \n_* + \m_-$. But then $$\k = \u_* +
\b_{\m}\eqno (50)$$ is a linear direct sum by (29) and (43). 

{\bf 1.7.} Now simply order
$\Delta_+$ so that $\Delta_+ = \{\varphi_1,\ldots,\varphi_r\}$. Then
$\{e_{-\varphi_1},\ldots,e_{-\varphi_r}\}$ is a basis of $\u_-$. Recall (27) and
(28). If
$\varphi\in \Delta(\n)$ put $z_{-\varphi} = e_{-\varphi} +
\theta(e_{-\varphi})$ and if $\varphi\in \Delta(\m_+)$ put $z_{-\varphi} =
e_{-\varphi}$. \vskip 1pc {\bf Lemma 1.12.} {\it The set
$\{z_{-\varphi_1},\ldots,z_{-\varphi_r}\}$ is a basis of $\u_*$. }\vskip 1pc
{\bf Proof.} This is obvious from (32), (45) and the fact that $\m_-\s\k$.
QED\vskip 1pc 1.7. Now simply order the roots in $\Delta(\m_+)$ so that we can
write $\Delta(\m_+) = \{\mu_1,\ldots, \mu_d\}$. Hence
$\{e_{\mu_1},\ldots,e_{\mu_d}\}$ is a basis of $\m_+$. Consequently, using the
notation of Proposition 1.9 the following lemma follows from (34) and (50).
\vskip 1pc {\bf Lemma 1.13.} {\it The set $\{e_{-\varphi_1},\ldots,e_{-\varphi_r},
y_1,\ldots,y_k,e_{\mu_1},\ldots,e_{\mu_d}\}$ is a basis of $\v$ and the set
$\{z_{-\varphi_1},\ldots,z_{-\varphi_r},
y_1,\ldots,y_k,e_{\mu_1},\ldots,e_{\mu_d}\}$ is a basis of $\k$.} \vskip 1pc
 The elements $y_1,\ldots,y_k,$ are algebraically independent generators of
the polynomial ring $U(\hh_m)$. The same is true for elements $w_i =
y_i-c_i,\,i=1,\ldots,k$, where $c_i\in \Bbb C,\,i=1,\ldots,k,$ are arbitrary
fixed constants. Let
$n = r + k +d$ and regard $\Bbb Z_+^n = \Bbb Z_+^r \times \Bbb Z_+^k\times \Bbb
Z_+^d$. Furthermore if we write $(p,q,s)\in \Bbb Z_+^n$ it will mean that
$p=(p_1,\dots,p_r)\in \Bbb Z_+^r ,\,q=(q_1,\dots,q_k)\in \Bbb Z_+^k$, and 
$s=(s_1,\dots,s_d)\in \Bbb Z_+^d$. Furthermore if
$(p,q,s)\in
\Bbb Z_+^n$ let
$$e(p,q,s) = e_{-\varphi_1}^{p_1}\cdots e_{-\varphi_r}^{p_r}
w_{1}^{q_1}\cdots w_k^{q_k} e_{\mu_1}^{s_1}\cdots e_{\mu_d}^{s_d}$$ and let 
$$z(p,q,s) = z_{-\varphi_1}^{p_1}\cdots z_{-\varphi_r}^{p_r}
w_{1}^{q_1}\cdots w_k^{q_k} e_{\mu_1}^{s_1}\cdots e_{\mu_d}^{s_d}$$
 so that by the Poincar\'e-Birkhoff-Witt theorem one has \vskip 1pc {\bf
Proposition 1.14.} {\it The set $\{e(p,q,s)\},\,(p,q,s)\in \Bbb Z_+^n$, is a basis
of $U(\v)$ and the set $\{z(p,q,s)\},\,(p,q,s)\in \Bbb Z_+^n$, is a basis
of $U(\k)$.} \vskip 1pc {\bf 1.8.} Let $R$ be the span of $\{z(p,0,0)\},\,p\in \Bbb
Z_+^r$. Then it is immediate from Proposition 1.14 that (by abuse of notation) we
may write
$$U(\k) = R\otimes U(\b_{\m})\eqno (51)$$ where $\otimes$ identifies with
multiplication. \vskip 1pc {\bf Remark 1.15.} Note that the set $\{e(p,0,0)\},p\in
\Bbb Z_+^r$, is a basis of $U(\u_-)$. \vskip 1pc
Now let $h\in \hh$ be the unique element such that $\alpha(h) = 1$ for all
$\alpha\in \Pi$. Then the eigenvalues of the diagonalizable operator $ad\,h$ on
$\g$ are integers. Extend the action of $ad\,h$ to $U(\g)$ by commutation. Then
the eigenvalues of
$ad\,h$ on
$U(\g)$ are still integers and $U(\g)$ is still completely reducible under the
action of $ad\,h$. Thus one has the direct sum $$U(\g) =
\sum_{i\in \Bbb Z} U_i(\g)\eqno (52)$$ where $U_i(\g)$ is the eigenspace
of $ad\,h$ corresponding to the eigenvalue $i$. For any $u\in U(\g)$ let $u_i$
be the component of $u$ in $U_i(\g)$ with respect to the decomposition (52). If
$u\neq 0$ let $\delta(u)\in \Bbb Z$ be the minimal value of $i$ such that
$u_i\neq 0$ and put $u_{\delta} = u_{\delta(u)}$ so that $u_{\delta} \neq 0$.
Put $0_{\delta} = 0$ and $\delta(0) = \infty$. One readily verifies the
following properties: If $u,v\in U(\g)$ then $$\eqalign{(uv)_{\delta} &=
u_{\delta}v_{\delta}\cr \delta(uv) &= \delta(u) + \delta(v)\cr}\eqno (53)$$ and if
$\delta(u) <\delta(v)$ one has $$\eqalign{(u + v)_{\delta} &= u_{\delta}\cr
\delta(u+v)&=\delta(u)\cr}\eqno (54)$$ If $\delta(u) =\delta(v)$ and $u_{\delta}
+v_{\delta} \neq 0$, then $$\eqalign{(u+v)_{\delta} &= u_{\delta} + \v_{\delta}\cr
\delta(u + v)&= \delta(u) = \delta(v)\cr}\eqno (55)$$ If $\delta(u)
=\delta(v)< \infty $ and $u_{\delta}
+v_{\delta} = 0$, then $$\delta(u + v) >\delta(u) = \delta(v)\eqno (56)$$ Now any
root vector in $\g$ is of course an eigenvector of $ad\,h$. It follows then
from (47) that for $i=1,\ldots,r$, $$(z_{-\varphi_i})_{\delta} =
e_{-\varphi_i}\eqno (57)$$ Also for $i=1,\ldots,d$, $$(e_{\mu_i})_{\delta}=
e_{\mu_i}$$ Any element in $U(\hh_m)$ is an $ad\,h$ eigenvector with
eigenvalue $0$. Thus for $i=1,\ldots,k,$ one has $$(w_i)_{\delta} = w_i\eqno
(58)$$ But now (53) implies \vskip 1pc {\bf Proposition 1.16.} {\it Let
$(p,q,s)\in \Bbb Z_+^n$  be arbitrary. Then $$z(p,q,s)_{\delta} =
e(p,q,s)\eqno (59)$$}\vskip 1pc Now $U(\u_-)$ is clearly stable under
$ad\,h$. Let $U_i(u_-) = U_i(\g)\cap U(\u_-)$. One notes that $U_i(\u_-) = 0$ if
$i$ is positive so that one has the  direct sum $$U(\u_-) =
\sum_{i=0}^{\infty} U_{-i}(\u_-)\eqno (60)$$  More explicitly, let $E =
\sum_{i=1}^{\ell}\Bbb C e_{-\alpha_i}$. Obviously $E$ generates $U(\u_-)$. Put
$E^0 = \Bbb C\,1$, and for
$j\in
\Bbb Z_+$ we define $E^{-j}$ inductively so that if $j>0$, then $E^{-j} =
E\,E^{1-j}$. The following proposition is immediate.\vskip 1pc {\bf Proposition
1.17.} {\it One has $E^{-j} = U_{-j}(\u_-)$ for any $j\in \Bbb Z_+$. In particular 
 $U_{-j}(\u_-)$ is finite dimensional. Also 
the maximal possible value of $\delta(u)$ for $0\neq u\in U(\u_-)$ is $0$.}
\vskip 1pc Recall the subspace $R\s U(\k)$. See (51). \vskip 1pc {\bf
Proposition 1.18.} {\it One has $v_{\delta}\in U(\u_-)$ for any $v\in R$.
Furthermore if $j\in \Bbb Z_+$ and $u\in U_{-j}(u_-)$ there exists $v\in R$ such
that $v_{\delta} = u$.} \vskip 1pc {\bf Proof.} Let
$0\neq v
\in R$. Then there exists a nonempty finite subset $P\s \Bbb Z_+^r$ and a
set of nonzero scalars $\{c_p\},\,p\in P$, such that $v = \sum_{p\in P}\,c_p
z(p,0,0)$. Let $m$ be the minimal value of $\delta(z(p,0,0))$ for $p\in P$ and
let $P' = \{p\in P\mid \delta(z(p,0,0)) = m\}$. Obviously $P'$ is not empty
and hence $0\neq v'\in U(u_-)$ by Remark 1.15 where $v'= \sum_{\p\in P'}
c_p\,e(p,0,0)$. But $v_{\delta} = v'$ by (54), (55) and (59). This proves the
first statement of the proposition.

For the second statement we can assume $0\neq u\in U_{-j}(u_-)$. Since any
$e(p,0,0)$ is an
 $ad\,h$-eigenvector it follows from Remark 1.15 and Proposition 1.17 that there
exists a finite subset $P_j\s \Bbb Z_+^r$ such that $\{e(p,0,0)\},\,p\in P_j,$
is a basis of $U_{-j}(u_-)$. But then there exists a nonempty subset $P_j'\s P_j$
and nonzero scalars $\{c_p\},\,p\in P_j'$, such that $u=\sum_{p\in P_j'} c_j
e(p,0,0)$. Let $v = \sum_{p\in P_j'} c_p\,z(p,0,0)$ so that $0\neq v\in R$. But
then $v_{\delta} = u$ by (55) and (59). QED\vskip 1pc {\bf Proposition 1.19.} {\it
Let
$\Gamma$ be some countable index set. Assume 
$\{b_{\gamma}\},\gamma\in \Gamma$, is a subset of $R$ with the property that
$\{(b_{\gamma})_{\delta}\},\gamma\in \Gamma$, is a basis of
$U(u_-)$. Then $\{b_{\gamma}\},\gamma\in \Gamma$, is a basis of $R$.} \vskip 1pc
{\bf Proof.} We first prove linear independence. Let $\Gamma_o$ be a finite
nonempty subset of $\Gamma$ and let $\{c_{\gamma}\},\,\gamma\in \Gamma_o$, be a
set of nonzero scalars. Put $v=\sum_{\gamma\in \Gamma_o}c_{\gamma} b_{\gamma}$.
We must show that
$v\neq 0$. Let $m$ be the minimal value of $\delta(b_{\gamma})$ for $\gamma\in
\Gamma_o$. (Clearly $m$ is finite since all $(b_{\gamma})_{\delta}$ are not
zero.) Let
$\Gamma_o' = \{\gamma\in \Gamma_o\mid \delta(\gamma) = m\}$. Obviously
$\Gamma_o'$ is not empty. Put $u = \sum_{\gamma\in \Gamma_o'}
c_{\gamma}(b_{\gamma})_{\delta}$. Then $u\neq 0$. But $v_{\gamma} = u$ by (54)
and (55). Thus $v\neq 0$. 

Let $R^s$ be the span of $\{b_{\gamma}\},\gamma\in
\Gamma$. Put $R_{\infty} = \{0\}$ and for $j\in \Bbb Z_+$ let $R_{-j}
= \{v\in R\mid \delta(v) = -j\}$. But then $$R = R_{\infty} \cup \bigcup_{j\in
\Bbb Z_+} R_{-j}\eqno (61)$$ is a disjoint union by (60) and Proposition 1.18.
Obviously $R_{\infty}\s R^s$. Assume $j\in \Bbb Z_+$ and assume inductively that
$R_i\s R^s$ for all $i> -j$ in $-\Bbb Z_+ \cup \infty$. Let $v\in R_{-j}$.
But now since all the elements of the basis
$\{(b_{\gamma})_{\delta}\},\gamma\in \Gamma$ of $U(\u_-)$ are, by definition,
eigenvectors of $ad\,h$, there exists a nonempty finite subset
$\Gamma_{j}\s\Gamma$ such that $\{(b_{\gamma})_{\delta}\},\gamma\in \Gamma_{j}$
is a basis of $U_{-j}(\u_-)$. But then there exists a nonempty subset
$\Gamma_j'$ of $\Gamma_j$ and nonzero scalars $\{c_{\gamma}\},\,\gamma\in
\Gamma_j'$ such that if $u= \sum_{\gamma\in \Gamma_j'}
c_{\gamma}(b_{\gamma})_{\delta}$, then $u = v_{\delta}$. Put $w =
\sum_{\gamma\in \Gamma_j'} c_{\gamma}b_{\gamma}$. But $w\in R^s$ and $w_{\delta} =
u$ by (55). Hence $w_{\delta} = v_{\delta}$. Put $v-w= z\in R$. But then
$\delta(z) >-j$ by (56). But then $z\in R^s$ by induction. Hence $w + z = v\in
R^s$. Thus
$R^s = R$. QED \vskip 1pc {\bf 1.9.} Recall the definition of
$z_{-\varphi}$ for $\varphi\in \Delta_+$. See the paragraph preceding Lemma 1.12.
Let
$\Pi_{s} =
\{\alpha\in
\Pi\mid
\theta(-\alpha) = \alpha\}$. One has $$\Pi_{s}\s \Pi_{\n}\eqno (62)$$ recalling (47)
and the notation of (36).
\vskip 1pc {\bf Proposition 1.20.} {\it Let
$\alpha\in
\Pi$. Then $z_{-\alpha}$ is either a semisimple element or
a nilpotent element. In addition in the latter case $z_{-\alpha}$ is
also nilpotent as an element of $\k$ (i.e. $z_{-\alpha}\in [\k,\k]$). Furthermore
the following conditions are equivalent:$$\eqalign{& (a)\,\,
\alpha\in \Pi_{s}\cr & (b)\,\, z_{-\alpha}\,\,\hbox{is semisimple}\cr & (c)\,\,
(z_{-\alpha},z_{-\alpha})\neq 0\cr}$$} \vskip 1pc {\bf Proof.} Recall (36). If
$\alpha\in \Pi_{\m}$, then $\alpha\notin \Pi_{s}$ by (62) and $z_{-\alpha} =
e_{-\alpha}$ is nilpotent. The element $z_{-\alpha}$ is also nilpotent as an
element of
$\k$ since $\g_i\s [\k,\k]$ where $i\in \{1,\ldots,\ell,\}$ is such that
$\alpha = \alpha_i$. Clearly
$z_{-\alpha}$ does not satisfy either (a), (b) or (c). Next assume $\alpha\in
\Pi_{\n}$ but $\alpha\notin \Pi_{s}$. Then $e_{-\alpha}$ commutes with
$\theta(e_{-\alpha})$ by Proposition 1.11 and hence $z_{-\alpha}= e_{-\alpha} +
\theta(e_{-\alpha})$ is nilpotent. Thus $(ad\,z_{-\alpha})^j = 0$. To
show that $z_{\alpha}$ is nilpotent as an element of $\k$, it suffices
then to prove $z_{-\alpha}\in [\k,\k]$. But the Killing form $B_{\g}$ is
nonsingular on $\k$ and certainly nonsingular on the semisimple Lie algebra
$[\k,\k]$. It follows immediately that $B_{\g}|Cent\,\k$ is nonsingular and
$Cent\,\k$ is the $B_{\g}$-orthocomplement of $[\k,\k]$ in $\k$. But if $w\in
Cent\,\k$, then
$(w,z_{-\alpha}) = tr\,ad\,z_{-\alpha}\,ad\,w$. But $(ad\,z_{-\alpha}\,ad\,w)^j
= (ad\,z_{-\alpha})^j\,(ad\,w)^j = 0$. Thus $(w,z_{-\alpha})=0$ so that
$z_{-\alpha}\in [\k,\k]$. In this present case where $\alpha\in
\Pi_{\n}$ but $\alpha\notin \Pi_{s}$,
it is obvious that (a), (b) and (c) are not satisfied. Finally assume
$\alpha\in
\Pi_{s}$. Then (c) is satisfied since
$(e_{-\alpha},e_{\alpha}) \neq 0$. But $\alpha = \alpha_i$ for some $i\in
\{1,\ldots,\ell\}$, using the notation of
\S 1.1. But then (c) implies that
$z_{-\alpha}$ is a semisimple element of the TDS $\g_i$. QED \vskip 1pc
We have not normalized our choice of root vectors. We wish to do so now to the
extent that
$\theta(e_{-\alpha_i}) = e_{\alpha_i}$ for $\alpha_i\in \Pi_{s}$. Hence for
$\alpha_i\in \Pi_{s}$ we may assume 
$$z_{-\alpha_i} = e_{-\alpha_i} + e_{\alpha_i}\eqno (63)$$ Furthermore by
replacing
$e_{-\alpha}$ by a suitable complex multiple of $e_{-\alpha_i}$, if
necessary, we can also assume that
$$(z_{-\alpha_i},z_{-\alpha_i}) = (h_i,h_i)\eqno (64)$$ using the notation of
(4). \vskip 1pc

{\bf 1.10}. Let $\lambda\in \Lambda$. For $i\in \{1,\ldots,\ell\}$ we will
now define a polynomial $q_{\lambda,i}(t)\in \Bbb C[t]$ of degree
$n_i(\lambda)+1$. See \S1.1. If $\alpha_i\notin \Pi_{s}$, put
$q_{\lambda,i}(t) = t^{n_i(\lambda) + 1}$. If $i\in \Pi_{s}$, let
$q_{\lambda,i}(t)$ be the monic polynomial of degree
$n_i(\lambda)+1$ with (multiplicity 1) roots $\{n_i(\lambda) -
2j\},\,j=0,\ldots,n_i(\lambda)$. Thus in this case
$$q_{\lambda,i}(t) = (t-n_j(\lambda))(t-n_j(\lambda)+2)\cdots (t+
n_i(\lambda))\eqno (65)$$ Now recall that $L_{\lambda}(\k)$ is the
(left ideal) annihilator in $U(\k)$ of the $U(\k)$-cyclic vector
$v_{\lambda}\in V_{\lambda}$. \vskip 1pc {\bf Proposition 1.21.} {\it
For any $i\in \{1,\ldots,\ell\}$ one has
$$q_{\lambda,i}(z_{-\alpha_i})\in L_{\lambda}(\k)\eqno (66)$$
Furthermore $$q_{\lambda,i}(z_{-\alpha_i})_{\delta} =
e_{-\alpha_i}^{n_i(\lambda) + 1}\eqno (67)$$}\vskip 1pc {\bf
Proof.} Recall (36). If $\alpha_i\in \Pi_{\m}$, then by
definition $z_{-\alpha_i} = e_{-\alpha_i}$ and hence (66) and (67)
follow from Proposition 1.1. Next assume $\alpha_i\in \Pi_{\n}$ but
$\alpha_i\notin \Pi_{s}$. Then $e_{-\alpha_i}$ commutes with
$\theta(e_{-\alpha_i})$ by Proposition 1.11 and $\theta(-\alpha_i)\in
\Delta_+$ by (47). Thus $\theta(e_{-\alpha_i})v_{\lambda}= 0$. This
proves (66) by Proposition 1.1. The equation (67) follows from
(57) and (53). Finally assume $\alpha_i\in \Pi_{s}$. The ring of
polynomial adjoint invariants on $\g_i$ is generated by the
quadratic form on $\g_i$ defined by $B_{\g}|\g_i$. Any element of
this ring therefore takes the same value on the semisimple
elements $z_{-\alpha_i}$ and $h_i$ by (64). Consequently
$z_{-\alpha_i}$ and $h_i$ are conjugate in $\g_i$. But. as one
knows, $q_{\lambda,i}(t)$ is the characteristic polynomial of
$h_i$ operating on the irreducible $U(\g_i)$-module
$U(\g_i)v_{\lambda}$. See Proposition 1.1. By conjugacy this proves
(66). The equation (67) follows from the factorization (65), with
$z_{-\alpha_i}$ replacing $t$, together with (53) and (57).
QED\vskip 1pc Recall (51). By definition $R\s
U(\k)$ is the span of
$\{z(p,0,0)\},\,p\in \Bbb Z_+^r$. On the other hand $\{e(p,0,0)\},\,p\in
\Bbb Z_+^r$ is a clearly a basis of $U(\u_-)$ as noted in Remark 1.15.
Let $I = \{1,\ldots,\ell\}$. For each $(p,i)\in \Bbb Z_+^r\times I$, let
$z(p,i) = z(p,0,0)\,q_{\lambda,i}(z_{-\alpha_i})$ and let $e(p,i) =
e(p,0,0)\,e_{-\alpha_i}^{n_i(\lambda) + 1}$. Obviously $z(p,i)\in U(\k)$
and $e(p,i)\in U(\u_-)$. Recall that $L_{\lambda}(\u_-)$ is the annihilator
of
$v_{\lambda}$ in $U(\u_-)$. See (7). \vskip 1pc {\bf Proposition 1.22.}
{\it Let
$(p,i) \in \Bbb Z_+^r\times I$. Then $z(p,i)\in L_{\lambda}(\k)$ and
$e(p,i)\in L_{\lambda}(\u_-)$.  Furthermore $$z(p,i)_{\delta} = e(p,i)\eqno
(68)$$} \vskip .5pc {\bf Proof.} The first statement follows from (66) and
(8). The second statement follows from (59), (67) and (53). QED \vskip 1pc
Now recalling (38) and (51) one has $$U(\k) = R \oplus
R\,Ker\,\omega_{\lambda}\eqno (69)$$ Note that (69) implies
$$R\,Ker\,\omega_{\lambda}= U(\k)\,Ker\,\omega_{\lambda}\eqno (70)$$ But now
obviously
$$R\,Ker\,\omega_{\lambda}\s L_{\lambda}(\k)$$ Hence if $R_{\lambda} =
R\cap L_{\lambda}(\k)$ one has $$L_{\lambda}(\k) = R_{\lambda} \oplus
R\,Ker\,\omega_{\lambda}\eqno (71)$$ For any $(p,i)\in \Bbb Z_+^r\times
I$ let (recalling Proposition 1.22) $b(p,i)\in R_{\lambda}$ and
$d(p,i)\in R\,Ker\,\omega_{\lambda}$ be the components of $z(p,i)$ with
respect to the decomposition (71) so that $$z(p,i) = b(p,i) +
d(p,i)\eqno (72)$$ \vskip .5pc {\bf Lemma 1.23.} {\it Let 
$(p,i)\in \Bbb Z_+^r\times I$. Then $$b(p,i)_{\delta} = e(p,i)\eqno
(73)$$} \vskip 1pc {\bf Proof.} If $d(p,i) = 0$ the result follows from
directly from (68). Hence we can assume $d(p,i) \neq 0$. In \S 1.7
following Lemma 1.13 we defined  $ w_i= y_i-c_i,\,i=1,\ldots,k$, where 
$\{y_i\},\,i=1,\ldots,k$, is a basis of $\hh_{\m}$ (see paragraph
preceding (37)) and
$\{c_i\},\,i=1,\ldots,k$, are arbitrary complex scalars. We now fix the
$c_i$ so that $c_i = \lambda(y_i)$. One then immediately has
$$\{z(0,q,s)\},\,(q,s)\in \Bbb Z_+^k \times \Bbb Z_+^d - (0,0)\,\,\hbox{
is a basis of $Ker\,\omega_{\lambda}$}\eqno (74)$$ But then if $(\Bbb
Z_+^n)^* = \{(p',q,s)\in \Bbb  Z_+^n\mid (q,s)\neq (0,0)\}$ one has  
$$\{z(p',q,s)\},\,(p',q,s)\in (Z_+^n)^*,\,\,\hbox{is a basis of}\,\,
R\,\,Ker\,\omega_{\lambda}\eqno (75)$$ But then there exists a finite
nonempty subset $F\s (Z_+^n)^*$ and a set $\{c_f\},\,f\in F$, of nonzero
scalars, such that $d_{p,i} = \sum_{f\in F}c_f z(f)$. Let $m$ be the
minimal value of $\delta(z(f))$ for $f\in F$ and let $F_m=\{f\in F\mid
\delta(f) = m\}$. Then $F_m$ is not empty and 
$$(d(p,i))_{\delta} = \sum_{f\in F_m} c_f e(f)\eqno (76)$$ by (54), (55)
and (59). But if $f\in F_m$ and $f = (p',q,s)$ one has $(q,s)\neq (0,0)$.
Thus $0\neq (d(p,i))_{\delta}\in U(u_-)\,Ker\,\omega_{\lambda}$. But
$(b(p,i))_{\delta}\in U(u_-)$ by Proposition 1.18. However one has
the direct sum $$U(\v) = U(\u_-) \oplus
U(u_-)\,\,Ker\,\omega_{\lambda}\eqno (77)$$ by (35). Thus
$$(b(p,i))_{\delta} + (d(p,i))_{\delta}\neq 0 $$ However
$$(b(p,i) + d(p,i))_{\delta} = z(p,i)_{\delta} = e(p,i)\in U(\u_-)$$ Thus
we cannot have $m< \delta(b(p,i))$ or even $m = \delta(b(p,i))$.
Consequently $m > \delta(b(p,i))$. Thus $(b(p,i))_{\delta}=
z(p,i)_{\delta}$. This together with (68) proves (73). QED\vskip 1pc
{\bf 1.11.} We can now state and prove our key result. \vskip 1pc
{\bf Theorem 1.24.} {\it Let $\k$ be any symmetric Lie subalgebra of a 
complex semisimple Lie algebra $\g$ (i.e. $\k$ is the set of fixed elements
in $\g$ for some involutory automorphism of $\g$.) Let $V_{\lambda}$ be a
finite dimensional irreducible $\g$-module with arbitrary highest weight
$\lambda$ and let $0\neq v_{\lambda}\in V_{\lambda}$ be a highest weight
vector. Then $V_{\lambda}$ is a cyclic
$U(\k)$-module with cyclic generator $v_{\lambda}$. Let $L_{\lambda}(\k)$ be the
annihilator (and hence left ideal) of $v_{\lambda}$ in $U(\k)$.
Then in the proceding notation $$L_{\lambda}(\k) = \sum_{i=1}^{\ell}
U(\k)\,q_{\lambda,i}(z_{-\alpha_i}) + \sum_{j=1}^k
U(\k)(y_i-\lambda(y_i)) + \sum_{\alpha\in
\Pi_{\m}}U(\k)\,e_{\alpha}\eqno (78)$$} \vskip 1pc {\bf Proof.} The
statement about cyclicity is just Proposition 1.3. For notational
convenience let $J_{\lambda} = L_{\lambda}(u_-)$ so that $J_{\lambda}$ is
the annihilator of $v_{\lambda}$ in $U(\u_-)$ (see \S 1.1) Obviously
$J_{\lambda}$ is stable under $ad\,h$ so that, recalling (60), one has
$$J_{\lambda} =
\sum_{j=0}^{\infty} (J_{\lambda})_{-j}\eqno (79)$$ where
$(J_{\lambda})_{-j}= J_{\lambda}\cap U_{-j}(\u_-)$. If $Y$ is a vector
space and $X\s Y$ is a subspace of Y we will write
$codim(X,Y)$ for the codimension of $X$ in $Y$. Thus $codim(X,Y)$ will
take values in $\Bbb Z_+\cup \{\infty\}$. Let $d_{\lambda} =
dim\,V_{\lambda}$. By (7) one has $$\sum_{j=0}^{\infty}
codim((J_{\lambda})_{-j}, U_{-j}(\u_-)) = d_{\lambda}\eqno (80)$$
Obviously all but a finite number of summands on the left side of (80)
are equal to 0. Let
$(p,i)\in \Bbb Z_+^r\times I$. Then 
$e(p,i)\in J_{\lambda}$ by Proposition 1.2 and furthermore $e(p,i)$ is an
$ad\,h$ eigenvector. Thus $$e(p,i) \in (J_{\lambda})_{-j}\eqno (81)$$ for
some $j$. Let $(Z_+^r\times I)_j = \{(p,i)\in \Bbb Z_+^r\times I\mid e(p,i)
\in (J_{\lambda})_{-j}\}$. But $\delta(e(p,i)) = \delta(e(p,0,0)) +
n_i(\lambda) +1$ by (53). Thus $(Z_+^r\times
I)_j$ is a finite set by Proposition 1.17 and Remark 1.15. But the
set $\{e(p,i)\mid (p,i)\in (Z_+^r\times I)_j\}$ spans
$(J_{\lambda})_{-j}$ by Proposition 1.2 and Remark 1.15. Let $\Gamma'_j$ be
a subset of $(Z_+^r\times I)_j$ such that $\{e(p,i)\mid (p,i)\in
\Gamma'_j\}$ is a basis of $(J_{\lambda})_{-j}$. Now let $\Gamma_j''$ be an
index set of cardinality $codim((J_{\lambda})_{-j}, U_{-j}(\u_-))$ and
choose for each $\gamma\in \Gamma_j''$ an element $e(\gamma) \in
U_{-j}(\u_-)$ such that if $\Gamma_j$ is the disjoint union of $\Gamma_j'$
and
$\Gamma_j''$, then
$\{e(\gamma)\mid \gamma\in \Gamma_j\}$ is a basis of $U_j(\u_-)$. Let
$\Gamma,\Gamma'$ and $\Gamma''$, be, respectively, the (disjoint) union
over all $j\in \Bbb Z_+$ of $\Gamma_j, \Gamma_j'$ and $\Gamma_j''$. Thus
$$\Gamma = \Gamma' \cup \Gamma''\eqno (82)$$ is disjoint union. Also, by
(80), $$card\,\,\, \Gamma'' = d_{\lambda}\eqno (83)$$ and
$$\{e(\gamma)\},\,\gamma\in \Gamma,\,\,\hbox{is a basis of
$U(\u_-)$}\eqno (84)$$ and $$\Gamma'\s \Bbb Z_+^r\times I\eqno (85)$$

Now let $L_{\lambda}(\k)^{(1)}$ be the left ideal of $U(\k)$ defined by the
right side of (78). Clearly $L_{\lambda}(\k)^{(1)}\s L_{\lambda}(\k)$ by (66).
However $codim(L_{\lambda}(\k),U(\k)) = d_{\lambda}$ by (15). Thus to prove
Theorem 22 it suffices to prove $$codim(L_{\lambda}(\k)^{(1)},U(\k))\leq
d_{\lambda}\eqno (86)$$ But now the sum of the last two sums on the right
side of (78) is clearly $U(\k)\,Ker\,\omega_{\lambda}$. Thus
$$L_{\lambda}(\k)^{(1)} = 
\sum_{i=1}^{\ell} U(\k)\,q_{\lambda,i}(z_{-\alpha_i}) +
R\,Ker\,\omega_{\lambda}\eqno (87)$$ by (70). But then if
$$L_{\lambda}(\k)^{(2)} = \sum_{i=1}^{\ell} R\,q_{\lambda,i}(z_{-\alpha_i}) +
R\,Ker\,\omega_{\lambda}$$ one has $L_{\lambda}(\k)^{(2)}\s L_{\lambda}(\k)^{(1)}$,
since of course $R\s U(\k)$, and hence to prove the theorem it suffices to
prove
$$codim(L_{\lambda}(\k)^{(2)},U(\k))\leq d_{k}\eqno (88)$$ But
$\sum_{i=1}^{\ell} R\,q_{\lambda,i}(z_{-\alpha_i})$ is the span of
$\{z(p,i)\},\,(p,i)\in \Bbb Z_+^r\times I$ by definition of $R$ (see
\S 1.8). Recall (72). One has $b(p,i)\in R$. Let $R^{(2)}$ be the span of
$\{b(p,i)\},\,(p,i)\in \Bbb Z_+^r\times I$. Since $d(p,i)\in
R\,Ker\,\omega_{\lambda}$ it then follows that $$L_{\lambda}(\k)^{(2)} = R^{(2)}
\oplus R\,Ker\,\omega_{\lambda}$$ Hence by (69) to prove the theorem it
suffices to prove $$codim\,(R^{(2)}, R) \leq d_{\lambda}\eqno (89)$$
Recall (85). Let $R^{(3)}$ be the span of $\{b(\gamma)\},\,\gamma\in
\Gamma'$. One has $R^{(3)}\s R^{(2)}$ by (85). Thus it suffices to prove
$$codim\,(R^{(3)},R) = d_{\lambda}\eqno (90)$$ But for each $\gamma\in
\Gamma''$ there exists $b(\gamma)\in R$ such that $b(\gamma)_{\delta} =
e(\gamma)$ by Proposition 1.18. Thus
$b(\gamma)\in R$ is defined for all
$\gamma\in \Gamma = \Gamma'\cup \Gamma''$. But also for any $\gamma\in
\Gamma$ one has $b(\gamma)_{\delta} = e(\delta)$ by (73). But then
$\{b(\gamma)\},\,\gamma\in \Gamma$, is a basis of $R$ by Proposition 1.19
and (84). But then if $R''$ is the span of all $b(\gamma)$ for $\gamma\in
\Gamma''$, one has that $dim\,R'' = d_{\lambda}$ and $$R = R'' \oplus
R^{(3)}$$  This proves (90). QED\vskip 1.5pc\centerline{\bf
2. The branching law for $\k$ and the precise structure of $M$} \vskip 1.5pc
{\bf 2.1.} We may regard the dual $\hh_{\m}^*$ to $\hh_{\m}$ as the subspace
of $\hh^*$ that is orthogonal to $\a$. In particular $\Delta(\m)\s
\hh_{\m}^*$. Let
$\Lambda_{\m}\s
\hh_{\m}^*$ be the set of dominant (with respect to $\b_{\m})$ integral
linear forms on
$\m$. That is, if
$\nu\in
\hh_{\m}^*$ then $\nu\in \Lambda_{\m}$ if and only if $$2(\nu,\alpha)/
(\alpha,\alpha)\in \Bbb Z_+,\,\,\,\forall \alpha\in \Pi_{\m}$$ For each
$\nu\in \Lambda_{\m}$ there exists an irreducible
representation $$\tau_{\nu}:U(\m)\to End\,Y_{\nu}\eqno (91)$$ with
highest weight $\nu$ (with respect to $\b_{\m}$). If $X$ is any
completely reducible finite dimensional $\m$-module and $\nu\in
\Lambda_{\m}$, we will denote the primary $\tau_{\nu}$-submodule of $X$ by
$X[\nu]$. The space of highest weight vectors in $X[\nu]$ is of course
$X[\nu]^{\m_+}$. Since $\m$ is reductive in $\k$ and in $\g$, any finite
dimensional completely reducible
$\k$-module or any finite dimensional $\g$-module is completely reducible
with respect to the action of $\m$. Obviously $\lambda|\m \in
\Lambda_{\m}$. The following is fairly well known. For $\lambda\in \Lambda$
let $V_{\lambda}^{\n} = \{v\in V_{\lambda}\mid \n\,v= 0\}$.
\vskip 1pc {\bf Proposition 2.1.} {\it Let $\lambda\in \Lambda$. Then
$U(\m)v_{\lambda}$ is an $\m$-irreducible
component of $V_{\lambda}[\lambda|\m]$ with highest weight vector
$v_{\lambda}$ and
$$U(\m)v_{\lambda}= U(\m_-)v_{\lambda}\eqno (92)$$ Furthermore if one 
regards $V_{\lambda}$ as an
$\a$-module, then $U(\m)v_{\lambda}$ is the
weight space for
$\a$ of weight $\lambda|\a$. In addition $$U(\m)\,v_{\lambda}=
V_{\lambda}^{\n}
\eqno (93)$$}\vskip 1pc {\bf Proof.} The
first statement (including (92)) follows from the equation 
$u\,v_{\lambda} = \omega_{\lambda}(u) v_{\lambda}$ for all $u\in
U(\b_{\m})$. Furthermore if
$V_{\lambda}^{\lambda|\a}$ is the $\a$-weight space in $V_{\lambda}$
with weight $\lambda|\a$ one has $$U(\m)v_{\lambda}\s
V_{\lambda}^{\lambda|\a}$$ since
$\m$ centralizes $\a$. But
$$U(\u_-) =
U(\m_-)\oplus U(\u_-)\n_-\eqno (94)$$ by (32) so that $$V_{\lambda} =
U(\m)v_{\lambda} +  U(\u_-)\n_-\,v_{\lambda}\eqno (95)$$ by (7) and
(92). Recalling the definition of
$w\in
\a$ in \S 1.2 one notes that the spectrum of $ad\,w$ on $U(\u_-)\n_-$ is
strictly negative. Thus the $\a$-weight space for the weight $\lambda|\a$
in the second summand on the right side of (95) vanishes. This proves the
second statement of the proposition 
and incidentally the fact that (95) is a direct sum. 

Let $V$ equal the right side of (93). Obviously $v_{\lambda}\in V$. But
since $\n$ is normalized by
$\m$ it follows that $V$ is stabilized by $\m$ and in particular one must
have that $U(\m)v_{\lambda}\s V$. But, by (26), any highest vector $v$ for
$\b_{\m}$ in
$V$ is annihilated by $\u$. Thus $v\in \Bbb C v_{\lambda}$. Hence $V$ is
$\m$-irreducible. This proves (93). QED \vskip 1pc Recall
$I = \{1,\ldots,\ell\}$. Let $I_{\m} = \{i\in I\mid \alpha_i\in \Pi_{\m}\}$.
Note that if $i\in I_{\m}$ then $h_i\in \hh_{\m}$ (see (4)) since
$h_i \in \Bbb C\,[e_{\alpha_i},e_{-\alpha_i}]$. in particular 
$\g_i\s
\m$. (See Proposition 1.1.)	\vskip 1pc {\bf Remark 2.2.} Note that if $X$
is any finite dimensional completely reducible $\m$-module and
$\lambda\in \Lambda$, then for any $i\in I_{\m}$
$$e_{-\alpha_i}^{n_i(\lambda)+1}(X[\lambda|\m]^{\m_+}) = 0\eqno
(96)$$ Indeed, as in the case of Proposition 1.1, this follows from the
representation theory of $Sl(2,\Bbb C)$ since $h_i$ reduces to the
scalar $\langle\lambda,h_i\rangle$ on $X[\lambda|\m]^{\m_+}$.
\vskip 1pc Let $I_{\n}$ be the complement of $I_{\m}$ in $I$ so that
$j\in I_{\n}$ if and only if $\alpha_j\in \Pi_{\n}$ (see (36)). We recall
that $\Pi_{s}=\{\alpha\in \Pi\mid \theta(\alpha)= -\alpha\}$ and that
$\Pi_{s}\s \Pi_{\n}$ (see (62)). Let 
$I_{s}=\{i\in I\mid \alpha_i\in \Pi_{s}\}$ and let $I_{nil}$ be the
complement of $I_{s}$ in $I_{\n}$. For any $i\in I$ and $\lambda\in \Lambda$
we defined, in
\S 1.9, a polynomial $q_{\lambda,i}(t)\in \Bbb C[t]$ of degree $n_i(\lambda)
+ 1$. The branching law, Theorem 2.3, below, will require knowledge of the
polynomials $q_{\lambda,i}$ only for $i\in I_{\n}$. Also since the
definition of 
$z_{-\alpha_i}\in \k$ for $i\in I_{\n}$ is different from its definition
when $i\in I_{\m}$, it is convenient to simplify the notation and put
$z_i = z_{-\alpha_i}$ for $i\in I_{\n}$.  Thus for $i\in I_{\n}$
$$z_i = e_{-\alpha_i} + \theta(e_{-\alpha_i})\eqno (97)$$ where, if $i\in
I_s\s I_{\n}$, then $e_{-\alpha_i}$ is normalized so that  
$\theta(e_{-\alpha_i}) = e_{\alpha_i}$ and hence $$z_i = e_{-\alpha_i} +
e_{\alpha_i}\eqno (98)$$ and, in addition, $e_{-\alpha_i}$ is 
normalized in this case so that $$(z_i,z_i)
= (h_i,h_i)\eqno (99)$$ Of course if $i\in I_{s}$ then $z_i\in \g_i$ and
in fact $$\g_i\cap \k =\Bbb C z_i\eqno (100)$$ using the notation of
Proposition 1.1. \vskip 1pc
{\bf 2.2.} If $i\in I$ and $\lambda\in \Lambda$, then we have defined a
polynomial
$q_{\lambda,i}(t)\in \Bbb C[t]$ of degree $n_i(\lambda) +1$. If $i\in
I_{s}$ then, recalling (65), $$q_{\lambda,i}(t) =
(t-n_j(\lambda))(t-n_j(\lambda)+2)\cdots (t+ n_i(\lambda))$$ 
and if
$i$ is in the complement $I_{nil}$ of $I_{s}$ in $I_{\n}$, then
$q_{\lambda,i}(t)= t^{n_i(\lambda) + 1}$. 

If $Z$ is any finite dimensional irreducible $\k$-module and $\lambda\in
\Lambda$, let
$$Z^{\lambda} = \{v\in Z[\lambda|\m]^{\m_+}\mid
q_{\lambda,i}(z_i)v = 0,\,\,\forall i\in I_{\n}\}\eqno (101)$$ Also
let $mult_{V_{\lambda}}\,(Z)$ be the multiplicity of the irreducible
representation $Z$ in $V_{\lambda}$, regarded as a $\k$-module. The following
branching law is one of the main theorems.
\vskip 1pc {\bf Theorem 2.3.} {\it Let $\k$ be any symmetric Lie
subalgebra of a complex semisimple Lie algebra $\g$ (i.e. $\k$ is the set
of fixed elements in $\g$ for some involutory automorphism of $\g$.)
Using the notation of \S 1.1, let $\lambda\in \Lambda$ and let $V_{\lambda}$
be an irreducible (necessarily finite dimensional) 
$\g$-module with highest weight $\lambda$. Let $Z$ be any
finite dimensional irreducible $\k$-module. Then in the notation of (101)
one has $$ mult_{V_{\lambda}}\,(Z) = dim\,Z^{\lambda}\eqno (102)$$} \vskip
1pc {\bf Proof.} Recalling the definition of $L_{\lambda}(\k)$ (see \S 1.2) one
has 
$$mult_{V_{\lambda}}(Z) = dim\,Z^{L_{\lambda}(\k)}\eqno (103)$$ by (18). But
now using the notation of Theorem 1.23 one obviously has
$$Z[\lambda|\m]^{\m_+} =
\{v\in Z\mid (y_i-\lambda(y_i))v = 0,\,
i=1,\,\ldots,k,\,\,\hbox{and}\,\,e_{\alpha}v=0,\,\forall \alpha\in
\Pi_{\m}\}$$ However, by (96), one automatically has
$q_{\lambda,i}(z_{-\alpha_i})v= 0$ for any $i\in I_{\m}$ and $v\in
Z[\lambda|\m]^{\m_+}$, using the notation of Theorem 1.24. But then
$$Z^{L_{\lambda}(\k)} = Z^{\lambda}\eqno(104)$$ by (78). The theorem then
follows from (103). QED \vskip 1pc {\bf 2.3.} Now let $G_{\Bbb C}$ be a
simply connected semisimple Lie group such that $\g = Lie\,G_{\Bbb C}$ and
let $G$ be the real form and subgroup of $G_{\Bbb C}$ corresponding to
$\g_o$. Thus the subgroup $K$ of
$G$ corresponding to $\k_o$ is a maximal compact subgroup of $G$. Let $A$
be the subroup of $G$ corresponding to $\a_o$ and let $M$ be the 
centralizer of $A$ in $K$. Let $\ell_o$ be the (real) dimension of
$A$ so that $\ell_o$ is the split rank of $G$. The compact group
$M$ is not connected in general and we denote the identity component of $M$
by $M_e$. Since $G_{\Bbb C}$ is simply connected, the involutory
automorphism
$\theta$ lifts to an automorphism of $G_{\Bbb C}$ which we continue to
denote by $\theta$. Of course $\theta$ is the identity on $K$. Let $A_{\Bbb
C}$ be the subgroup of $G_{\Bbb C}$ which corresponds to $\a$. One thus has 
$$A_{\Bbb C}\cong (\Bbb C^*)^{\ell_o}$$ In particular if $T$ is
the group of all elliptic elements in $A_{\Bbb C}$, then one has the 
group product structure $$A_{\Bbb C} = A\times T$$ In
particular if $F = \{a\in A_{\Bbb C}\mid a^2 = 1\}$ then $F\s T$ and
$$F\cong \Bbb Z_2^{\ell_o}\eqno (105)$$ Since $\theta=-1$ on $\a$ one has
$\theta(a) = a^{-1}$ for all $g\in A_{\Bbb C}$. The finite 2-group $F$ can
then also be characterized by $$F=\{a\in A_{\Bbb C}\mid \theta(a) = a\}\eqno
(106)$$ Let $F_M = F\cap K$. One notes that $$F_M = A_{\Bbb C}\cap K\eqno
(107)$$ Indeed any element on the right side of (107) is fixed under
$\theta$ so that (107) follows from (106). The following is well known and is
proved for completeness. Let
$$\kappa:M\to M/M_{e}\eqno (108)$$ be the quotient homomorphism \vskip 1pc
{\bf Lemma 2.4.} {\it  One has
$F_{M}\s M$ and in fact $$F_{M}\s
Cent\, M\eqno (109)$$ Furthermore $\kappa\vert F_{M}$ is surjective so
that $$M= F_{M}\,M_e\eqno (110)$$} \vskip 1pc {\bf Proof.} Clearly $F_{M}$
commutes with $A_o$, since $A_{\Bbb C}$ is commutative, so that $F_{M}\s M$. 
Furthermore if $b\in M$, then $Ad\,b$ fixes all elements in $\a_0$ and,   
by complexification, fixes all elements in $\a_{\Bbb C}$. Thus $F_M$ 
centralizes $M$ so that $F_{M}\s Cent\,M$.

Let $\g_u = \k_o +i\p_o$ so that $\g_{u}$ is a
compact real form of $\g$. Let $G_u\s G_{\Bbb C}$ be
the group corresponding to $\g_{u}$ so that
$G_u$ is a maximal compact subgroup of $G_{\Bbb C}$
and furthermore $G_{u}$ is simply connected. Of course $K\s G_{\u}$. Let
$g\in M$. Then the centralizer $G_u^g$ of $g$ in $G^u$ is a
reductive connected (since $G_u$ is simply connected)
compact subgroup of $G_u$ and $\ell = rank\,G_u^g$
since $g $ lies in some maximal torus of $G_u$. But $T\s
G_u^g$ since clearly $$Lie\,T = i\a_o\eqno (111)$$ Thus  
$T$ is contained in a maximal torus $T_u$ of $G_u^g$. Let
$\tt_{u} = Lie\,T_u$ and let $\tt_{\k} = \tt_{u}\cap \k_o$ so that
$$\tt_{\k}\s
\m_o\eqno (112)$$ But now if $z\in \tt_u$ we may write $z = x + y$ where
$x\in
\k_o$ and $y\in i\p_o$. But then $[x,i\a_o]\s ip_o$ and
$[y,i\a_o] \s k_o$. However $[z,i\a_o] = 0$. Thus we
must have $[y,i\a_o] = 0$. However $i\a_o$ is maximal
abelian in $i\p_o$. Thus $y\in i\a_o\s \tt_u$. Thus $x\in \tt_u\cap k =
\tt_{\k}$. Consequently one has $$\tt_{u} = \tt_{\k} + i\a_o\eqno (113)$$ and
hence $T_u = T_{\k}\,T$ where $T_{\k}$ is the group corresponding to
$\tt_{\k}$. But $g\in T_u$ since of course $g\in Cent\,G _u^g$ and
consequently $g$ lies in every maximal torus of $G_u^g$. Consequently we
may write $g = g_e\,a$ where $g_e\in T_{\k}$ and $a\in T$. But $T_{\k}\s
M_e$ by (112) so that $g$ is congruent to $a$ modulo $M_e$. But $a\in
T\cap K$ so that
$a\in F_M$ by (107). QED\vskip 1pc Let $j\in I$ and let $(G_{\Bbb C})_j$ be the
subgroup of $G$ which corresponds to the TDS $\g_j$ (see \S1.1). One
immediate question is whether $(G_{\Bbb C})_j$ is isomorphic to $Sl(2,\Bbb C)$ or
$PSl(2,\Bbb C)$. In fact, as noted in the next proposition, the former is
always the case. Recalling the notation of \S1.1 let $\pi_j =
\pi_{\lambda_j},\,\,V_j = V_{\lambda_j}$ and
$v_j = v_{\lambda_j}$ so that $\pi_j:\g\to End \,V_j$ is a fundamental
representation of $\g$. If $i\in
I$ and $j\neq i$ then, as a consequence of (4) and Proposition 1.1, 
$$U(\g_i)v_j = \Bbb C v_j\eqno (114)$$ whereas
$$dim\,U(\g_j)v_j = 2\eqno (115)$$
\vskip .5pc {\bf Proposition 2.5.} {\it Let $j\in I$. Then $(G_{\Bbb C})_j$
is isomorphic to 
$Sl(2,\Bbb C)$.} \vskip 1pc {\bf Proof.} It suffices to show that there
exists a representation $\pi$ of $G$ such that one of the
irreducible components of $\pi|(G_{\Bbb C})_j$ is 2-dimensional. But if $\pi =
\pi_j$ this is immediate from (115) and Proposition 1.1. QED \vskip 1pc
Of course $Cent\,Sl(2,\Bbb C)$ is isomorphic to $\Bbb Z_2$. Thus
$$Cent\,(G_{\Bbb C})_j\cong \Bbb Z_2\eqno (116)$$ for any $j\in I$ by
Proposition 2.5. Consequently for any $j\in I$  there exists a unique element
$\varepsilon_j\in (G_{\Bbb C})_j$ such that $1\neq
\varepsilon_j\in Cent\, (G_{\Bbb C})_j$. Also of course $$\varepsilon_j^2 =1\eqno
(117)$$ We are really only concerned about the case where $j\in I_{s}$.\vskip
1pc {\bf Lemma 2.6.} {\it Let
$j\in I_{s}$. Then using the notation of (4) and (98)
$$\eqalign{\varepsilon_j&= exp\,\pi\,i\,h_j\cr &=exp\,\pi\,i\,z_j\cr}\eqno
(118)$$} \vskip .5pc {\bf Proof.} Under any isomorphism of $(G_{\Bbb C})_j$ with
$Sl(2,\Bbb C)$ the element
$\varepsilon_j$ corresponds to minus the identity. But the eigenvalues of
$h_j$ on $U(\g_j)v_j$ are $\pm 1$. This proves the first line in (118). But
as observed in the proof of Proposition 1.21, the element $z_j$ is conjugate
to $h_j$ in $\g_j$. This proves the second line in (118) since
$\varepsilon_j$ is central in $(G_{\Bbb C})_j$. QED\vskip 1pc {\bf Lemma 2.7.} {\it
Let $j\in I_{s}$. Then $\varepsilon_j\in F_M$.} \vskip
1pc {\bf Proof.} It is obvious from the first line in (118) that
$\varepsilon_j\in F$. It suffices then only to prove that $\varepsilon_j\in
K$. But from the first line in (118) this would follow if one proves
$i\,z_j\in \k_o$. Oviously $i\,z_j\in \k$. Let $\theta_o$ be the conjugation
involution of $\g$ whose fixed set is $\g_o$. It suffices then to prove that
$i\,z_j$ is fixed by $\theta_o$ or equivalently $$\theta_o(z_j) = - z_j\eqno
(119)$$ But $\theta_o$ stabilizes the Cartan subalgebra $\hh = \hh_{\m} +
\a$. Hence $\theta_o(e_{-\alpha_j})$ is a root vector corresponding to some
$\beta\in \Delta$. But $\theta_o$ stabilizes both $\hh_{\m}$ and $\a$.
However $-\alpha_j$ vanishes on $\hh_{\m}$ since
$[e_{-\alpha_j},e_{\alpha_j}]\in \hh$. But $\theta_o = 1$ on $\a_o$ and
since $-\alpha_j|\a_o$ is a real linear functional one must have $\beta =
-\alpha_j$. Thus $e_{-\alpha_j}$ is an eigenvector
for $\theta_o$. But $\pm 1$ are the only eigenvalues of
$\theta_o$. Thus $\theta_o(e_{-\alpha_j}) = c\,e_{-\alpha_j}$ where $c\in
\{\pm 1\}$. But $\theta_o$ clearly commutes with $\theta$ so that
$\theta_o(e_{\alpha_j}) = c\,e_{\alpha_j}$ (see (98)). Hence $\theta_o(z_i)
= c z_i$. However one cannot have $c=1$ since otherwise $z_i\in \k_o$, and
this contradicts (99) because $B_{\g}|\k_o$ is negative definite. Hence
(119) is established. QED \vskip 1pc Let $H$ be the (Cartan) subgroup of $G$
corresponding to $\hh_o = \hh\cap \g_o$ and let $H_{\Bbb C}$ be the
(Cartan) subgroup of $G_{\Bbb C}$ corresponding to $\hh$. For $j\in
I$, let $H_j$ and $(H_j)_{\Bbb C}\,\, (\cong \Bbb C^*)$ be the subgroups of
$H_{\Bbb C}$ corresponding, respectively, to $\Bbb R\,h_j$ and $\Bbb
C\,h_j$. In particular
$(H_j)_{\Bbb C}$ is a Cartan subgroup of $(G_{\Bbb C})_j$. The equations
(114) and (115) immediately imply that $H_{\Bbb C}$ is a direct product
$$H_{\Bbb C} = (H_1)_{\Bbb C}\,\times\,\cdots \,\times (H_{\ell})_{\Bbb
C}\eqno (120)$$ Let
$\ell_{s}$ be the cardinality of
$\Pi_{s}$ (see (62)). Let $\a_{s}$ be the complex span of $\{h_j\},\,j\in I_{s}$.
Let $A_{s}$ and $(A_{s})_{\Bbb C}$ be subgroups of $G$ and $G_{\Bbb C}$,
respectively, which correspond  to $\a_{s}\cap \a_o$ and $\a_{s}$. It
is clear that
$$(A_{s})_{\Bbb C}\cong (\Bbb C^*)^{\ell_{s}}$$ In particular if
$F_{s}= F\cap (A_{s})_{\Bbb C}$ then 
$$F_{s}\cong (\Bbb Z_2)^{\ell_{s}}\eqno (121)$$

For $j\in I_{s}$ let $C_j = Cent\,(G_{\Bbb C})_j$ so that $$C_j=
\{1,\varepsilon_j\}\,\,\hbox{and}\,\,C_j\s  (H_j)_{\Bbb C}\eqno (122)$$
Clearly $C_j\s F_{s}$.
\vskip 1pc {\bf Proposition 2.8. } {\it $F_{s}$ is the group generated
by $\{C_j\},\,j\in I_{s}$. In fact $$F_{s} = \prod_{j\in \Pi_{s}}
C_j\eqno (123)$$ and (123) is a direct product.} \vskip 1pc {\bf
Proof.} Let $F_{s}'$ be the group generated by $\{C_j\},\,j\in 
I_{s}$. But
$F_{s}'$ is a direct product
 $$ F_{s}' = \prod_{j\in \Pi_{s}} C_j \eqno (124)$$ by (120).
However (124) implies that $F_{s}'$ and $F_{s}$ have the 
same order and hence
$F_{s}' = F_{s}$. QED\vskip 1pc {\bf 2.4.} Let $\varphi\in \Delta$. Then,
(see (19)), 
$\varphi|\a \neq 0$ if and only if $\varphi\in \Delta(\ss)$ An element
$\nu\in
\a^*$ is called a restricted root if there exists $\varphi\in
\Delta(\ss)$ such that $\nu = \varphi|\a$. Let $\Delta^{res}\s \a^*$ be the
set of restricted roots. By definition one has a surjection
$$p:\Delta(\ss)\to \Delta^{res}\eqno (125)$$ where $p(\varphi)=
\varphi|\a$. It is well known that $\Delta^{res}$ is a root system (see e.g.
[A]) and a system, $\Delta^{res}_+$, of positive restricted roots is
chosen (see (23) and (24)) by putting $$\Delta^{res}_+ = p
(\Delta(\n))\eqno (126)$$ For any $\gamma \in \Delta^{res}$ let
$\g(\gamma)$ be the $ad\,\a$ weight space for the weight $\gamma$. It is
obvious that
$\g(\gamma)$ is stable under $ad\,h$ and $$\g(\gamma) = \sum_{\varphi\in
\Delta(\g(\gamma))} \Bbb C e_{\varphi}\eqno (127)$$ where clearly
$$\Delta(\g(\gamma)) = p^{-1}(\gamma)\eqno (128)$$ Furthermore clearly
$[\m,\g(\gamma)]\s
\g(\gamma)]$ so that the adjoint action of $\m$ defines a representation
$$\sigma_{\gamma}:U(\m)
\to End\,\g(\gamma)\eqno (129)$$ of $\m$ on $\g(\gamma)$. The space
$\g(\gamma)$ is clearly also stable under $Ad\,M$. We will
use $\sigma_{\gamma}$ to also denote the corresponding group
representation of $M$ on $\g(\gamma)$. \vskip 1pc {\bf Proposition 2.9.} {\it
Let
$\gamma\in
\Delta^{res}$. Then
$e_{\varphi}$, for $\varphi\in p^{-1}(\gamma)$, is an $\m$-weight vector of
$\sigma_{\gamma}$ with respect to $\hh_{\m}$ and $\varphi|\hh_{\m}$ is the
corresponding weight. Furthermore all weights of $\sigma_{\m}$ have
multiplicity one.} \vskip 1pc {\bf Proof.} If $\varphi,\varphi'\in
p^{-1}(\gamma)$ and $\varphi|\hh_{\m} = \varphi'|\hh_{\m}$ then one has
$\varphi = \varphi'$ since of course $\varphi|\a = \varphi'|\a = \gamma$.
This proves the last statement of the proposition. The remainder is
obvious. QED \vskip 1pc The Killing form $B_{\g}$ is negative definite on
$\k_o$ and positive definite on $\p_o$. Thus if $B_{\theta}$ is
the symmetric bilinear form $(x,y)_{\theta}$ on $\g$ defined by putting
$$(x,y)_{\theta} = -(\theta(x),y)$$ then $B_{\theta}$ is positive definite
on $\g_o$. Furthermore the adjoint action of $K$ on $\g$ is orthogonal
with respect to $B_{\theta}$ since this action commutes with $\theta$.

Let $\gamma\in \Delta^{res}$. Since $ad\,x|\g(\gamma)$ has a real
spectrum for any $x\in \a_o$, it follows that $\g(\gamma)$ is stable under
the conjugation involution (over $\g_o$) $\theta_o$. Thus $\g(\gamma)$ is the
complexification of $\g(\gamma)_o$ where $\g(\gamma)_o = \g(\gamma)\cap
\g_o$. Let
$O(\g(\gamma)_o)$ be the orthogonal group on
$\g(\gamma)_o$ with respect to the positive definite
bilinear form $B_{\theta}|\g(\gamma)_o$. Obviously 
$$\sigma_{\gamma}(M)\,\,\hbox{stabilizes $\g(\gamma)_o$ and
$\sigma_{\gamma}(M)|\g(\gamma)_o\s O(\g(\gamma)_o)$}\eqno (130)$$ Note that
since $\theta = -1$ on $\a$ it follows that $\theta(\g(\gamma)) =
\g(-\gamma)$. But since $\theta$ commutes with $\theta_o$ one also has
$$\theta(\g(\gamma)_o) =
\g(-\gamma)_o\eqno (131)$$ Clearly $[\g(\gamma)_o,\g(-\gamma)_o]$ is
centralized by $\a_o$ so that one has $$[\g(\gamma)_o,\g(-\gamma)_o]\s
\m_o +\a_o\eqno (132)$$ We may regard $\a^*$ as the subspace of $\hh^*$
orthogonal to $\hh_{\m}$ so that $\Delta^{res}\s \hh^*$.
Recalling (1) let $w_{\gamma} = \eta^{-1}(\gamma)$ so that
$(w_{\gamma},x) = \gamma(x)$ for any $x\in \a$. Since
$B_{\g}|\a_o$ is positive definite clearly $$w_{\gamma}
\in \a_o\,\,\hbox{and}\,\,\gamma(w_{\gamma}) >0\eqno (133)$$ Let
$S(\gamma)$ be the unit sphere in $\g(\gamma)_o$ with respect to
$B_{\theta}|\g(\gamma)_o$. 
\vskip .5pc {\bf Proposition 2.10.} {\it Let
$\gamma\in \Delta^{res}$ and let $e\in S(\gamma)$. Then
$$[\theta(e),e] = w_{\gamma}\eqno (134)$$} \vskip .5pc {\bf Proof.} Let
$w= [\theta(e),e]$. Then $\theta(w) = -w$ and hence $w\in \a_o$ by (132).
On the other hand, if $x\in \a$ then since $e\in
S(\gamma)$, $$\eqalign{(w,x)&= ([\theta(e),e],x)\cr &= -(\theta(e), [x,e])\cr
&= -(\theta(e),e)\,\gamma(x)\cr &= \gamma(x)\cr}$$ This proves $w =
w_{\gamma}$. QED \vskip 1pc {\bf Remark 2.11.} Note that, as a consequence
of Proposition 2.10, if $\gamma\in \Delta^{res}$ and $e\in S(\gamma)$, then
$e,\theta(e)$ and $w_{\gamma}$ span a TDS.\vskip 1pc {\bf Lemma 2.12.} {\it
Let $\gamma\in \Delta^{res}$ and $e\in S(\gamma)$. Assume $0\neq f\in
\g(\gamma)_o$ is $B_{\theta}$-orthogonal to $e$. Then $[e,\theta(f)]$ is
a nonzero element of $\m_o$}. \vskip 1pc {\bf Proof.} Let $z=
[e,\theta(f)]$. Then $z\in \m_o + \a_o$ by (132). However if $x\in \a_o$
then $(z,x) = (\theta(f),[x,e])$. But $[x,e] = \gamma(x)\,e$. Hence
$(z,x)=0$ since $(\theta(f),e) = 0$. Thus one must have $z\in\m_o$. But
$[w_{\gamma},\theta(f)] = - \gamma(w_{\gamma})\theta(f)$. Thus if $z=0$, then
$\theta(f)$ is a highest weight vector for the adjoint action, on $\g$, of
the TDS in Remark 2.11. This contradicts a well-known result in the
representation theory of a TDS since the corresponding weight is negative by
(133). QED \vskip 1pc 
{\bf 2.5.} The following result in the split rank 1 case
($\ell_o = 1$) was proved as Theorem 2.1.7 in [K1]. The proof given here in
the general case is essentially the same as the proof in [K1]. In fact
even in the general case the statement really reduces to the split rank
1 case. Using the notation of Proposition 2.10 note that, by (130), one
has an action of
$M$ on the sphere
$S(\gamma)$ $$M\times
S(\gamma)\to S(\gamma)\eqno (135)$$ via the representation
$\sigma_{\gamma}$\vskip 1pc  {\bf Theorem 2.13.} {\it Let $\gamma\in
\Delta^{res}$. Assume
$dim\,\g(\gamma)_o > 1$. Then
$M_e$ operates transitively on $S(\gamma)$ with respect to the action
(135).}

\vskip 1pc {\bf Proof.} 
Let $e\in S(\gamma)$. Let $W$ be the
$B_{\theta}$ orthocomplement of $\Bbb R e$ in $\g(\gamma)_o$. Then  
 $[\m_o,e]\s W$ by (130). To prove the theorem it suffices
to prove $$[\m_o,e] = W\eqno (136)$$ Indeed (136) implies
that the $M_e$-orbit of $e$ is open in $S(\gamma)$. But the orbit is
certainly closed since $M_e$ is compact. Thus (136) implies
$\sigma_{\gamma}(M_e)(e) = S(\gamma)$. Assume (136) is false. Then there
exists $0\neq f\in W$ such that $(f,[y,e])_{\theta} = 0$ for any $y\in
\m_o$. But $(f,[y,e])_{\theta} = -(\theta(f),[y,e])$ and
$-(\theta(f),[y,e])= -([e,\theta(f)],y)$. But $[e,\theta(f)]$ is a
nonzero element of $\m_o$ by Lemma 2.12. Hence $-([e,\theta(f)],y)$
cannot vanish for all $y\in \m_o$ since $B_{\g}|\m_o$ is negative definite.
QED

\vskip 1pc

Let $\gamma\in \Delta^{res}$. Then a root $\psi\in
\Delta(\g(\gamma))$ will be called a lowest weight (relative to (29))
in case $\sigma_{\gamma}(\m_-)e_{\psi} = 0$. Let
$\Delta_{low}(\gamma)$ be the set of lowest weights in $\Delta(\g(\gamma))$. Since
it is understood that $\g(\gamma)$ is a $U(\m)$-module with respect to
$\sigma_{\gamma}$ we will, when convenient, suppress the symbol
$\sigma_{\gamma}$. \vskip 1pc {\bf Proposition 2.14.} {\it Let
$\gamma\in
\Delta^{res}$ and let
$\psi\in
\Delta_{low}(\gamma)$. Then $U(\m_+)e_{\psi}$ is an irreducible
component of $\sigma_{\gamma}$ and $$\g(\gamma) = \sum_{\psi\in
\Delta_{low}(\gamma)} U(\m_+)e_{\psi} \eqno (137)$$ is a unique decomposition of
$\g(\gamma)$ into irreducible components of $\sigma_{\gamma}$. In particular
$card\,\Delta_{low}(\gamma)$ is the number of such components. Furthermore all the
components define inequivalent irreducible representations of $\m$.}\vskip 1pc {\bf
Proof.} This is an immediate consequence of Proposition 2.9, and especially the
multiplicity one statement in  Proposition 2.9. QED \vskip 1pc As noted in the proof of
Lemma 2.7, the conjugation involution $\theta_o$ carries root vectors to root vectors
and consequently induces a bijection map $\theta_o:\Delta\to \Delta$. On the other hand,
clearly $\theta_o = -\theta$ on the space $i(\hh_{\m})_o + \a_o$ of hyperbolic elements in
the Cartan subalgebra $\hh$. Thus $$\theta_o = -\theta\,\,\hbox{on}\,\,\Delta\eqno
(138)$$\vskip .5pc {\bf Proposition 2.15.} {\it The conjugation involution
$\theta_o$ interchanges
$\m_+$ and $\m_-$. Also $\theta_o$ stabilizes $\n$. }\vskip 1pc {\bf Proof.} The first
statement follows from (138) since $\theta$ is clearly the identity on $\Delta(\m)$.
The second statement follows from (138) and (47). QED \vskip 1pc We will say that an
element $\gamma\in \Delta^{res}$ is irreducible in case $\sigma_{\gamma}$ is an
irreducible representation of $\m$, i.e. in case (see Proposition
2.14) $\Delta_{low}(\gamma)$ only has one element. Otherwise $\gamma$ will be called
reducible. \vskip 1pc {\bf Remark 2.16.} If $W$ is a real finite dimensional module for
a group $\Gamma$ and $W$ is real irreducible, then the commuting ring $End_{\Gamma}W$
is isomorphic to either $\Bbb R, \Bbb C$ or the quaternions $\Bbb H$. One knows that
in the first case ($\Bbb R$) the complexification
$W_{\Bbb C}$ is complex irreducible. In the second and third cases ($\Bbb C$ or $\Bbb H$)
$W_{\Bbb C}$ decomposes into a direct sum of two complex irreducible components, 
$$W_{\Bbb C} = W_{\Bbb C}' \oplus W_{\Bbb C}''\eqno (139)$$ The cases $\Bbb C$ and $\Bbb
H$ are distinguished by the fact that, in the $\Bbb H$, case the two
components of (139) are equivalent and the decomposition (139)
is not unique. In the $\Bbb C$ case the two components in (139)
are inequivalent, the decomposition is unique, and $W_{\Bbb
C}''$ is the conjugate (with respect to $W$) image of $W_{\Bbb
C}'$. Furthermore in case there exists a nonsingular symmetric $
\Gamma$-invariant bilinear form on $W_{\Bbb C}$ which is real on $W$ then of course, in
the $\Bbb R$ case, the irreducible $\Gamma$-module is self-contragredient. In the
$\Bbb C$ and $\Bbb H$ cases both components of (139) are isotropic with respect to
$B_W$ and $W_{\Bbb C}''$ has the structure of the $\Gamma$-module which
is contragredient to $W_{\Bbb C}'$. 
\vskip 1pc {\bf Theorem 2.17.} {\it Let
$\gamma\in
\Delta^{res}$. Then the real representation of 
$\sigma_{\gamma}|M_e$ on $\g(\gamma)_o$ is real irreducible so that (see Remark 2.16) if
$\gamma$ is reducible, then $\g(\gamma)$ decomposes into a sum of
only two irreducible components. That is, in any case,
$\Delta_{low}(\gamma)$ has either one or two elements according as
$\gamma$ is irreducible or reducible. Furthermore if $\gamma$ is reducible then, 
in the notation of Remark 2.16, one is always in the $\Bbb C$ case
(i.e. the $\Bbb H$ case never occurs). In particular if 
$\Delta_{low}(\gamma)= \{\psi,\psi'\}$ then $$\theta_o (U(\m_+)e_{\psi})=
U(\m_+)e_{\psi'} \eqno (140)$$ (see (137)). In addition if $\gamma$ is
irreducible, then $$\Delta(\g(\gamma))|\hh_{\m} = -
\Delta(\g(\gamma))|\hh_{\m}\eqno (141)$$ and if $\gamma$ is reducible
$$\Delta(U(\m_+)e_{\psi'})|\hh_{\m} = - \Delta(U(\m_+)e_{\psi})|\hh_{\m}
\eqno(142)$$ in the notation of (140).}\vskip 1pc {\bf Proof.} The transitivity
of $M_e$ on $S(\gamma)$, given in Theorem 2.13, obviously implies the real representation
$\sigma_{\gamma}|\g(\gamma)_o$ is real irreducible. If $\gamma$ is reducible we are in
the $\Bbb C$ case of Remark 2.16 by the inequivalence statement (last statement) in
Proposition 2.14. The existence of the $M$-invariant bilinear form
$B_{\theta}|\g(\gamma)$ puts into effect the last two statements in Remark
2.16. This yields (141) and (142) since the set of weights of a representation of a
reductive group is the negative of the set of weights of its contragredient
representation. QED \vskip 1pc {\bf Remark 2.18.} It is suggestive from (140) that
$\theta_o(\psi)$ might be $\psi'$. This, however, is not true in general. In fact by the
first statement in Proposition 2.15, $\theta_o(e_{\psi})$ is in fact the highest (and
not the lowest weight vector) in $U(\m_+)e_{\psi'}$.\vskip 1pc
 {\bf 2.6.} Let $\Pi^{res}$ be the
set of simple positive restricted roots. One knows (see [A]) that $\Pi^{res}$ is a
basis of $\a^*$ and hence we can write $$\Pi^{res} =
\{\beta_1,\ldots,\beta_{\ell_o}\}$$ Let $J = \{1,\ldots,\ell_o\}$. \vskip 1pc {\bf
Theorem 2.19.} {\it Recall the map (126). Then $$p(\Pi_{\n}) = \Pi^{res}\eqno (143)$$
In fact $$\Pi_{\n} = \bigcup_{j\in J} \Delta_{low}(\beta_j)\eqno (144)$$ where of
course (144) is a disjoint union.}
\vskip 1pc {\bf Proof.} Let $j\in J$ and let $\psi\in \Delta_{low}(\beta_j)$ so that
$\psi\in \Delta(\n)$ and $p(\psi)= \beta_j$. Assume that $\psi$ is not simple. Then
there exists $\varphi_1,\varphi_2\in \Delta_+$ such that $\psi = \varphi_1 +\varphi_2$.
However we cannot have $\varphi_1,\varphi_2\in \Delta(\n)$ since otherwise $\beta_j =
p(\varphi_1) + p(\varphi_2)$ and this contradicts the simplicity of $\beta_j$. Hence
without loss we may assume $\varphi_2\in \Delta(\m_+)$. But then obviously 
$\varphi_1\in \Delta(\n)$ and $p(\varphi_1) = \beta_j$. Also one has 
$$e_{\psi}\in \Bbb C\,\, \sigma_{\beta_j}(e_{\varphi_2})(e_{\varphi_1})\eqno (145)$$
Hence $e_{\varphi_1}$ and $e_{\psi}$ lie in the same irreducible $\m$-submodule of
$\g(\beta_j)$ by Proposition 2.9. But then (145) contradicts the minimality of
$\psi$. We have proved that the right side of (144) is contained in the left side.

 Now
let $i\in I_{\n}$ (see \S 2.1) and let $\beta = p(\alpha_i)$. Assume
$\beta\notin
\Pi^{res}$. Then we may write $\beta = \gamma +\gamma_1$ for some $\gamma,\gamma_1\in
\Delta^{res}_+$. But then there must exist $\varphi\in \Delta(\g(\gamma))$ and
$\varphi_1\in \Delta(\g(\gamma_1))$ such that $\varphi +\varphi_1 = \varphi_2$ is a root.
Otherwise, if $\g_{\gamma}$ is the TDS of Remark 2.11, one must have $[e,\g(\gamma_1)]=
0$ and
$[\theta(e), \gamma(\beta)] = 0$. But $\beta(w_{\gamma}) > \gamma_1(w_{\gamma})$ by
(133). But this contradicts either the statement that a lowest weight vector for
this TDS corresponds to a nonpositive weight or the statement that a highest weight
vector corresponds to a nonnegative weight. Hence choices can be made so
that $\varphi_2\in
\Delta(\n)$. Clearly $p(\varphi_2)= \beta$. Without loss we may assume that
$e_{\varphi_2}$ and
$e_{\alpha_i}$ lie in the same $\m$-irreducible component of $\g(\beta)$. Otherwise,
recalling (140), (138) and the last statement of Proposition 2.15, we replace $\varphi$
and $\varphi_1$ by $\theta_o(\varphi)$ and
$\theta_o(\varphi_1)$. But then $e_{\varphi_2}\in \sigma_{\beta}(U(\m))e_{\alpha_1}$ so
that $\varphi + \varphi_1$ differs from $\alpha_1$ by a linear combination of elements
in $\Pi_{\m}$. This is a contradiction of the simplicity of $\alpha_i$ since both
$\varphi$ and
$\varphi_1$ have positive coefficients of some elements in $\Pi_{\n}$ when they are
expressed as a positive integral combination of elements of $\Pi$. This proves (143) and
we can assume $\beta= \beta_j$ for some $j\in J$. We assert $\alpha_i\in
\Delta_{low}(\beta_j)$, thereby establishing (144). Indeed this is immediate from the
fact that
$\alpha_i-\alpha_{i'}$ is not a root for any $i'\in I$ and in particular for $i'\in
\Pi_{\m}$. QED \vskip 1pc 
{\bf 2.7.} Let $\cc = Cent\,\m$ so that $\cc\s \hh_{\m}$ and if
$\hh_{\m}^s = [\m,\m]\cap \hh_{\m}$, then clearly $$\hh = \hh_{\m}^s + \cc + \a\eqno
(146)$$ is a $B_{\g}$-orthogonal direct sum. Let $p_{\cc}$ and $p_{\a}$ be the
 respective projections of $\hh$ on $\cc$ and $\a$ with respect to the
decomposition (146). Let $p_{\cc + \a} = p_{\cc} + p_{\a}$. Obviously
$card\,\,\Pi_{\m} = dim\,\hh_{\m}^s$ and $\{h_i\},\,i\in I_{\m}$, is a basis of
$\hh_{\m}^s$ (see \S 2.1). This observation and of course (36) clearly proves
\vskip 1pc {\bf Proposition 2.20.} {\it One has 
$$card\,\,\Pi_{\n} = \ell_o + dim\,\cc\eqno (147)$$ and in fact (see \S 2.1)
$$\{p_{\cc + \a}(h_i)\},\,i\in I_{\n},\,\,\hbox{is a basis of}\,\,\cc + \a\eqno
(148)$$}\vskip 1pc  Let $J = J_1\cup J_2$ be the disjoint union defined so that
$\Pi^{res}_1 = \{\beta_j\mid j\in J_1\}$ is the set of irreducible restricted simple
roots and
$\Pi^{res}_2=\{\beta_j\mid j\in J_2\}$ is the set of reducible restricted simple roots.
Let
$I_{\n} = I_1\cup I_2$ be the disjoint union defined so that (recalling (144))
$p(\alpha_i)\in \Pi^{res}_1$ for $i\in I_1$ and $p(\alpha_i)\in \Pi^{res}_2$ for
$i\in I_2$. \vskip 1pc {\bf Lemma 2.21.} {\it Let $i\in I_1$. Then $p_{\cc}(h_i) = 0$.
Furthermore (see \S 1.9) $I_{s}\s I_1$.} \vskip 1pc {\bf Proof.} Let $j\in J_1$ be defined
so that
$p(\alpha_i) =\beta_j$. Then $\sigma_{\beta_j}$ is irreducible. But then
$\sigma_{\beta_j}(x)$ is a scalar operator if $x\in \cc$. However the scalar must
be 0 since it equals its negative by (141). Thus $\alpha_i\vert\cc = 0$. But then
$p_{\cc}(\hh_i) = 0$ by (3). 

Now let $i\in I_{s}$. Then $\theta_o({\alpha_i}) = \alpha_i$ by (138) and the definition
of $\Pi_{s}$ (see\S 1.9). Hence if $\beta_j = p(\alpha_j)$ then, recalling (140),
$\beta_j$ must be irreducible. QED \vskip 1pc For any $i\in I_2$ let $i'\in I_2$ be the
unique element where $i\neq i'$ and $p(\alpha_i) =p(\alpha_{i'})$. From
the inverse point of view given $j\in J_2$, let $i_j\in I_2$ be the unique element such
that $i_j<i_j\,\kern-1.8pt'$ and such that $\{\alpha_{i_j},\alpha_{i_j\,\kern-1.8pt'}\} =
p^{-1}(\beta_j)\cap \Pi_{\n}$. One has the disjoint union $$I_2 = \cup_{j\in J_2}
\{i_j,i_j\,\kern-1.8pt'\}\eqno (149)$$ \vskip .5pc {\bf Theorem 2.22.} {\it One has
$$card\,\, \Pi^{res}_2 = dim\,\cc\eqno (150) $$ Furthermore for any $j\in J_2$,
$$p_{\cc}(h_{i_j}) = -p_{\cc}(h_{i_j\,\kern-1.8pt'})\eqno (151)$$  and
$$\{p_{\cc}(h_{i_j}-h_{i_j\,\kern-1.8pt'})\},\,j\in J_2,\,\,\hbox{is a basis of}\,\,\cc\eqno (152)$$
In fact $$\{h_i\mid i\in I_{\m}\} \cup (\{h_{i_j}-h_{i_j\,\kern-1.8pt'}\},\,j\in J_2)\,\,\hbox{ is
a basis of}\,\,\hh_{\m}\eqno (153)$$}\vskip 1pc {\bf Proof.} One has $card\,\, I_{\n} =
card\,\, I_1 + card\,\, I_2$. But $card\,\,I_1 = card\,\,J_1$ and $card\,\,I_2 =
2\,card\,\,J_2$ by (149). However $card\,\,J_1 + card\,\,J_2 = \ell_o$. Thus
$card\,\,I_{\n} = \ell_o + card\,\,J_2$. But then (150) follows from (147). 

Let $j\in J_2$ so that $\g(\beta_j) = U(\m_-)e_{\alpha_{i_j}} +
U(\m_-)e_{\alpha_{i_j\,\kern-1.8pt'}}$ is the decomposition of $\g(\beta_j)$ into a direct sum
of irreducible
$\m$-modules. But if $x\in \cc$, then $\sigma_{\beta_j}(x)$ is a scalar operator in each
of these two modules. However one scalar is the negative of the other by (142).
Thus $$\alpha_{i_j}\vert \cc = -\alpha_{i_j\,\kern-1.8pt'}\vert \cc \eqno (154)$$ Next we assert
that $\alpha_{i_j}$ and $\alpha_{i_j\,\kern-1.8pt'}$ have the same length. That is
$$(\alpha_{i_j},\alpha_{i_j}) = (\alpha_{i_j\,\kern-1.8pt'},\alpha_{i_j\,\kern-1.8pt'}) \eqno (155)$$ It
suffices by (138) to show that $\theta_o(\alpha_{i_j})$ and $\alpha_{i_j\,\kern-1.8pt'}
$ have the same length. But the restriction of these two roots to $\a$ are the same,
namely $\beta_j$. It suffices then to show that the restriction of these roots to
$\hh_{\m}$ have the same length. But these restrictions are the highest and lowest
weights of the same irreducible $\m$-module by Remark 2.18.  Consequently these
restrictions have the same length since they are Weyl group conjugate for the Weyl
group of the pair $(\hh_{\m},\m)$. But then (154) and (3) proves (151). 

Since $\cc\s \hh$ it follows that $\cc$ is spanned by $\{p_{\cc}(h_i)\},\,i\in I$. But
obviously
$p_{\cc}(h_i)= 0$ for $i\in I_{\m}$. But also $p_{\cc}(h_i)= 0$ for $i\in I_1$ by
Lemma 2.21. Hence $\cc$ must be spanned by all $p_{\cc}(h_i)$ for $i\in I_2$. But, by
(151), $\cc$ must be spanned by $\{p_{\cc}(h_{i_j})\},\,j\in J_2$. But then 
$\{p_{\cc}(h_{i_j})\},\,j\in J_2$, must be a basis of $\cc$ by (150). Consequently (152)
follows from (151). 

But now if $j\in J_2$ then $\alpha_{i_j}\vert \a =\alpha_{i_j\,\kern-1.8pt'}\vert \a =\beta_j$.
Thus $p_{\a}(h_{i_j}) = p_{\a}(h_{i_j\,\kern-1.8pt'})$ by (3) and (155). This proves
$$h_{i_j}-h_{i_j\,\kern-1.8pt'}\in \hh_{\m}\eqno (156)$$ But obviously $\{h_i\mid i\in I_{\m}\}$ is a
basis of $\hh_m^s$ (see (146)). But then (153) follows from (152), (156) and the obvious
fact that $\hh_{\m} = \hh_{\m}^s + \cc$ is an orthogonal direct sum. QED\vskip 1pc 
Let
$\lambda\in \Lambda$. We will say that $\lambda$ is $\m$-trivial in
case $\lambda|\hh_{\m}=0$. \vskip 1pc {\bf Remark 2.23.} Note that Proposition 2.1
implies that $\lambda$ is $\m$-trivial if and only if $$a\,v_{\lambda} = v_{\lambda}\eqno
(157)$$ for all $a\in M_e$. \vskip 1pc {\bf Theorem 2.24.} {\it Let $\lambda\in \Lambda$.
Then, in the notation of \S1.1 and Theorem 2.22, $\lambda$ is $\m$-trivial if and only if
the following two conditions are satisfied:
$$\eqalign{&(a)\,\,n_i(\lambda) = 0,\,\,\,\forall i\in I_{\m}\cr
&(b)\,\,n_{i_j}(\lambda)= n_{i_j\,\kern-1.8pt'}(\lambda),\,\,\,\forall j\in J_2\cr}\eqno
(158)$$ where we recall that if $\ell_o$ is the split rank of $\g_o$, $J=
\{1,\ldots,\ell_o\}$, and $\Pi^{red}= \{\beta_1,\ldots,\beta_{\ell_o}\}$ is the set
of simple positive restricted roots, then $J_2= \{j\in J\mid \beta_j\,\,\hbox{is
reducible}\}$. We recall also that if $j\in J_2$, then $\{i_j,i_j\,\kern-1.8pt'\}\s I_2$
where
$I_{\n} = I_1\cup I_2$ is the disjoint union defined, after (148), and $I_{\n}$ is the
subset of $I= \{1,\ldots,\ell\}$ which indexes the set, $\Pi_{\n}$, of simple positive
roots in $\Delta(\n)$.} 
\vskip 1pc {\bf Proof.} This is immediate from (4) and (153). QED \vskip 1pc 
{\bf 2.8.} We can now deal with the full group $M$. Recall the subset $I_{s}\s I$ defined
in \S 2.1. \vskip 1pc  {\bf Lemma 2.25.} {\it Let $i\in I_{s}$. Then $\lambda_i$ is
$\m$-trivial. } \vskip 1pc {\bf Proof.} By Theorem 2.24 we have only to show that
$I_{s}$ is disjoint from $I_2$ and $I_{\m}$. But $I_{s}\s I_1$ by Lemma 2.21 and hence
$I_{s}$ is disjoint from $I_2$. But $I_1\s I_{\n}$ and hence $I_1$ is certainly disjoint
from $I_{\m}$. QED \vskip 1pc We have put
$\ell_{s} = card\,\,I_{s}$ (see (120)). Recall the 2-group
$F_{s}$ of order
$2^{\ell_{s}}$, defined in (121). One has $F_{s}\s F_M$ by Lemma
2.1 where $F_M\s M$ is defined by (107). Finally recall the quotient map $\kappa:M\to
M/M_e$ (see (108)). \vskip 1pc {\bf Lemma 2.26.} {\it The restriction $$\kappa:F_{s}\to
M/M_e\eqno (159)$$ is injective. }\vskip 1pc {\bf Proof.} Let $\varepsilon\in F_{s}\cap
M_e$. We must prove that $\varepsilon = 1$. Assume not. Then there exists a nonempty
subset $I'\s I_{s}$ such that $\varepsilon =\prod_{i\in I'}\varepsilon_i$. Let $j\in I'$.
Then using the notation of (114) and (115) one has, by (114), (115) and (118), that
$\varepsilon\,v_j = -v_j$. However $\lambda_j$ is $\m$-trivial by Lemma 2.25. Thus
$\varepsilon\, v_j = v_j$ by Remark 2.23. This is a contradiction. QED\vskip 1pc Let
$j\in I_{s}$ (see (62) and \S 2.1. Then as noted in the proof of Lemma 2.7, one has
$iz_j\in \k_o$. If $S_j\s K$ is the subgroup corresponding to $\Bbb R\,i\,z_j$, then
$S_j$ is isomorphic to the unit circle in $\Bbb C^*$ by (117) and (118).
\vskip 1pc {\bf Lemma 2.27.} {\it Assume that $Z$ is a finite dimensional $K$- (and hence
$\k$-) module. Let $j\in I_{s}$ and assume that $0\neq v\in Z$ is a
$\varepsilon_j$-eigenvector (see (117) and (118)) so that $\varepsilon_j\,v = c\,v$ where
$c\in
\{-1,1\}$. Let
$Z_j(v) = U(\Bbb C\,z_j)v$. Then $z_j$ is diagonalizable on $Z_j(v)$. Furthermore the set
of eigenvalues of $z_j|Z_j(v)$ is a set of odd integers or a set of even integers
according as $c = -1$ or $c=1$.} \vskip 1pc {\bf Proof.} Since $S_j$ is compact the
$U(\Bbb C z_j)$-module $Z_j(v)$ is completely reducible and hence $z_j|Z_j(v)$ is
diagonalizable and the eigenvalues are all integers by (117) and (118). On the other
hand, by commutativity, $\varepsilon$ must operate as the scalar operator defined by $c$
on $Z_j(v)$. The result then follows from (117) and (118). QED \vskip 1pc The following
result does not appear to be in the literature. It was, however, known to David Vogan who
points out that it is implicit in [V].
\vskip 1pc {\bf Theorem 2.28.} {\it The
 restriction $$\kappa:F_{s}\to M/M_e\eqno (160)$$ of the quotient map $\kappa$ (see
(108)) is an isomorphism. Furthermore $F_{s}\s Cent\,M$ and in fact $M$ is a
group direct product $$M =F_{s}\,\times \,M_e\eqno(161)$$ of groups so that
$$M\,\cong\,
\Bbb Z_2^{\ell_{s}}\,\times\,M_e \eqno (162)$$ (see (120) and (121)).}\vskip 1pc {\bf
Proof.} Let
$M' = F_{s}\,M$. But $$F_{s}\s F_M\eqno (163)$$ by Lemma 2.7. Hence $F_{s}\s Cent\, M$ by
(109). But then
$M' = F_{s}\,\times\, M_e$ by (159). Recalling (121) we have only to prove that $M' =
M$. But $M'$ is normal in $M$ and by Lemma 2.4. Assume $M'\neq M$. Then $M/M'\cong \Bbb
Z_2^d$ for some positive integer $d$ by Lemma (2.4). But then there exists a nontrivial
1-dimensional representation $\tau:M\to \Bbb C^*$ such that $\tau|M'$ is trivial. Since
$M\s K$ and both $M$ and $K$ are compact there exists an irreducible representation
$\sigma: K\to Aut\, V$ such that $\sigma|M$ contains $\tau$ with positive multiplicity.
Let $0\neq v\in Z$ be such that $$a\,v = \tau(a)\,v\eqno (164)$$ for all $a\in M$. But
then if $j\in I_{s}$ one has $\varepsilon_j\,v = v$. Hence if $Z_j(v)$ is defined as in
Lemma 2.27 the spectrum of $z_j|Z_j(v)$ is a set $X_j$ of even integers. Let $n_j\in 2\Bbb
Z_+$ be such that $n_j\geq |q|$ for all $q\in X_j$. But then if $\lambda\in \Lambda$ is
such that $n_j = n_j(\lambda)$, for all $j\in I_{s}$, using the notation of \S1.1, one has
$$f_{\lambda,j}(z_j)v = 0\eqno (165)$$ for all $j\in I_{s}$, using the notation of (65).
Recall (see \S 2.1) $I_{nil}$ is, by definition, the complement of $I_{s}$
in
$I_{\n}$. But then, recalling (97), if $i\in I_{nil}$, $z_i$ is a nilpotent element of
$\k$ by Proposition 1.20. We further fix
$\lambda$ by choosing
$n_i(\lambda)$, for $i\in I_{s}$ so that $$z_i^{n_i(\lambda) +1}v= 0\eqno(166)$$
Recall that
$I_2
\s I_{nil}$ since $I_{s}\s I_1$ by Lemma 2.21. For $i\in I_2$ we further choose
$n_i(\lambda)$ so that in addition to (166) one has $$n_i(\lambda) =
n_{i'}(\lambda)\eqno (167)$$ To completely determine $\lambda$ it remains only to fix
$n_i(\lambda)$ for $i\in I_{\m}$. We do so by choosing $$n_i(\lambda) = 0$$
for $i\in I_{\m}$. But then, by (4), (153), and (167) one has
$\lambda|\hh_{\m} =0$. 
But
$v\in Z[0]^{\m_+}$. Hence $v\in Z^{\lambda}$ by (101). But by (104) one has
$Z^{\lambda} = Z^{L_{\lambda}(\k)}$ so that $v\in Z^{L_{\lambda}(\k)}$. But then, by
Proposition 1.3, there exists a unique $$s\in  Hom_{K}(V_{\lambda}, Z)\eqno (168)$$
such that $s(v_{\lambda}) = v$. On the other hand, $\tau |F_M$ is nontrivial by
(110) since $\tau|M$ is nontrivial but $\tau|M_e$ is trivial. But $F_M$ is
contained in the Cartan subgroup $H_{\Bbb C}$ by (107).  Hence $$a\,v_{\lambda} =
\tau(a)\,v_{\lambda}\eqno (169)$$ for all $a\in F_M$. Now let $Y= \Bbb C\,y$ be the
1-dimensional trivial $K$-module. Then, in the notation of (101), where $Y$ replaces
$Z$, one clearly has $Y= Y[0]^{\m_+}$. But in fact $$Y=Y^{\lambda}\eqno (170)$$ That is,
$q_{\lambda,i}(z_i)v= 0$ for all $i\in I_{\n}$. This is immediate if $i\in I_{nil}$.
But if $i\in I_{s}$, then $n_i(\lambda)$ is even and hence $t$ is one of the factors
of
$q_{\lambda,i}(t)$. This proves (170). But then the argument yielding (168) implies that
there exists $s_o\in Hom_K(V_{\lambda},Y)$ such that $s_o(v_{\lambda})= y$. But
$F_M$ operates trivially on $y$ and nontrivially on $v_{\lambda}$, by (169). This is a
contradiction. Hence $M' = M$. QED\vskip 1.5pc \centerline{\bf 3. The minimal
$G$-module $V_{\lambda_{\zeta,\nu}}$ and a}\centerline
 {\bf generalization of the Cartan-Helgason Theorem}\vskip 1.5pc {\bf 3.1.} Recalling the
notation of (91), let
$\Lambda_{M_e}$ be the set of all
$\nu\in \Lambda_{\m}$ such that $\tau_{\nu}|\m_o$ exponentiates to an
(irreducible) representation of
$M_e$. Extending the use of the notation of (91), let
$\tau_{\nu}$, for
$\nu\in
\Lambda_{M_e}$, also denote the representation of $M_{e}$ which arises by
exponentiating $\tau_{\nu}|\m_o$. Let
$\widehat{F}_{s}$ be the character group of
$F_{s}$ so that
$\widehat{F}_{s}$ is a 2-group of order $2^{\ell_{s}}$. Then, by
Theorem 2.28, $\widehat{F}_{s}\,\times\, \Lambda_{M_e}$ parameterizes the set
$\widehat {M}$ of all equivalences classes of irreducible representations (the unitary
dual) of
$M$. For each
$(\zeta,\nu)\in \widehat{F}_{s}\,\times\, \Lambda_{M_e}$, let
$$\tau_{\zeta,\nu}:M_e\to Aut\, Y_{\zeta,\nu}$$ be the irreducible representation where,
as a vector space, $Y_{\zeta,\nu} = Y_{\nu}$, but with the $M$-module structure defined so
that
$$(\varepsilon,a)\,y = \zeta(\varepsilon)\,\tau_{\nu}(a)\,y\eqno (171)$$ for
$(\varepsilon,a)\in F_{s}\,\times\, M_e = M$ and $y\in Y_{\zeta,\nu}$. By abuse of
notation we will take $$\widehat {M} =\{\tau_{\zeta,\nu}\},\,(\zeta,\nu)\in
\widehat{F}_{s}\,\times\,
\Lambda_{M_e}\eqno (172)$$ For any $\lambda\in \Lambda$ let 
$\lambda|F_{s}\in
\widehat{F}_{s}$ be (uniquely) defined so that, for $j\in I_{s}$, one has  
$(\lambda|F_{s})(\varepsilon_j) = 1$ if $n_j(\lambda)$ is even and
$(\lambda|F_{s})(\varepsilon_j) = -1$ if
$n_j(\lambda)$ is odd.
\vskip 1pc {\bf Proposition 3.1.} {\it Let
$\lambda\in \Lambda$. 
 Then $\lambda|\hh_{\m} \in \Lambda_{M_e}$ so that
$(\lambda|F_{s},\lambda|\hh_m)\in 
\widehat{F}_{s}\,\times\, \Lambda_{M_e}$ and hence one has a map $$\Lambda\to
\widehat{F}_{s}\,\times\, \Lambda_{M_e},\qquad \lambda\mapsto
(\lambda|F_{s},\lambda|\hh_{\m})\eqno (173)$$ Furthermore (see Proposition
2.1) the irreducible $\m$-module $U(\m)v_{\lambda}\s V_{\lambda}$ is stable under $M$
and as an irreducible $M$-module, $U(\m)v_{\lambda}\s V_{\lambda}$, transforms according
to $\tau_{\lambda|F_{s},\lambda|\hh_{\m}}\in \widehat{M}$.}\vskip 1pc {\bf
Proof.} Since
$M\s G$ it is obvious, by Proposition 2.1, that $U(\m)v_{\lambda}$ is an irreducible
$M_e$-module with highest weight $\lambda|\hh_{\m}$. Thus $\lambda|\hh_{\m}\in
\Lambda_{M_e}$. On the other hand, if $j\in I_{s}$, then $\varepsilon_j\,v_{\lambda}
= (-1)^{n_j(\lambda)}\,v_{\lambda}$ by (4) and the first line in (118). Since $F_{s}$
centralizes $M_e$ it follows that $\varepsilon_j$ stabilizes $U(\m)v_{\lambda}$ and
reduces to the scalar operator $(\lambda|F_{s})(\varepsilon_j)$ on
$U(\m)v_{\lambda}$. QED \vskip 1pc 
{\bf 3.2.} Let $H_{\m}$ be the subgroup of $G$ corresponding
to
$(\hh_{\m})_o$ (see \S 1.2) so that $H_{\m}$ is a maximal torus of $M_e$. Let

Let $J_2$ be as in \S 2.7 and for any $j\in J_2$, let $h_{(j)} =
h_{i_j}-h_{i_j\,\kern-1.8pt'}$, using the notation of (153), so that $h_{(j)}\in \hh_{\m}$ by (153).
The element $i\, h_{(j)}$ is clearly an elliptic element of $\hh_{\m}$, by (120), and
hence
$i\, h_{(j)}\in (\hh_{\m})_o$. Let
$H_{(j)}$ be the one parameter subgroup of
$H_{\m}$ corresponding to $\Bbb R\,i\,h_{(j)}$. Let $(\hh_{(\m)})_o\s (\hh_{\m})_o$ be
the real span of
$\{h_{(j)}\},\,j\in J_2$, and let $H_{(\m)}$ be the subgroup of $H_{\m}$ corresponding
to $(\hh_{(\m)})_o$. One has $$dim\,\hh_{(m)} = dim\,\cc\eqno (174)$$ For
notational simplicity let $q = dim\,\cc$. \vskip .5pc {\bf Proposition 3.2.} {\it For
any
$j\in J_2$ the group
$H_{(j)}$ is isomorphic to the unit circle in $\Bbb C^*$. In fact if $t\in \Bbb R$, then
$exp\,t\,i\,h_{(j)} =  1$ if and only if $t\in 2\,\pi\,\Bbb Z$. Furthermore $H_{(\m)}$ is
a closed and hence compact subgroup of the torus $H_{\m}$. In fact $H_{(\m)}$ has the
product structure $$H_{(\m)} = H_{(1)}\times\cdots\times H_{(q)}\eqno (175)$$}
\vskip 1pc {\bf Proof.} It is immediate from Proposition 2.5 that if $j\in I$, then
$exp\,t\,i\,h_{j} =  1$ if and only if $t\in 2\,\pi\,\Bbb Z$. But then same statement is
true for if $j\in J_2$ and $(j)$ replaces $j$, by the product formula (120). But then,
recalling (149), the product formula (120) also yields (175). QED \vskip 1pc For any
$\lambda\in \Lambda$ let $\chi_{\lambda}: H_{\Bbb C}\to \Bbb C^*$ be the character
defined by $\lambda$ so that if $a\in H_{\Bbb C}$ then $$a\,v_{\lambda} =
\chi_{\lambda}(a)\,v_{\lambda}\eqno (176)$$ Let
$\widehat {H}_{(\m)}$ be the character group of $H_{(\m)}$. \vskip 1pc {\bf Lemma 3.3.}
{\it The group $\widehat {H}_{(\m)}$ may be parameterized by $\Bbb Z^{q}$ in such
a fashion that if
$m \in \Bbb Z^{q}$ and $m = (m_1,\ldots,m_{q})$
then $$\Phi_{m} (exp\,\sum_{j=1}^{q} t_j\,i\,h_{(j)})=
e^{i\,\sum_{j=1}^{q}t_j m_j}\eqno (176)$$ where $\Phi_{m}$ is the
character corresponding to $m$. Moreover given $m \in \Bbb Z^{q}$ there
exists $\lambda\in \Lambda$ such that $$\chi_{\lambda}|H_{(\m)} = \Phi_{m}\eqno (177)$$
In fact this is the case if and only if $$n_{i_j}(\lambda) -
n_{i_j\,\kern-1.8pt'}(\lambda) = m_j\eqno (178)$$ for all $j\in J_2$.}
\vskip 1pc {\bf Proof.} This is immediate from (4), Proposition 3.2 and the definition
above of $h_{(j)}$.  QED \vskip 1pc Recall (see above (91)) the definition of the set
$\Lambda_{\m}\s \hh_{\m}^*$ of integral dominant linear forms on $\hh_{\m}$.
That is, if $\nu\in\hh_{\m}^*$ then, by (3), $\nu \in \Lambda_{\m}$ if and only if
$$\nu(h_i)\in
\Bbb Z_+,\,\,\forall \,\,i\in I_{\m}\eqno (179)$$  The subset $\Lambda_{M_e}$ of
$\Lambda_{\m}$ was defined in \S 3.1 by the condition that $\nu \in \Lambda_{M_e}$
if the representation $\tau_{\nu}|\m_o$ (see (91)) integrates to a representation (also
denoted by $\tau_{\nu}$) of
$M_e$. Recall the basis of $\hh_{\m}$ given by (153). \vskip 1pc {\bf Theorem 3.4.} {\it
Let
$\nu\in
\hh_{\m}^*$. Then 
$\nu\in
\Lambda_{M_e}$ if and only if the following 2 conditions are satisfied:
$$\eqalign{&(1)\,\,\nu(h_i)\in \Bbb Z_+,\,\, \forall\,\,i\in I_{\m}\cr
&(2)\,\,\nu(h_{i_j}-h_{i_j\,\kern-1.8pt'})\in \Bbb Z,\,\,\forall j\in J_2\cr}\eqno (180)$$  
Furthermore if $\lambda\in \Lambda$, then $\lambda|\hh_{\m}\in \Lambda_{M_e}$ and the map
$$\Lambda\to \Lambda_{M_e},\,\,\lambda\mapsto \lambda|\hh_{\m}\eqno (181)$$ is surjective
and in fact if $\nu \in \Lambda_{M_e}$ satisfies (180), then $\lambda|\hh_{\m} = \nu$
if and only if $$\eqalign{&(1)\,\,n_i(\lambda) = \nu(h_i),\,\,\forall\,\,i\in
I_{\m}\,\,\hbox{and}\cr&(2)\,\, n_{i_j}(\lambda)-n_{i_j\,\kern-1.8pt'}(\lambda) =
\nu(h_{i_j}-h_{i_j\,\kern-1.8pt'})\in \Bbb Z,\,\,\forall\,\,j\in J_2\cr}\eqno (182)$$  In
particular if
$\lambda^a,\lambda^b\in \Lambda$, then
$\lambda^a|\hh_{\m}= \lambda^b|\hh_{\m}$ if and only if
$$\eqalign{&(1)\,\,\lambda^a(h_i)=\lambda^b(h_i),\,\,\forall\,\,i\in
I_{\m}\,\hbox{and}\cr &(2)\,\,n_{i_j}(\lambda^a)-n_{i_j\,\kern-1.8pt'}(\lambda^a)=
n_{i_j}(\lambda^b)-n_{i_j\,\kern-1.8pt'}(\lambda^b),\,\,\forall\,\,j\in J_2\cr}\eqno
(183)$$}\vskip 1pc {\bf Proof.} Assume $\nu\in
\Lambda_{M_e}$. Then $\nu\in \Lambda_{\m}$ so that (1) of (180) is satisfied. But
$\nu|(\hh_{\m})_o$ integrates to $H_{\m}$ since (recalling (91)) if $0\neq y_{\nu}$ is the
highest weight vector of $Y_{\nu}$, then one must have $exp\,x\,y_{\nu} =
e^{\nu(x)}\,y_{\nu}$ for any
$x\in (\hh_{\m})_o$. But then (2) of (180) is satisfied by Lemma 3.3 since
$(\hh_{(\m)})_o\s
(\hh_{\m})_o$. Conversely if (1) and (2) of (180) are satisfied there clearly exists
$\lambda\in \Lambda$ which satisfies (1) and (2) of (182). But then
$\pi_{\lambda}|M_e$ restricted to $U(\m)v_{\lambda}$ (see Proposition 2.1) defines an
integration of
$\tau_{\nu}|\m_o$ to $M_e$. Hence $\nu\in \Lambda_{M_e}$. 

If $\lambda\in \Lambda$ then $\lambda|\hh_{\m}$ satisfies (1) and (2) of (180) 
 so that  $\lambda\in \Lambda \in
\Lambda_{M_e}$. The argument above establishes that the map (181) is surjective and 
additionally (182) and (183) are immediate. QED \vskip 1pc 
{\bf 3.3.} We now consider the full
group $M= F_{s}\times M_e$. If $\widehat {G}$ here denotes the set of all
equivalence classes of finite dimensional irreducible representations of $G$ then
(abusing notation) we may identify $\widehat {G}$ with $\{\pi_{\lambda}\mid \lambda\in
\Lambda\}$. Thus by Proposition 3.1 one has a map $$\widehat {G} \to \widehat {M},\quad
\pi_{\lambda}\mapsto \tau_{\lambda|F_{s},\lambda|\hh_{\m}}\eqno (184)$$\vskip 1pc
{\bf Theorem 3.5.} {\it The map (185) is surjective. Furthermore if (recalling the
notation of (172)) if $(\zeta,\nu) \in \widehat {F}_{s}\times \Lambda_{M_e}$ is
arbitrary, then $\pi_{\lambda}$ is in the inverse image of $\tau_{\zeta,\nu}$ in (184) if
and only if $\lambda$ satisfies (1) and (2) of (182) and in addition
$$(-1)^{n_i(\lambda)} = \zeta(\varepsilon_i),\,\,\forall\,\,i\in I_{s}\eqno (185)$$
(recalling the notation of \S 2.1). Furthermore if $\lambda^a,\lambda^b\in \Lambda$,
then $\pi_{\lambda^a}$ and $\pi_{\lambda^b}$ have the same image under the map (184) if
and only if (1) and (2) of (183) are satisfied and in addition $$n_i(\lambda^a) \equiv
n_i(\lambda^b)\,\,\,mod\,\, 2, \,\,\forall\,\,i\in I_{s}\eqno (186)$$} \vskip 1pc
{\bf Proof.} This follows immediately from (118) and the fact that $I_{s}$ is
disjoint from both $I_2$ (see (149) and Lemma 2.21) and $I_{\m}$ (see (136)). QED\vskip
1pc Let
$\zeta_1$ be the the trivial character of $F_{s}$ so that $\tau_{\zeta_1,0}$ is the
1-dimensional trivial representation of $M$. Let $$\Lambda_{[sph]} = \{\lambda\in
\Lambda\mid \pi_{\lambda}\mapsto \tau_{\zeta_1,0},\,\hbox{with respect to}\,\, (184)\}$$
\vskip 1pc {\bf Remark 3.6.} The Cartan-Helgason theorem asserts that, for
$\lambda\in \Lambda$, $\pi_{\lambda}|K$ has the trivial representation of $K$ as a
component, i.e. $\pi_{\lambda}$ is a spherical representation of $G$, if and only if
$\lambda\in
\Lambda_{[sph]}$. The following result characterizes
$\Lambda_{[sph]}$ using the notation of \S1.1. As noted in Remark 3.8 below this
characterization yields the Cartan-Helgason theorem. It will later yield a
generalization of the Cartan-Helgason theorem wherein 
the trivial representation of
$M$ is replaced by any irreducible representation of $M$. \vskip 1pc As an immediate
corollary of Theorem 3.5 one has \vskip 1pc {\bf Theorem 3.7.} {\it Let $\lambda\in
\Lambda$. Then $\lambda\in \Lambda_{[sph]}$ if and only if the following three
conditions are satisfied: $$\eqalign{&(a)\,\,n_i(\lambda) \in 2\,\Bbb
Z_+,\,\,\forall\,\,i\in I_{s}\cr &(b)\,\,n_i(\lambda) = 0,\,\,\forall\,\,i\in
I_{\m}\cr &(c)\,\,n_{i_j}(\lambda) = n_{i_j\,\kern-1.8pt'}(\lambda),\,\,\forall\,\,j\in
J_2\cr}\eqno (187)$$}\vskip 1pc {\bf Remark 3.8.} Let $\lambda\in \Lambda$ and let $Z$
be the 1-dimensional trivial representation of $K$. To see that Theorem 3.7 yields an
algebraic proof of the Cartan-Helgason theorem it suffices by Theorem 2.3 and (103) to
show that $$Z^{\lambda} = Z\eqno (188)$$ if and only if $\lambda$ satisfies (a),(b) and
(c) of Theorem 3.7. If $\lambda$ satisfies (a), (b) and (c) of Theorem 3.7 then
$\lambda||\hh_{\m} = 0$ by Theorem 3.7 and and hence obviously $Z[0]^{\m_+} = Z$. But,
then recalling (101), 
$Z^{\lambda} = Z[0]^{\m_+}$ since $q_{\lambda,i}(z_i)$ clearly vanishes on
$Z[0]^{\m_+}$ for all $i\in I_{\n}$. The only really significant observation to be made
here is that, for $i\in I_{s}$, $n_i(\lambda)$ is even and hence $t$ is a factor of
$q_{\lambda,i}(t)$. This proves (188). Conversely if (188) is satisfied one must have
$Z[\lambda]^{\m_+}= Z$ by (101). But then obviously $\lambda|\hh_{\m}= 0$ so that
$\lambda$ satisfies (b) and (c) of Theorem 3.7. However the equality (188) implies that
$q_{\lambda,i}(z_i)$ vanishes on $Z$ for $i\in I_{s}$. Hence $t$ must be a factor
of $q_{\lambda,i}(t)$ for $i\in I_{s}$. Consequently $\lambda$ also satisfies (a).
\vskip 1pc

{\bf 3.4.} If $V^a$ and $V^b$ are two finite dimensional
$K$-modules we will say that the $K$ spectrum of $V^b$ dominates the $K$-spectrum of
$V^a$  in
case
$$\hbox{multiplicity of $Z$ in $V^a \leq\,\,$  multiplicity of
$Z$ in $V^b$}$$ for all irreducible $K$-modules $Z$. That is, if there
exists an injective map $$r\in Hom_{K}(V^a,V^b)\eqno (189)$$
Theorem 3.5 asserts that the map (184) and hence equivalently
the map (173) are surjective. We will refer to the inverse
images of elements in $\widehat {M}$, in the case of (184), or
the inverse images of elements in $\widehat {F}_{s} \times
\Lambda_{M_e}$, in the case (173), as fibers of the respective maps
(184) and (173).  \vskip .5pc {\bf Theorem 3.9.} {\it Let
$\lambda^a\in
\Lambda$ and
$\lambda'\in \Lambda_{sph}$. Put $\lambda^b = \lambda^a +
\lambda'$. Then $\lambda^a$ and $\lambda^b$ lie in the same fiber of
(173). Furthermore, recalling Proposition 1.3, one has
$$L_{\lambda^b}(\k)\s L_{\lambda^a}(\k)\eqno (190)$$ so that there
exists a unique surjective map $K$-map $$s:V_{\lambda^b}\to
V_{\lambda^a}\eqno (191)$$ such that
$s(v_{\lambda^b})= v_{\lambda^a}$. In particular the $K$-spectrum of
$V_{\lambda^b}$ dominates the $K$-spectrum of $V_{\lambda^a}$.}
\vskip 1pc {\bf Proof.} The statement about lying in the same fiber
is an immediate consequence of Theorems 3.5 and 3.7. But now (190)
follows from Theorem 1.24 since the polynomial
$q_{\lambda^a,i}(t)$ is clearly a factor of the polynomial
$q_{\lambda^b,i}(t)$ for any $i\in I$. By reductivity the kernel of
the map $s$ has a $K$-stable complement which is clearly
$K$-equivalent to $\pi_{\lambda^a}$. This proves the final statement
of the theorem. QED\vskip 1pc {\bf Remark 3.10.} The final
statement of Theorem 3.9 can be proved much more easily using the
fact that $\pi_{\lambda^b}$ is the Cartan product of
$\pi_{\lambda'}$ and $\pi_{\lambda^a}$. The spherical vector in
$V_{\lambda'}$ tensor $V_{\lambda^a}$ has a faithful
projection on the Cartan product. \vskip 1pc 
{\bf 3.5.} We will now find that
the ``fibration" (173) (or equivalently (184)) has a natural
cross-section. We recall $I= \{1,\ldots,\ell\}$. From (36) (and recalling \S 2.1) one
has the disjoint union $I = I_{\m}\cup I_{\n}$. Recalling \S 2.7 one has the disjoint
union $I_{\n} = I_1\cup I_2$. By Lemma 2.21 one has $I_{s}\s I_1$. Let
$I_{s}'$ be the complement of $I_{s}$ in $I_1$ so that $$I = I_{\m}\cup
I_{s}\cup I_{s}'\cup I_2\eqno (192)$$ is a disjoint union. We recall also
that $J = \{1,\ldots,\ell_o\}$ where $\ell_o$ is the split rank of $\g_o$ and that $J$
parameterizes the restricted simple positive roots (see \S 2.6). Recall also the
disjoint union $J= J_1\cup J_2$ (see \S 2.7) and that $J_2$ parameterizes the
``pair" decomposition (149) of $I_2$. Now for any 
$(\zeta,\nu)\in
\widehat {F}_{s}
\times
\Lambda_{M_e}$, let $\lambda_{\zeta,\nu}\in \Lambda$ be defined so that
$$\eqalign{&\,[1]\,\,\hbox{For}\,i\in I_{s},\,\,\, n_i(\lambda_{\zeta,\nu}) =
0,\,\,\hbox{if}\,\,\zeta(\varepsilon_i) = 1,\,\,\hbox{and}\,\,n_i(\lambda_{\zeta,\nu})
= 1,\,\,\hbox{if}\,\,\zeta(\varepsilon_i) = -1\cr &\,[2]\,\,\hbox{For}\,i\in
I_{\m},\,\,\,n_i(\lambda_{\zeta,\nu}) = \nu(h_i)\cr &\,[3]\,\,\hbox{For}\,j\in
J_2\,\,\,\hbox{and}\,\,\nu(h_{i_j}-h_{i_j\,\kern-1.8pt'})
\geq 0,\,\,\,n_{i_j}(\lambda_{\zeta,\nu}) =
\nu(h_{i_j}-h_{i_j\,\kern-1.8pt'})\,\,\hbox{and}\,\,\,n_{i_j\,\kern-1.8pt'}
(\lambda_{\zeta,\nu}) = 0\cr &\,[4]\,\,\hbox{For}\,j\in
J_2\,\,\,\hbox{and}\,\,\nu(h_{i_j}-h_{i_j\,\kern-1.8pt'})
\leq 0,\,\,\,n_{i_j}(\lambda_{\zeta,\nu}) = 0
\,\,\hbox{and}\,\,\,n_{i_j\,\kern-1.8pt'}(\lambda) = -
\nu(h_{i_j}-h_{i_j\,\kern-1.8pt'})\cr
 &\,[5]\,\,\hbox{For}\,i\in
I_{s}',\,\,\,n_i(\lambda_{\zeta,\nu}) = 0\cr}\eqno (193)$$ For any $(\zeta,\nu)
\in \widehat {F}_{s} \times \Lambda_{M_e}$, let ${\cal F}_{\zeta,\nu}$ be the fiber
of the map (173) over $(\zeta,\nu)$. Obviously $$\Lambda = \cup_{(\zeta,\nu)\in \widehat
{F}_{s} \times \Lambda_{M_e}}{\cal F}_{\zeta,\nu}\eqno (194)$$ is a disjoint union.
It is immediate from (193) and Theorem 3.5 that, for all $(\zeta,\nu)\in \widehat
{F}_{s} \times \Lambda_{M_e}$, 
$$\lambda_{\zeta,\nu}\in {\cal F}_{\zeta,\nu}\eqno (195)$$ We also recall that (see
\S 3.3) for the trivial pair $(\zeta_1,0)\in \widehat {F}_{s} \times
\Lambda_{M_e}$, $${\cal F}_{\zeta_1,0} = \Lambda_{[sph]}\eqno (196)$$ and that
$\Lambda_{[sph]}$ is characterized in Theorem 3.7. Note that $$\lambda_{\zeta_1,0} =
0\eqno (197)$$ Introduce a partial ordering in $\Lambda$ by defining, for
$\lambda^a,\lambda^b\in \Lambda$,  $\lambda^a <
\lambda^b$, in case $\lambda^b-\lambda^a\in \Lambda$. The following theorem
asserts that the most general fiber in (194) is just a translate of
$\Lambda_{[sph]}$. \vskip 1pc {\bf Theorem 3.11.} {\it Let $(\zeta,\nu)\in \widehat
{F}_{s} \times \Lambda_{M_e}$. Then $\lambda_{\zeta,\nu}$ is the unique
minimal element in ${\cal F}_{\zeta,\nu}$ and that furthermore $${\cal F}_{\zeta,\nu} =
\lambda_{\zeta,\nu} + \Lambda_{[sph]}\eqno (198)$$}\vskip 1pc {\bf Proof.} The proof is
an immediate observation using (193), (195), Theorem 3.5 and Theorem 3.7. QED \vskip
1pc As a corollary of Theorem 3.11 one has \vskip 1pc {\bf Theorem 3.12.} {\it Let
$(\zeta,\nu)\in \widehat {F}_{s} \times \Lambda_{M_e}$ so that $\tau_{\zeta,\nu}$
is the most general (up to equivalence) irreducible representation of $M$. Let
$\lambda\in {\cal F}_{\zeta,\mu}$ so that $\pi_{\lambda}$, up to equivalence, is the
most general finite dimensional irreducible representation of $G$ such that
$U(\m)v_{\lambda}$ transforms under $M$ according to $\tau_{\zeta,\nu}$. Then there
exists a unique $K$-map $$s:V_{\lambda}\to V_{\lambda_{\zeta,\nu}}\eqno (199) $$ such
that $s(v_{\lambda}) = v_{\lambda_{\zeta,\nu}}$. In particular, as $K$-modules,
$V_{\lambda}$ dominates $V_{\lambda_{\zeta,\nu}}$. Conversely given $\lambda\in
\Lambda$ and a $K$-map (199) such that $s(v_{\lambda}) = v_{\lambda_{\zeta,\nu}}$, then
$\lambda\in {\cal F}_{\zeta,\nu}$.}
\vskip 1pc {\bf Proof.} The first conclusion follows from Theorems 3.9 and 3.11.
The last statement follows from the fact that the condition on $s$ implies the
equivalence of $U(\m)\,v_{\lambda}$ and $U(\m)\,v_{\lambda_{\zeta,\nu}}$ as $M$-modules. 
QED\vskip 1pc {\bf Remark 3.13.} Note that Theorem 3.12 is a
generalization of the Cartan-Helgason theorem. The latter is the
case where
$\tau_{\zeta,\nu}$ is the trivial representation of $M$. Indeed if $\tau_{\zeta,\nu}$ is
the trivial representation of
$M$ then $\pi_{\lambda_{\zeta,\nu}}$ is the trivial representation of $G$, by
(197), and hence $V_{\lambda_{\zeta,\nu}}$ as a $K$-module is the trivial representation
of $K$. \vskip 1.5pc \centerline{\bf 4. The noncompact analogue of the
Borel-Weil theorem}\vskip 1.5pc {\bf 4.1.} If
$G_1$ is any Lie subgroup of $G$ let $C^{\infty}(G_1)$ be the space of all
infinitely differentiable
$\Bbb C$-valued functions on $G_1$. If $Y$ is any finite dimensional $\Bbb C$-vector
space then $C^{\infty}(G_1)\otimes Y$ naturally identifies with the space of all 
infinitely differentiable $Y$-valued functions on $G_1$. Let $\g_1= (Lie\,G_1)_{\Bbb
C}$. The space
$C^{\infty}(G_1)\otimes X$ is a
$G_1$ and a
$\g_1$-module where if $g\in G_1,\,z\in \g_1$ and $\phi\in C^{\infty}(G_1)\otimes
X$, then $g\cdot \phi$ is defined so that for any $h\in G_1$
$$(g\cdot \phi)(h) = \phi(g^{-1}\,h)\eqno (200)$$ The map $z\mapsto z\cdot \phi$ is
complex linear and for (real) $z\in Lie\,G_1$, $$(z\cdot \phi)(h) = {d\over
dt}\phi(exp\,(-t)z\,h)|_{t=0}\eqno (201)$$ 

For notational simplicity let
$D=
\widehat {F}_{s}
\times
\Lambda_{M_e}$ so that
$D$ parameterizes $\widehat M$ in that we can write $\widehat M = \{\tau_{\delta}\mid
\delta\in D\}$. See (172). Let $\delta\in D$ and let $C^{\infty}(G,
Y_{\delta})$ be the space of all $C^{\infty}$-functions $\phi$ on
$G$ with values in $Y_{\delta}$ such that for any $g\in G$ and any $m\in M$ one has 
$$\phi(g\,m) = \tau_{\delta}(m^{-1})\,\phi(g)\eqno (202)$$ The space 
$C^{\infty}(G,Y_{\delta})$ is a
$G$- and a $\g$-module where the action is defined by (200) and (201). 

The (vector) subgroup $A$ of $G$ has been
defined in \S 2.3. Let $N$ be the (unipotent) subgroup of $G$ corresponding to
$\n_o$. For any $(\delta,\xi)\in D\times \a^*$ let $H(Y_{\delta},\xi)$
be the space of all $\phi\in C^{\infty}(G,Y_{\delta})$ such that
$$\eqalign{&(a)\,\,\hbox{the span of $K\cdot \phi$ is finite dimensional}\cr
&(b)\,\,\phi(g\,a\,n) = e^{-\xi(log\,a)}\phi(g),\,\,\forall\, g\in
G, a\in A\,\,\hbox{and}\,\,n\in N\cr}\eqno (203)$$ 
where $log\,a$ is the
unique element in $\a_o$ such that $a = exp\,\log\,a$. It is simple to verify and is
well known that, with respect to (200) and (201), $H(Y_{\delta},\xi)$ is stable
under the action of $K$ and $\g$. Furthermore using the Iwasawa decomposition $G =
K\,A\,N$ one notes that any $\phi\in H(Y_{\delta},\xi)$ is determined by its restriction
to $K$. In addition as a
$K$-module the Frobenius multiplicity theorem yields the equivalence
$$H(Y_{\delta},\xi)\cong \sum_{Z\in {\widehat K}} dim\,Z[\delta]\,\,Z\eqno(204)$$
where ${\widehat K}$ is a set of representatives of the unitary dual of $K$ and,
extending the notation of \S 2.1 so as to apply to $M$ and not just
$\m$, $X[\delta]$, for $\delta\in D$, is the primary $M$-component corresponding to
$\tau_{\delta}$ in any finite dimensional $M$-module $X$. \vskip 1pc Let $\Bbb C[K]$ be 
the abstract group algebra of $K$. Then the adjoint action of $K$ on
$U(\g)$ defines a Hopf algebra structure on the smash product $${\cal A} = \Bbb
C[K]\,\#\,U(\g)$$ A 
 Harish-Chandra module is a module $H$ for ${\cal A}$ which, as a $K$-module, admits an
equivalence $$H\cong \sum_{Z\in \widehat {K}} m_Z\,Z $$ where $m_Z\in \Bbb Z_+$, for any 
$Z\in \widehat {K}$, and is such that the $\k$-action of $H$ is the differential of the
$K$-action. (This is well defined since $K\cdot v$ spans a finite dimensional subspace
of $H$ for any $v\in H$.) It is clear that, in the notation of (204), that
$H(Y_{\delta},\xi)$ is a Harish-Chandra module. Conforming to common parlance we will
also refer to a Harish-Chandra module as a $(\g,K)$-module. \vskip 1pc
{\bf Remark 4.1.}
The $(\g,K)$-module $H(Y_{\delta},\xi)$ is the Harish-Chandra module
(with respect to $K$) of the principal series representation of $G$ corresponding to
$(\delta,\xi)\in D\times \a^*$. As a $(\g,K)$-module it is a well known result of
Harish-Chandra that
$H(Y_{\delta},\xi)$ has a finite composition series and (the subquotient theorem) any
irreducible Harish-Chandra module with respect to $K$ is equivalent to a subquotient of
$H(Y_{\delta},\xi)$ for some choice of $(\delta,\xi)\in D\times \a^*$. \vskip 1pc Now
it is much more convenient for our purposes to replace $H(Y_{\delta},\xi)$ with
scalar valued functions on $G$. This may be done as follows using the Borel-Weil
theorem. Assume $G_1$ is a Lie subgroup of $G$ and let $\g_1 = (Lie\,G_1)_{\Bbb C}$. If
$z\in \g_1$ and
$f\in C^{\infty}(G_1)$ let $f\cdot z\in C^{\infty}(G_1)$ be defined so that $z\mapsto
f\cdot z$ is complex linear and if $z\in Lie\,G_1$ (i.e. $z$ is real) then 
for any $g\in G_1$, $$(f\cdot z) (g) =
{d\over dt}f(g\,exp\,t\,z)|_{t=0}\eqno (205)$$ 

Let $\delta\in D$ and let
$Y_{\delta}^*$ be the $M$- (and also $\m$-) module which is contragredient to
$Y_{\delta}$. The pairing of $y\in Y_{\delta}$ and $y^*\in Y_{\delta}^*$ is denoted by
$\langle y,y^*\rangle$. Let $0\neq y_{\delta}^*\in Y_{\delta}^*$ be a highest weight
vector with respect to the Borel subalgebra $\b_{\m}$ (see \S 1.4) of $\m$. We will
regard any $\nu\in \hh_{\m}^*$ as an element in $\b_{\m}^*$ with the property that
$\nu|\m_+ = 0$ (see (29)). Thus if
$\delta = (\zeta,\nu)\in \widehat {F}_{s}\times \Lambda_{M_e}$ (see (172)) and
$\nu^c\in \Lambda_{M_e}$ is the highest weight of $Y_{\delta}^*$, then $$
u\,y_{\delta}^* = \nu^c(u)\,y_{\delta}^*,\,\,\forall u\in \b_{\m}\eqno (206)$$ Now if
$G_1$ is any Lie subgroup of $G$ that contains $M$, let
$C^{\infty}(G_1,\delta)$ be the space of all $f \in C^{\infty}(G_1)$
such that $$\eqalign{&(a)\,\,f\cdot u = \nu^c(u)\,f,\,\,\,\forall u\in
\b_{\m}\,\,\hbox{and}\cr &(b)\,\,f(g\,\varepsilon) =
\zeta(\varepsilon)\,f(g),\,\,\forall g\in G_1,\,\varepsilon\in F_{s}\cr}\eqno
(207)$$ 
It is immediate that 
$C^{\infty}(G_1,\delta)$ is an $G_1$- and $\g_1$-module with respect to the action
defined by (200) and (201). Using an obvious extension of the Borel-Weil theorem to the
(possibly disconnected) group $M = F_{s}\times M_e$, 
 the Borel-Weil theorem asserts that the map $$Y_{\delta}\to C^{\infty}(M,\delta),\quad
y\mapsto \psi_y$$ is an equivalence of $M$ and $\m$-modules where $$\psi_y(m) = \langle
y,m\,y_{\delta}^*\rangle\eqno (208)$$\vskip .5pc
{\bf 4.2.} Let
$$Q_{\delta}:C^{\infty}(G,Y_{\delta})\to C^{\infty}(G)$$ be the map defined by
$$(Q_{\delta}(\phi))(g) = \langle \phi(g),y_{\delta}^*\rangle\eqno (209)$$ where
$g\in G$ and $\phi\in C^{\infty}(G,Y_{\delta})$. \vskip 1pc {\bf Proposition 4.2.} {\it
One has $Q_{\delta}(\phi)\in C^{\infty}(G,\delta)$ for any $\phi \in
C^{\infty}(G,Y_{\delta})$ and the map $$Q_{\delta}:C^{\infty}(G,Y_{\delta})\to
C^{\infty}(G,\delta)\eqno (210)$$ is an isomorphism of $G$- and $\g$-modules with
respect to the action defined by (200) and (201).} \vskip 1pc {\bf Proof.} Let $\phi\in
C^{\infty}(G,Y_{\delta})$ and let $g\in G,\, m\in M$. Then
$$\eqalign{Q_{\delta}(\phi)(g\,m) &= \langle \phi(g\,m),y_{\delta}^*\rangle\cr
&= \langle \tau_{\delta}(m^{-1})\phi(g),y_{\delta}^*\rangle\cr
&= \langle \phi(g),m\,y_{\delta}^*\rangle\cr}\eqno (211)$$ That is, if we put $f=
Q_{\delta}(\phi)$, then $f(g\,m) = \langle \phi(g),m\,y_{\delta}^*\rangle$. But then by
differentiation and (206) one has $(f\cdot u)(g) =
\nu^c(u)\,f(g)$ for any $u\in \b_{\m}$ where $\delta= (\zeta,\nu)$ in the
notation of (206) and (207). Furthermore (211) implies that if $\varepsilon \in
F_{s}$, then $$\eqalign{f(g\,\varepsilon) &= \langle
\phi(g),\varepsilon\,y_{\delta}^*\rangle\cr &= \zeta(\varepsilon)\,\langle
\phi(g),\,y_{\delta}^*\rangle\cr &= \zeta(\varepsilon) f(g)\cr}$$ so that $f$ satisfies
(a) and (b) of (207). This proves that $Q_{\delta}(\phi)\in C^{\infty}(G,\delta)$. 

Now let $0\neq \phi\in C^{\infty}(G,Y_{\delta})$. Then there exists $g\in G$ such that
$0\neq \phi(g)\in Y_{\delta}$. But, for $m\in M$, $\phi(g\,m) =
\tau_{\delta}(m^{-1})\phi(g)$. But then the set $\{\phi(g\,m)\mid m\in M\}$ spans
$Y_{\delta}$ since $\tau_{\delta}$ is an irreducible representation of $M$. Thus
$\phi(G)$ is not contained in the orthocomplement of $y_{\delta}^*$. Consequently 
$Q_{\delta}(\phi)\neq 0$. Thus (210) is injective. Now let
$f\in C^{\infty}(G,\delta)$. For any $g\in G$, let $\Psi_f(g)\in C^{\infty}(M)$ be
defined by putting $\Psi_f(g)(m) = f(g\,m)$ for any $m\in M$. It is immediate from the
definition of $C^{\infty}(G,\delta)$ that $\Psi_f(g)\in C^{\infty}(M,\delta)$. By (208)
there exists a unique
$y\in Y_{\delta}$ such that
$$\Psi_f(g) =
\psi_y\eqno (212)$$ Let $\phi:G\to Y_{\delta}$ be defined so that $\phi(g) = y$. Since
the function
$G\times M\to \Bbb C$ defined by $(g,m)\mapsto f(g\,m)$ is $C^{\infty}$, it is obvious
that $\phi$ is $C^{\infty}$. On the other hand, if $g\in G,\,m\in M$, then by (208)
and (212)
$$\eqalign{(\Psi_f(g))(m)&= \psi_y(m)\cr
&= \langle y,m\,y_{\delta}^*\rangle\cr
&= \langle \phi(g),m\,y_{\delta}^*\rangle\cr}\eqno (213)$$ On the other hand, if $m=
m_1\,m_2$ for
$m_1,\,m_2\in M$, then 
$$\eqalign{(\Psi_f(g))(m)&= f(g\,m)\cr
&=(\Psi_f(g\,m_1))(m_2)\cr &=\langle \phi(g\,m_1),m_2\,y_{\delta}^*\rangle\cr}$$ Hence
$$\eqalign{\langle \phi(g\,m_1),m_2\,y_{\delta}^*\rangle &= \langle
\phi(g),m_1\,m_2\,y_{\delta}^*\rangle\cr
&=\langle \tau_{\delta}(m_1^{-1}) \phi(g),m_2\,y_{\delta}^*\rangle \cr}\eqno (214)$$ But
$m_1,\,m_2\in M$ are arbitrary and, by irreducibility, $Y_{\delta}^*$ is spanned by
$m_2\,y_{\delta}^*$ over all $m_2\in M$.  Hence (213) implies $\phi(g\,m_1) =
\tau_{\delta}(m_1^{-1})\phi(g)$. Thus $\phi\in C^{\infty}(G,Y_{\delta})$. But, for
$g\in G$, $(Q_{\delta}(\phi))(g) = \langle \phi(g),y_{\delta}^*\rangle $ by definition
(see (209)) of $Q_{\delta}$. But then putting $m=1$ in (213) one has
$(Q_{\delta}(\phi))(g) =
\Psi_f(g)(1) = f(g)$. Thus $Q_{\delta}(\phi) = f$. Hence (210) is bijective. It is
immediate from (209) that (210) is a $G$- and $\g$-module map using the the action
defined by (200) and (201). QED \vskip 1pc {\bf 4.3.} Let $\xi\in \a^*$. Assume $\psi$ is
any scalar or vector space valued function on $G$. We will say that $\psi$ satisfies the
$A\,N-\xi$ condition if $$\psi(g\,a\,n) = e^{-\xi(\log\,a)}\,\psi(g)\eqno (215)$$ for
any $g\in G,\, a\in A$ and $n\in N$. Let $\delta\in D$ and let
$C^{\infty}(Y_{\delta},\xi) =
\{\phi\in C^{\infty}(G,Y_{\delta})\mid \phi\,\,\hbox{satisfies
the}\,\,A\,N-\xi\,\,\hbox{condition}\}$. Similarly, for scalar valued functions, let 
 $C^{\infty}(\delta,\xi) =\{f\in C^{\infty}(G,\delta)\mid
f\,\,\hbox{satisfies the}\,\,A\,N-\xi\,\,\hbox{condition}\}$. It is obvious that both
$C^{\infty}(Y_{\delta},\xi)$ and $C^{\infty}(\delta,\xi)$ are stable under the $G$ and
$\g$ action defined by (200) and (201). \vskip 1pc  {\bf Proposition 4.3.} {\it Let
$(\delta,\xi)\in D\times \a^*$. Then $Q_{\delta}(\phi)\in C^{\infty}(\delta,\xi)$ for
any $\phi\in C^{\infty}(Y_{\delta},\xi)$ so that the isomorphism (see
Proposition 4.2) $Q_{\delta}$ restricts to a
$G$- and
$\g$-module map $$Q_{\delta}: C^{\infty}(Y_{\delta},\xi)\to C^{\infty}(\delta,\xi)\eqno
(216)$$ with respect to the action defined by (200) and (201). Furthermore (216) is an
isomorphism of $G$- and $\g$-modules.}
\vskip 1pc {\bf Proof.} Let $\phi \in C^{\infty}(G,Y_{\delta})$ and put $f =
Q_{\delta}(\phi)$ so that $f\in C^{\infty}(G,\delta)$ by Proposition 4.2. If
$g\in G$ then by definition $$f(g) = \langle \phi(g),y_{\delta}^*\rangle\eqno (217)$$ It
is immediate from (217) that $f\in C^{\infty}(\delta,\xi)$ in case $\phi\in
C^{\infty}(Y_{\delta},\xi)$.  This establishes the map (216). Also (216) is
injective by Proposition 4.2. Now assume that
$f\in C^{\infty}(\delta,\xi)$. Let $g\in G,\,m\in M,\,a\in A$ and $n\in N$. Then since
$m$ normalizes $N$ and commutes with $A$ there exists $n'\in N$ such that
$$\eqalign{\langle
\phi(g\,a\,n),m\,y_{\delta}^*\rangle &= 
\langle \tau_{\delta}(m^{-1})\phi(g\,a\,n),m\,y_{\delta}^*\rangle\cr
&= \langle \phi(g\,a\,n\,m),y_{\delta}^*\rangle\cr
&= \langle \phi(g\,m\,a\,n'),y_{\delta}^*\rangle\cr
&= f(g\,m\,a\,n')\cr
&= e^{-\xi(\log\,a)}\,f(g\,m)\cr
&= e^{-\xi(\log\,a)}\langle \phi(g\,m),y_{\delta}^*\rangle\cr
&= e^{-\xi(\log\,a)}\langle \phi(g),m\,y_{\delta}^*\rangle\cr}$$ Thus $\phi(g\,a\,n)-
e^{-\xi(\log\,a)}\phi(g)$ is orthogonal to $m\,y_{\delta}^*$ for all $m\in M$. By the
$M$-irreducibility of $Y_{\delta}^*$ one has $\phi(g\,a\,n)= e^{-\xi(\log\,a)}\phi(g)$.
Hence $\phi\in C^{\infty}(Y_{\delta},\xi)$ so that (216) is bijective. QED\vskip 1pc
{\bf 4.4.} Let $(\delta,\xi)\in D\times \a^*$. The Harish-Chandra module
$H(Y_{\delta},\xi)$ for the principal series representation of $G$, corresponding to
$(\delta,\xi)$ has been defined in \S 4.1. In terms of the notation of
Proposition 4.3 one has
$$H(Y_{\delta},\xi) = \{\phi\in C^{\infty}(Y_{\delta},\xi)\mid K\cdot \phi\,\,\hbox{
spans a finite dimensional vector space}\} $$ Without loss we can now replace the
vector (in general) valued functions on $G$ in $H(Y_{\delta},\xi)$ by scalar valued
functions on $G$. Let $$H(\delta,\xi) = \{f \in C^{\infty}(\delta,\xi)\mid K\cdot
f\,\,\hbox{ spans a finite dimensional vector space}\} $$ In summary $H(\delta,\xi)$
is the set of all $f\in C^{\infty}(G)$ satisfying the following condtions:
$$\eqalign{&(a)\,\,f\cdot u = \nu^c(u)\,f,\,\,\forall u\in \b_{\m},\,\,\hbox{using the
notation of (207)}\cr &(b)\,\,f(g\,\varepsilon) =
\zeta(\varepsilon)\,f(g),\,\,\forall\,g\in G,\,\,\varepsilon \in
F_{s},\,\,\hbox{using the notation of (207)}\cr &(c)\,\,f\,\,\hbox{satisfies
the}\,\, A\,N-\xi\,\,\hbox{condition of (215)}\cr
&(d)\,\,K\cdot f\,\, \hbox{ spans a finite dimensional vector space}\cr}\eqno (218)$$
It is immediate that $H(\delta,\xi)$ is a $(\g,K)$-submodule of 
$C^{\infty}(\delta,\xi)$ and one readily has the following consequence of Proposition
4.3. \vskip 1pc {\bf Proposition 4.4.} {\it Let $(\delta,\xi)\in D\times \a^*$. Then
the restriction of (216) to $H(Y_{\delta},\xi)$ defines an isomorphism
$$Q_{\delta}:H(Y_{\delta},\xi) \to H(\delta,\xi)\eqno (219)$$ of $(\g,K)$-modules. We
recall that if $\phi\in H(Y_{\delta},\xi),\, f=Q_{\delta}(\phi)$ and $g\in G$, then
$$f(g) = \langle
\phi(g),y_{\delta}^*\rangle\eqno (220)$$}\vskip 1pc Henceforth we will regard 
$H(\delta,\xi)$ as the Harish-Chandra module for the principal series representation of
$G$ which corresponds to
$(\delta,\xi)$. Now it follows immediately  that the Borel subalgebra (see \S 1.1)
$\b$ of $\g$ is given by (see \S 1.4) $$\b = \b_{\m} + \a + \n\eqno (221)$$ It is
immediate by Lemma 2.6 that $\b$ is normalized by $F_{s}$ and hence the smash groduct
${\cal B} =
\Bbb C[F_{s}]\# U(\b)$ is a Hopf subalgebra of ${\cal A} = \Bbb 
C[F_{s}]\,\#\,U(\g)$. 

Let $(\delta,\xi)\in D\times \a^*$. Write $\delta = (\zeta,\nu)\in \widehat
{F}_{s}\times \Lambda_{M_e}$. Let $\chi_{\delta,\xi}:{\cal B} \to \Bbb C$ be
the unique 1-dimensional character on ${\cal B}$ defined so that for
$\varepsilon\in F_{s},\,x\in \b_{\m}, y\in \a$ and $z\in \n$, one has
$$\eqalign{&(a)\,\,
\chi_{\delta,\xi}(\varepsilon) = \zeta(\varepsilon)\cr
&(b)\,\,\chi_{\delta,\xi}(x) = \nu^{c}(x)\cr
&(c)\,\,\chi_{\delta,\xi}(y) = -\xi(y)\,\cr
&(d)\,\, \chi_{\delta,\xi}(z) =0\cr}\eqno (222)$$

 The Hopf algebra ${\cal A}$ has an antipode $s$ where for
$k\in K,\, z\in \g$, one has $s(k) = k^{-1}$ and $s(z) = -z$. Let $H$ be a
$(g,K)$-module.  Then the (full algebraic ) dual space $H^*$ to $H$ has the structure of
an ${\cal A}$-module where for $h\in H^*,\,v\in H$ and $a\in {\cal A}$ one has $(a\cdot
h)(v) = h(s(a)\cdot v)$. Using the notation of (222) let $$H^*_{\delta,\xi} = \{h\in
H^*\mid b\cdot h = \chi_{\delta,\xi}(b)\,h\,\,\forall b\in {\cal B}\}\eqno (223)$$ 
\vskip 1pc {\bf Remark 4.5.} Note that any $0\neq h\in H^*_{\delta,\xi}$ spans a
1 dimensional space which is stable under $U(\b)$ and hence, by definition, is a highest
weight vector. If
$H$ is finitely generated as a
$U(\g)$-module (e.g. if $H$ is $(\g,K)$-irreducible) it follows easily from Proposition
5.2, p. 158 in [K2] that $H$ is finitely generated as a $U(\n)$-module. It is immediate
then that
$H_{\n} = H/(\n\cdot H)$ is finite dimensional and has the structure of a ${\cal
B}$-module and also a completely reducible $M$-module. The dual space
$H_{\n}^*$ to $H_{\n}$ naturally embeds as a finite dimensional $M$- and 
${\cal B}$-submodule of $H^*$. But clearly $H^*_{\delta,\xi}\s H_{\n}^*$ for any
$(\delta,\xi)\in D\times \a$ so that $H^*_{\delta,\xi}\neq 0$ for at most a finite
number of
$(\delta,\xi)\in D\times \a^*$ and in any case $H^*_{\delta,\xi}$ is finite
dimensional. On the other hand Casselman (see [C] or [BB]) has proved that if $H\neq 0$,
then
$H_{\n}^*\neq 0$ and Lie's theorem then implies that $H^*_{\delta,\xi}\neq 0$ for
some $(\delta,\xi)\in D\times \a^*$.\vskip 1pc 

{\bf 4.5.} Now let $H$ be a nontrivial Harish-Chandra module and assume $0\neq h\in
H^*_{\delta,\xi}$ for some $(\delta,\xi)\in D\times \a^*$ and $\delta= (\zeta,\nu)\in\ 
\widehat {F}_{s}\times \Lambda_{M_e}$. For any $v\in H$ one defines a
scalar valued function $f_v$ on $G$ by putting $$f_v(g) =
e^{-\xi(log\,a)}\,h(k^{-1}\cdot v)\eqno (224)$$ where $g = k\,a\,n$ is the Iwasawa
decomposition of $g$ with respect to $G = K\,A\,N $. Since any $K\cdot v$ spans a
finite dimesnional subspace of $H$, it is immediate that $f_v\in C^{\infty}(G)$ (in fact
$f_v$ is clearly analytic on $G$). Also one has $$f_v\,\,\hbox{satisfies
the}\,\,A\,N-\xi\,\,\hbox{condition}\eqno (225)$$ See (215). This follows easily from
the definition (224) and the fact that $A$ normalizes $N$. Clearly one has a linear map
$$H
\to C^{\infty}(G),\qquad v
\mapsto f_v\eqno (226)$$ On the other hand, if
$k\in K,\,g\in G$, and $v\in H$ then $(k\cdot f_v)(g) = f_v(k^{-1}\,g)$. Hence if
$g=k_1\,a \,n$ is the Iwasawa decompositon of $g$, then $$\eqalign{(k\cdot f_v)(g)&=
f_v(k^{-1}\,k_1\,a\,n)\cr & = e^{-\xi(log\,a)}\, f_v(k^{-1}\,k_1)\cr
&= e^{-\xi(log\,a)}\, h(k_1^{-1}\,k\cdot v)\cr
&= e^{-\xi(log\,a)}\,f_{k\cdot v}(k_1)\cr
&= f_{k\cdot v}(g)\cr}\eqno (227)$$ \vskip .5pc {\bf Lemma 4.6.} {\it Let the notation
be as in (224). Then the map (226) is a map of $\Bbb C[K]\# U(\k)$-modules.
Furthermore $K\cdot f_v$ spans a finite dimensional subspace of $C^{\infty}(G)$ for any
$v\in H$.} \vskip 1pc {\bf Proof.} The statement that (226) is a map of $\Bbb
C[K]$-modules is established in (227). But this also proves the last statement of Lemma
4.6. One then has $x\cdot f_v = f_{x\cdot v}$ for any $x\in \k$ since the derivatives
involve differentiation in a finite dimensional space. QED \vskip 1pc {\bf Lemma 4.7.}
{\it Let
$\varepsilon \in F_{s},\, g\in G$ and $v\in H$. Then $$f_v(g\,\varepsilon) =
\zeta(\varepsilon)\,f_v(g)\eqno (228)$$} \vskip 1pc {\bf Proof.} Let $g = k\,a\,n$ be
the Iwasawa decomposition of $g$. Then there exists $n_1\in N$ such that
$g\,\varepsilon = k\,\varepsilon\,a\,n_1$ is the Iwasawa decomposition of
$g\,\,\varepsilon$. But since $\varepsilon\cdot h = \zeta(\varepsilon)\,h$ one has 
$$\eqalign{f_v(g\,\varepsilon) &= e^{-\xi(log\,a)} f_v(k\,\varepsilon)\cr
&= e^{-\xi(log\,a)}\,h(\varepsilon\,k^{-1}\,v)\cr
&= e^{-\xi(log\,a)} (\varepsilon\cdot h)(k^{-1}\,v)\cr
&= \zeta(\varepsilon)\,e^{-\xi(log\,a)}\,h(k^{-1}\,v)\cr
&= \zeta(\varepsilon)\,f_v(g)\cr}$$
This establishes (228). QED \vskip 1pc Topologize $H^*$ with the weak topology defined
by the pairing of $H$ and $H^*$. \vskip 1pc {\bf Lemma 4.8.} {\it Let $h_1\in H^*$ and
$x\in \k_o$ (see (9)). Then ${d\over dt}\,exp\,t\,x\cdot h_1|_{t=0}$ exists and
$${d\over dt}\,exp\,t\,x\cdot h_1|_{t=0} = x\cdot h_1\eqno (229)$$} \vskip 1pc
{\bf Proof.} Let $v\in H$. Then $(exp\,t\,x\cdot h_1)(v) = h_1(exp\,(-t)\,x\cdot v)$. But
since there exists a finite dimensional $K$ and $U(\k)$-stable subspace $H_v$ of $H$
containing $v$, one has ${d\over dt}\,exp\,-t\,x\cdot v|_{t=0} = -x\cdot v$. But this
immediately implies (229). QED\vskip 1pc {\bf Lemma 4.9.} {\it Let $u\in \b_{\m}$ and
$v\in H$. Then $$f_v\cdot u =\nu^c(u)\,f_v\eqno (230)$$ (see (205)).} \vskip 1pc {\bf
Proof.} Let
$z\in\m_o$ and $g\in G$. Let $g= k\,a\,n$ be the Iwasawa decomposition of $g$. For $t\in
\Bbb R$ there exists $n_t\in N$ such that $g\,exp\,t\,z = k\,\,exp\,t\,z\,\,a\,n_t$. Thus
$f_v(g\,exp\,t\,z) =e^{-\xi(log\,a)}f_v(k\,exp\,t\,z)$. Hence $$f_v(g\,exp\,t\,z)=
e^{-\xi(log\,a)}\,h(exp\,(-t\,z)\,\,k^{-1}\cdot v)= e^{-\xi(log\,a)}\,(exp\,t\,z\cdot
h)(k^{-1}\cdot v)$$ But then by Lemma 4.8 one has $$(f_v\cdot z)(g) =
e^{-\xi(log\,a)}\,(z\cdot h)(k^{-1}\cdot v)\eqno (231)$$ By complexification (231) is
clearly valid for any $z\in \m$. In particular, one has (231) for $z= u\in \b_{\m}$. But
$u\cdot h = \nu^c(u)\,h$. Replacing $z\cdot h$ by $\nu^c(u)\,h$ in (231) yields
$(f_v\cdot u)(g) = \nu^c(u)\,f_v(g)$. This proves (230). QED \vskip 1pc {\bf 4.6.} 
The relation of the map (226) to submodules of the principal series is established in 
\vskip 1pc {\bf Proposition 4.10.} {\it Let
$H$ be a nontrivial Harish-Chandra module. Let $(\delta,\xi)\in D\times \a^*$. Assume
there exists $0\neq h\in H^*_{\delta,\xi}$. For any $v\in H$ let $f_v$ be the
scalar valued function on $G$ defined so that if $g\in G$ and $g=k\,a\,n$ is the Iwasawa
decomposition of $g$, then $f_v(g) = e^{-\xi(log\,a)}\,h(k^{-1}\cdot v )$. In addition
$f_v\in  H(\delta,\xi)$ (see (218)) so that one has a linear map $$T_h:H\to
H(\delta,\xi)\eqno (231)$$ where $T_v(h) = f_v$.}\vskip 1pc {\bf Proof.} This follows
immediately from the characterization of $H(\delta,\xi)$ given in (218) and the results
(225), (226), Lemma 4.6, (228) and (230). QED\vskip 1pc We have already seen that (231)
is a map of $\Bbb C[K]\# U(\k)$-modules (see Lemma 4.6). We now wish to show that it is a
map of $(\g,K)$-modules. We recall that in ${\cal A}$ one has the relation $$k\,u =
Ad\,k(u)\,k\eqno (232)$$ The existence of an element $h$ satisfying the condition of
Proposition 4.10, when $H$ is irreducible, is due to Casselman (see [C] or [BB]). Its
significance is the following submodule theorem, due to Casselman. \vskip 1pc {\bf Theorem
4.11.}{\it Let the notation be as in Proposition 4.10. Then
$T_h$ is a map of
$(\g,K)$-modules.}
\vskip 1pc {\bf Proof.} Let
$v\in H$ so that $f_v = T_h(v)$. Let $u\in \g$. Then we must show $$u\cdot f_v =
f_{u\cdot v}\eqno (233)$$ By complex linearity we can assume $u\in \g_o$. Let $g\in G$.
Then $$(u\cdot f_v)(g) = {d\over dt}\,f_v(exp\,(-t)u\,g)|_{t=0}\eqno (234)$$ Let $g =
k\,a\,n$ be the Iwasawa decomposition of $g$. Then $$f_v(exp\,(-t)u\,g) =
e^{-\xi(log\,a)}\,f_v(exp\,(-t)u\,k)$$ Let $w = Ad\,k^{-1}(u)$ so that, by Lemma 4.6,
$$f_v(exp\,(-t)u\,g) = e^{-\xi(log\,a)}\,f_{k^{-1}\cdot v}(exp\,(-t)w)\eqno (235)$$ Let
$exp\,(-t)w = k_t\,a_t\,n_t$ be the Iwasawa decomposition of $exp\,(-t)w$. Hence if we
write $w = x + y +z$ where $x\in \k_o,\,y\in \a_o$ and $z\in \n_o$, then the tangent
vectors to the curves $k_t,a_t$ and $n_t$ at $t=0$ are respecrtively $-x,\,-y$ and
$-z$. But $$f_{k^{-1}\cdot v}(exp\,(-t)w) = e^{-\xi(log\,a_t)}\,f_{k^{-1}\cdot
v}(k_t) = e^{-\xi(log\,a_t)} f_{k_t^{-1}\,k^{-1}\cdot v}(1)$$ However recall
that since $h\in H^*_{\delta,\xi}$, one has $z\cdot h = 0$ and $y \cdot h = -\xi(y)\,h$. 
Consequently differentiating (235) at
$t=0$ one has
$$(u\cdot f_v)(g) = e^{-\xi(log\,a)} (
\xi(y)\,f_{k^{-1}\cdot v}(1) + f_{x\,k^{-1}\cdot v}(1))$$  Thus $$\eqalign{
(u\cdot f_v)(g)&= e^{-\xi(log\,a)}(\xi(y) \,h(k^{-1}\cdot v) + h(x\,k^{-1}\cdot v))\cr
&= e^{-\xi(log\,a)}(k (\xi(y)-x) \cdot h)(v)\cr
&= e^{-\xi(log\,a)}(-k\,w\cdot h)(v)\cr
&= e^{-\xi(log\,a)}(-u\,k\cdot h)(v)\cr
&= e^{-\xi(log\,a)} h(k^{-1}u\cdot v)\cr
&= e^{-\xi(log\,a)}\,f_{u\cdot v}(k)\cr
&= f_{u\cdot v}(g)\cr}$$ This proves (233). QED\vskip 1pc Given the $T_h$ map (231) of 
$(\g,K)$-modules, one recovers the linear functional $h$ by composing $T_h$ with the Dirac
measure at the identity of $G$.\vskip 1pc {\bf Proposition 4.12.} {\it Let the notation
be as in Theorem 4.11. Then for any $v\in H$ one has $$h(v) = (T_h(v))(1)\eqno
(236)$$} \vskip 1pc {\bf Proof.} Let $v\in H$ so that in previous notation $T_h(v) = f_v$.
But by definition (see (224)) one has $f_v(1) = h(v)$. QED \vskip 1pc 4.7. We now deal
with the converse of Theorem 4.11 and Proposition 4.12. Let $(\delta,\xi)\in D\times
\a^*$ and let $0\neq H$ be a Harish-Chandra module and assume that $$T:H\to
H(\delta,\xi)$$ is a nontrivial $(\g,K)$-module map. Let $h\in H^*$ be defined by
putting $h(v) = f_v(1)$ where $v\in H$ and $f_v = T(v)$. If $g\in G$ and $g = k\,a\,n$
is the Iwasawa decomposition of $g$ then, since $f_v\in H(\delta,\xi)$, one has $$f_v(g) =
e^{-\xi(log\,a)}f_v(k)\eqno (237)$$ But since $T\neq 0$ there exists $v\in H$ such that
$f_v\neq 0$. Thus there exists $g\in G$ such that $f_v(g)\neq 0$. Using the notation of
(237) one has $f_v(k)\neq 0$. But $f_v(k) = (k^{-1}\cdot f_v)(1) = f_{k^{-1}\cdot
v}(1)$. Thus $h(k^{-1}\cdot v) \neq 0$ and hence $h\neq 0$. 

Now let $w\in \g_o$ and let $v\in H$ and let the notation be as in the proof of Theorem
4.11. Then (237) implies $$f_v(exp\,(-t)w) =
e^{-\xi(log\,a_t)}f_v(k_t)\eqno (238)$$ Since $w\cdot f_v= f_{w\cdot v}$ one has, by
differentiating (238) at $t=0$, $f_{w\cdot v}(1) = \xi(y)\,f_v(1) + f_{x\cdot v}(1)$.
That is $h(w\cdot v) = h(\xi(y)v + x\cdot v)$. Since $v\in H$ is arbitrary, this
implies $$-w\cdot h = \xi(y)h - x\cdot h\eqno (239)$$ By complexification one has (239)
for any $w\in \g$ and $w = x + y + z$ where $x\in\k,\,y\in \a$ and $z\in\n$. Now
choosing $w=-u\in \a + \n$, it follows from (239) and recalling (222) one has $$u\cdot
h = \chi_{\delta,\xi}(u)\,h\eqno (240)$$ for any $u\in \a + \n$. But now, recalling
(238), note that $(w\cdot f_v)(1) = (f_v\cdot (-w))(1)$. But then (239) implies that, 
for
$w= x\in \k$, $(f_v\cdot (-x))(1) = (-x\cdot h)(v)$. But
$f_v\cdot u =
\nu^c(u)\, f_v$ for
$u\in \b_{\m}$ and $\delta = (\zeta,\nu)\in {\widehat F}_{s}\times \Lambda_{M_e}$.
Thus for $-x= u\in \b_{\m}$ one has $u\cdot h = \nu(u)\,h$. Hence (240) is established
for all $u\in \b$ (see (221)). Finally let $\varepsilon \in F_{s}$. Then, by
(218), one has $f_v(\varepsilon)= \zeta(\varepsilon)\,f_v(1)$. But $f_v(\varepsilon)=
\varepsilon\cdot f_v(1)$ and $\varepsilon\cdot f_v= f_{\varepsilon\cdot v}$. Thus
$h(\varepsilon\,\cdot v) = \zeta(\varepsilon)\,h(v)$. But this implies $\varepsilon\cdot
h = \zeta(\varepsilon)\,h$. Hence, recalling (218), $h\in H^*_{\delta,\xi}$. This is
the essential step in proving the following converse to Theorem 4.11 and Proposition
4.12.
\vskip 1pc {\bf Theorem 4.13.} {\it Let $0\neq H$ be a Harish-Chandra module. Let
$(\delta,\xi)\in D\times \a^*$ and assume
$T:H\to H(\delta,\xi)$ is a nontrivial map of $(\g,K)$-modules. Let $h\in H^*$ be
defined by putting $h(v) = T(v)(1)$. Then $0\neq h\in H^*_{\delta,\xi}$. Furthermore
(recalling Theorem 4.11) $$T = T_h\eqno (241)$$}\vskip 1pc {\bf Proof.} We established
above that $0\neq h\in H^*_{\delta,\xi}$. For any $v\in H$ put $f_v = T_h(v)$ and
$f_v'= T(v)$. By Proposition 4.12 one has $f_v(1) = f_v'(1)$ for any $v\in H$. Let
$g\in G$ and let $g=k\,a\,n$ be the Iwasawa decomposition of $g$. Then since
$f_v,\,f_v'\in H(\delta,\xi)$ one has
$f_v(g) = e^{-\xi(log\,a)}\, f_v(k)$ and $f_v'(g) = e^{-\xi(log\,a)}\, f_v'(k)$. But
since $T$ and $T_h$ are $\Bbb C[K]$-module maps one has $f_v'(k) = f'_{k^{-1}\cdot
v}(1)$ and $f_v(k) = f_{k^{-1}\cdot v}(1)$. Thus $f_v(g) = f_v'(g)$. Hence $T = T_h$.
QED \vskip 1pc

{\bf 4.8.} Any $\mu\in \Lambda_{M_e}$ defines a
linear functional on the Cartan subalgebra $\hh_{\m}^s$ of the semisimple Lie algebra
$[\m,\m]$  (see (146)). The set of elements $\{h_i\},\,i\in I_{\m}$, is a basis of $\hh_{\m}^s$
(see \S 2.1 and \S 2.7). For $i\in I_{\m}$ let $n_i(\mu)\in \Bbb Z_+$ be defined by
putting $n_i(\mu) = \langle \mu,h_i\rangle $. Since $\lambda|\hh_{\m}\in \Lambda_{M_e}$ for
any $\lambda\in \Lambda$ note that $n_i(\lambda) = n_i(\lambda|\hh_{\m})$ using the
notation of \S 1.1 so that no ambiguity should arise. 

Let $(\delta,\xi) \in D\times \a^*$ and let write $\delta =
(\zeta,\nu) \in \widehat {F}_{s}\times \Lambda_{M_e}$. One has $\nu$ and $\nu^c\in
\Lambda_{M_e}$. See (206). We now observe that the functions in $H(\delta,\xi)\s
C^{\infty}(G)$ satisfy certain differential equations on $G$. \vskip 1pc {\bf Proposition
4.14.} {\it Let $f\in H(\delta,\xi)$ using the notation above. Then 
$$\eqalign{&(a)\,\,f\cdot u=
\nu^c(u)\, f\,\,\hbox{and}\cr
&(b)\,\,f\cdot e_{-\alpha_i}^{n_i(\nu^c) +1} = 0,\,\,\forall i\in I_{\m}\cr}\eqno
(242)$$}\vskip 1pc {\bf Proof.} The statement (a) has already been established. See (a) in
(218). Now recalling (216) and (218) there exists $\phi\in H(Y_{\delta},\xi)$ such that $f =
Q_{\delta}(\phi)$ where $Q_{\delta}$ is defined by (209). Thus if $g\in G$ and $m\in M$ one
has $$\eqalign{f(g\,m) &= \langle \phi(g\,m),y^*_{\delta}\rangle\cr
&= \langle \tau_{\delta}(m^{-1})\phi(g),y^*_{\delta}\rangle\cr
&= \langle \phi(g),m\cdot y^*_{\delta}\rangle\cr}\eqno (243)$$ Let $x\in \m_o$.
Replacing $m$ in (243) by $exp\,t\,x$ and differentiating at $t=0$ yields $$(f\cdot x)(g) =
\langle \phi(g),x\cdot y^*_{\delta}\rangle\eqno (244)$$ But then by complexification (244)
is valid for any $x\in \m$. However, for $i\in I_{\m}$, one has $e_{-\alpha_i}^{n_i(\nu^c)
+1}\cdot y_{\delta}^*=0$ by Proposition 1.1 applied to the case where $[\m,\m]$ replaces
$\g$. But then iterating the argument which yields (244) readily establishes
(242). QED\vskip 1pc  For any $\lambda\in \Lambda$ let $\lambda^c$ be the highest weight of
the dual $\g$-module $V_{\lambda}^*$ and let $0\neq v_{\lambda}^*$ be a highest weight
vector of $V_{\lambda}^*$. Let $\kappa$ be the longest element in the Weyl group of
$(\g,\hh)$ so that $\kappa\,\lambda$ is the lowest weight in $V_{\lambda}$. The pairing of
$V_{\lambda}$ and $V_{\lambda}^*$ immediately implies $$\lambda^c = -\kappa\,\lambda\eqno
(245)$$ Now obviously $(V_{\lambda})_{\n} = V_{\lambda}/(\n\,V_{\lambda})$ and
$(V_{\lambda}^*)^{\n}$ (see (93)) have the structure of $M$-modules. \vskip 1pc {\bf
Proposition 4.15.} {\it Let
$\lambda\in \Lambda$. Then $(V_{\lambda})_{\n}$ is an irreducible $M$-module. Moreover the
pairing of $V_{\lambda}$ and $V_{\lambda}^*$ induces a nonsingular $M$-invariant pairing
(see Proposition 3.1) of
$(V_{\lambda})_{\n}$ and $U(\m)\v^*_{\lambda}$.}\vskip 1pc {\bf Proof.} Obviously the
pairing of $V_{\lambda}$ and $V_{\lambda}^*$ induces a nonsingular
$M$-invariant pairing of
$(V_{\lambda})_{\n}$ and $(V_{\lambda}*)^{\n}$. But $(V_{\lambda}*)^{\n} =
U(\m)\,v_{\lambda}^*$ by (93) and $U(\m)\,v_{\lambda}^*$ is an irreducible $M$-module by
Proposition 3.1. Hence $(V_{\lambda})_{\n}$ must be an $M$-irreducible module. QED\vskip
1pc {\bf 4.9.} Let $\lambda\in \Lambda$. Then by Proposition 4.15 there exists $\delta\in D$
such that
$$(V_{\lambda})_{\n}\cong Y_{\delta}\eqno (246)$$ Of course $\delta = (\zeta,\nu)\in
\widehat{F}_{s}\times \Lambda_{M_e}$. It is easier to determine $\nu^c$  than
to determine $\nu$ itself. \vskip 1pc {\bf Proposition 4.16.} {\it Let $\lambda \in
\Lambda$ and let $\delta = (\zeta,\nu) \in \widehat{F}_{s}\times \Lambda_{M_e}$ be
defined by (246). Then in the notation of Proposition 3.1 one has $\zeta =
\lambda^c|F_{s}$ and $$\nu^c = \lambda^c|\hh_{\m}\eqno (247)$$ See (245).} \vskip 1pc
{\bf Proof.} Let $\delta^*\in \widehat {M}$ be defined so that $Y_{\delta}^*\cong
Y_{\delta^*}$. Obviously $\delta^* =(\zeta,\nu^c)$ as an element in
$\widehat{F}_{s}\times \Lambda_{M_e}$ since $\zeta$ is self-contragredient. But 
$Y_{\delta^*}\cong U(\m)\,v_{\lambda}^*$ by Proposition 3.1. The result then follows from
Proposition 3.1 since $\lambda^c$, given by (245), is the highest weight of
$V_{\lambda}^*$. QED\vskip 1pc Let $\lambda\in \Lambda$. Since $\a$ normalizes $\n$ it
follows that $(V_{\lambda})_{\n}$ has the structure of an $\a$-module. Furthermore 
since $\a$ commutes with $\m$, it follows that any element in $\a$ operates on
$(V_{\lambda})_{\n}$ by a scalar operator. Let $\xi\in \a^*$ be defined so that for
any $x\in \a$, $$x\,\,\hbox{operates as $\xi(x)$ on $(V_{\lambda})_{\n}$}\eqno
(248)$$ \vskip .5pc {\bf Proposition 4.17.} {\it Let
$\lambda\in \Lambda$. Let
$\xi\in \a^*$ be defined by (248). Then  $\xi = \kappa\,\lambda|\a$ where $\kappa$ is
the long element of the Weyl group of $(\g,\hh)$.} \vskip 1pc {\bf Proof.} The lowest
weight of
$V_{\lambda}$ is $\kappa\,\lambda$. Since any nonzero lowest weight vector clearly has a
nonzero projection in $(V_{\lambda})_{\n}$, the scalar defined by any $x\in \a$ on
$(V_{\lambda})_{\n}$ must be $\xi(x)$. QED\vskip 1pc We recover the following result of
Wallach. See \S 8.5 in [W]. \vskip 1pc {\bf Proposition 4.18.} {\it Let $\lambda\in \Lambda$
and let 
$\delta\in D$ and $\xi\in \a*$ be defined respectively by (246) and (248). Then regarding
$V_{\lambda}$ as the Harish-Chandra module $H$ of Proposition 4.10, Theorems 4.11 and
4.13 the pair
$(\delta,\xi)\in D\times
\a^*$ is the unique pair such that $H^*_{\delta,\xi} \neq 0$ (see (223)). Moreover
$$H^*_{\delta,\xi} = \Bbb C v_{\lambda}^*\eqno (249)$$ and hence, up to scalar multiple, 
$h =
v_{\lambda}^*$ is the unique element in $ H^*$ which satisfies Proposition 4.10. In
particular (recalling Theorems 4.11 and 4.13) $H(\delta,\xi)$ is the unique principal
series Harish-Chandra module which has a finite dimensional $(\g,K)$-submodule equivalent
to $V_{\lambda}$. Furthermore, the equivalence is defined by $T_h$ (see (231)) which, in this
case, is explicitly determined by
$$T_h(v)(g) = \langle v,g\,v^*_{\lambda}\rangle \eqno (250)$$ for any $v\in V_{\lambda}$ and
$g\in G$.} \vskip 1pc {\bf Proof.} If $(\delta',\xi')\in D\times \a^*$ and there exists
$0\neq h\in H^*_{\delta',\xi'}$ then $h$ is a highest weight vector in $H^*$ by (223). But
$H^* = V_{\lambda}^*$ is irreducible and hence any highest weight vector is a multiple of
$v^*_{\lambda}$. To prove (249) it therefore suffices, by (223), to show
that $$b\cdot v_{\lambda}^* = \chi_{\delta,\xi}(b)\,v_{\lambda}^*\eqno (251)$$ for any
$b\in {\cal B}$. As observed in (222) it suffices to take $b = \varepsilon,\,z$ and $x$
using the notation of (222). But $\delta = (\zeta,\nu)$ using the notation of
Proposition 4.16. But then (251) is satisfied for $b = \varepsilon$ or $b= z$ by
Propositions 3.1 and 4.16 since $v_{\lambda}^*$ is the highest weight vector of
$V_{\lambda}^*$ and the highest weight is $\lambda^c$. But now  $x\cdot v_{\lambda}^* =
\lambda^c(x)v_{\lambda}^*$ for $x\in \a$. However $\lambda^c = -\kappa\,\lambda$ by (245) and
hence
$$-\xi = \lambda^c|\a\eqno (252)$$ by Proposition 4.17. Recalling (222) this establishes
(249). But now for $h= v_{\lambda}^*,\,g\in G$ and $v\in V_{\lambda}$ one has $(T_h(v))(g) =
e^{-\xi(log\,a)}\,\langle k^{-1}\cdot v,v_{\lambda}^*\rangle$ by Proposition 4.10 where
$g=k\,a\,n$ is the Iwasawa decomposition of $g$. But $\langle k^{-1}\cdot
v,v_{\lambda}^*\rangle = \langle v,k\cdot v_{\lambda}^*\rangle$ and $a\,n\cdot
v_{\lambda}^* = e^{-\xi(log\,a)}\,v_{\lambda}^*$ by (252). This establishes (250). The
map $T_h$ is clearly nontrivial and hence is injective since $V_{\lambda}$ is
$\g$-irreducible. QED \vskip 1pc {\bf 4.10.} Let $\lambda\in \Lambda$ and let $(\delta,\xi)\in
D\times
\a^*$ be defined by (246) and (248). Let $h$ be as in Proposition 4.18. Let 
${\cal V}_{\lambda} = T_h(V_{\lambda})$ so that, by Proposition 4.18, ${\cal V}_{\lambda}$
is a finite dimensional $(\g,K)$-submodule, equivalent to $V_{\lambda}$, of the principal
series Harish-Chandra module $H(\delta,\xi)$. The functions in $H(\delta,\xi)$ are of
course determined by their restrictions to $K$. The following theorem characterizes the
functions in the finite dimensional space 
${\cal V}_{\lambda}$ as solutions in
$H(\delta,\xi)$ of differential equations which are
associated with certain elements in $U(\k)$. This can viewed as analogous to the Borel-Weil
theorem. In that case the differential equations are Cauchy-Riemann equations. 

Let $Diff\,\,\,C^{\infty}(G)$ be the algebra of all differential operators on
$C^{\infty}(G)$. One has an injective homomorphism $$U(\g)\to Diff\,C^{\infty}(G)\eqno (253)$$
where the image is the algebra of all left $G$-invariant differential operators on $G$. If
$f\in C^{\infty}(G)$ and $u\in U(\g)$ the action of $u$ on $f$ is denoted by $f\cdot u$. The
homomorphism (253) is determined by the restriction of (253) to $\g_o$. If $z\in \g_o$ and
$f\in C^{\infty}(G)$, then $f\cdot z$ is given in (205) where $G_1=G$. The following
noncompact analogue of Borel-Weil is one of the main results in the paper. \vskip 1pc {\bf
Theorem 4.19.} {\it Let $\lambda \in \Lambda$ and let $\delta \in D$ and $\xi\in a^*$ be
defined by (246) and (248). Recall (see Proposition 4.18) that $H(\delta,\xi)$ is the unique
principal series Harish-Chandra module which has a $(\g,K)$-submodule (also unique) 
${\cal V}_{\lambda}$ that is $(\g,K)$-equivalent to $V_{\lambda}$. Then ${\cal
V}_{\lambda}$ may be given by $${\cal V}_{\lambda} = \{f\in H(\delta,\xi)\mid f\cdot
q_{\lambda^c,i}(z_i)= 0,\,\,\forall i\in I_{\n}\}\eqno (254)$$ where $z_i\in
\k$, for
$i\in I_{\n}$, is defined by (97) and $q_{\lambda^c,i}(t)\in 
\Bbb C[t]$ is the polynomial defined in \S 2.2} \vskip 1pc {\bf Proof.} Let $f\in {\cal
V}_{\lambda}$. Then by (250) there exists $v\in V_{\lambda}$ such that, for any $g\in G$,
$f(g) = \langle v,g\,v_{\lambda}^*\rangle$. Clearly then $$(f\cdot u)(g)= \langle
v,g\,u\,v_{\lambda}^*\rangle\eqno (255)$$ for any $u\in U(\g)$. Now let ${\cal V}$ be the
right hand side of (254). To align the notation of (97) with that of Theorem 1.24 note that
(see paragraph preceding Lemma 1.12)
$$z_{-\alpha_i} = z_i,\,\,\hbox{for}\,\,i\in I_{\n}\eqno (256)$$ and $$z_{-\alpha_i} =
e_{-\alpha_i},\,\,\hbox{for}\,\,i\in I_{\m}\eqno (257)$$ But since $v_{\lambda}^*$ is the
highest weight vector in $V_{\lambda}^*$ and the highest weight is $\lambda^c$, it follows
that $$q_{\lambda^c,i}(z_i) v_{\lambda}^*= 0,\,\,\forall i\in I_{\n}$$ by Theorem 1.24.
Thus ${\cal V}_{\lambda}\s {\cal V}$ by (255). Hence to
prove (254) it suffices to prove that $$dim\,{\cal V}\leq d_{\lambda}\eqno (258)$$ where
$d_{\lambda} = dim\,V_{\lambda}$. Now let $\chi:U(\a +\n)\to \Bbb C$ be the character
defined so that $w\,v_{\lambda}^* = \chi(w)\,v_{\lambda}^*$ for all $w\in U(\a +\n)$. Since
$U(\g)$ is isomorphic, as a linear space, to $U(\k)\otimes U(\a +\n)$, it follows that
$U(\g) = U(\k) \oplus I$ where $I = U(\k)\,Ker\,\chi$. On the other hand, if
$L_{\lambda^c}(\k)$ is the left ideal in $U(\k)$ defined as in Proposition 1.3, it follows
from Proposition 1.3 that there exists a subspace $W\s U(\k)$ of dimension $d_{\lambda}$
such that $U(\k) = W\oplus L_{\lambda^c}(\k)$. Hence one has the direct sum $$U(\g) = W
\oplus L_{\lambda^c}(\k) \oplus I\eqno (259)$$ But now, by (218), any $f\in H(\delta,\xi)$
satisfies the $A\,N-\xi$ condition (see (215)). Thus in particular one has 
$$f\cdot I = 0,\,\,\forall f\in {\cal V}\eqno (260)$$ On the other hand we assert that $$f\cdot
L_{\lambda^c}(\k) = 0,\,\,\forall f\in {\cal V}\eqno (261)$$ To prove (261) first note that,
using of course Theorem 1.24, we can write $$L_{\lambda^c}(\k)= \sum_{i\in
I_{\n}}U(\k)\,q_{\lambda^c,i}(z_i) + \sum_{i\in
I_{\m}}U(\k)\,e_{-\alpha_i}^{n_i(\lambda^c)+1} + U(\k)\,Ker\,\nu^c\eqno (262)$$ The first
two sums in (262) replaces the first sum in (78) (recalling (255) and (256)) and the last
summand in (262) rewrites the last two sums in (78) using (247) and our extension of the
domain of $\nu^c$ to
$\b_{\m}$. But now (261) follows from Proposition 4.14 and the definition of ${\cal V}$. 

Now assume that $dim\,{\cal V}> d_{\lambda}$. Then there must exist $0\neq f\in {\cal V}$
such that $(f\cdot w)(1) = 0$ for all $w\in W$. But then by (259), (260) and (261) one has
$(f\cdot u)(1) = 0$ for all $u\in U(\g)$. This is a contradiction since $f$ is an analytic
function and hence is uniquely determined by the set of values $\{(f\cdot u)(1)\},\,\,u\in
U(\g)$. QED \vskip 1pc 

\centerline{\bf References}\vskip 1.3pc
\parindent=42pt

\item {[A]} S. Araki, On root systems and an infinitesmal classification of irreducible
symmetric spaces, {\it J. Math. Osaka City Univ.}, {\bf 13}(1962), 1-34
\item {[BB]} A. Beilenson and J. Bernstein, A generalization of Casselman's submodule 
theorem, {\it Representation Theory of Reductive Groups}, Birkhauser, Progress in
Mathematics, 1982, 35-67
\item {[C]} W. Casselman, Differential equations satisfied by matrix
coefficients,                                                                                                                                                                                                                                                                                                                                                                                                                                                                                                                                                                                                                                                                                                                                                                                                                                                                                                                                                                                                                                                                                                                                                                                                                                                                                                                                                                                                                                                                                                                                                                                                                                                                                                                                                                                                                                                                                                                                                                                                                                                                                                                                                                                                                                              
preprint
\item {[K1]} B. Kostant, On the existence and irreducibility of certain series
of representations, {\it Lie groups and their representations} edited by I.M.
Gelfand, Halsted Press, John Wiley, 1975, 231-329
\item {[K2]} B. Kostant, On Whittaker Vectors and Representation Theory, {\it Inventiones
Math.}, {\bf 48}(1978), 101-184.
\item {[K3]} B. Kostant, Clifford Algebra Analogue of the Hopf-Koszul-Samelson
Theorem, the $\rho$-Decomposition, $C(\g) = End\,V_{\rho}\otimes C(P)$, and the $\g$ - Module
Structure of $\wedge \g$, {\it Adv. in Math.}, {\bf 125}(1997), 275-350.
\item {[PRV]} K. R. Parthasarathy, R. Ranga Rao, and V. S. Varadarajan,
Representations of complex semisimple Lie groups and Lie algebras, Ann. of Math. (2)
Vol 85(1967), 383-429
\item {[V]} D. Vogan, Irreducible characters of semisimple Lie groups IV.
Character-multiplicity duality, {\it Duke Math. J.}, {\bf 49}(1982), 943-1073
\item {[W]} N. Wallach, {\it Harmonic analysis on homogeneous spaces}, Dekker, {\bf 19}(1973)

\smallskip
\parindent=30pt
\baselineskip=14pt
\vskip 1.9pc
\vbox to 60pt{\hbox{Bertram Kostant}
      \hbox{Dept. of Math.}
      \hbox{MIT}
      \hbox{Cambridge, MA 02139}}\vskip 1pc

      \noindent E-mail kostant@math.mit.edu

\end